%% file: thesis.tex
\documentclass[12pt,phd]{dms}

\usepackage[latin1]{inputenc}
\usepackage[T1]{fontenc}
\usepackage[dvips]{graphicx}
\usepackage{epsf, subfigure, verbatim, listings}
\usepackage{eucal}
\usepackage{psfrag}
\usepackage{epsfig}

\usepackage{amsmath, amssymb, amsthm}

\DeclareMathAlphabet{\mathscr}{OT1}{pzc}%
                                        {m}{it}

\newtheorem{theoreme}{Théorême}[section]
\newtheorem{conjec}{Conjecture}[section]

\theoremstyle{definition}
\newtheorem{defin}[theoreme]{Définition}

\numberwithin{equation}{section}
\numberwithin{table}{chapter} 
\numberwithin{figure}{chapter} 

\renewcommand{\baselinestretch}{1.5}

\begin{document}


\mdate{le $1^\text{er}$ mars 1997}      
\version{1} 


\title{GRAPHES EULÉRIENS ET COMPLÉMENTARITÉ LOCALE}      
\author{Fran\c{c}ois Genest}
\copyrightyear{2001}
\president{}
\directeur{M. Gert Sabidussi}
\membrejury{} 
\examinateur{M. Herbert Fleischner}   
\repdoyen{}
\dateacceptation{}
\sujet{Mathématiques}

\pagenumbering{roman}
\maketitle


\chapter*{Sommaire} 

\input{sommaire.tex}


\chapter*{Summary}

\input{summary.tex}


\chapter*{Dédicace}

À Marie, que j'admire.


\chapter*{Remerciements}

\input{remerciements.tex}


\tableofcontents
\listoffigures

%
%

\NoChapterPageNumber
\pagenumbering{arabic}


\chapter*{Introduction}

\input{introduction.tex}


\input{words.tex}

\input{transition.tex}


\chapter*{Conclusion}

\input{conclusion.tex}

\appendix


\chapter{Classe de parité des pentagones siamois}

\input{annexe1.tex}


\chapter{Classe de parité de $\operatorname{Cay}(\mathbb{Z}_{13},\pm \{1,3,4\})$}

\input{annexe2.tex}


\chapter{Programmes}

Cette annexe comprend les listings des différents logiciels qui ont
été écrits dans le cadre de cette thèse. Ces programmes ont été conçus
afin de générer des graphes et d'en tester la pureté.

Les structures de données standard utilisées sont tirées de
\cite{Aho}. La structure \flqq coloriage\frqq \ et ses méthodes sont basés sur l'algorithme de test d'isomorphisme que l'on retrouve dans \cite{Babai}.\\

\renewcommand{\baselinestretch}{1}
\footnotesize
\lstset{language=C++, breaklines, texcl, escapechar=ÿ}
\lstset{escapebegin=\normalsize\begin{itshape}}
\lstset{escapeend=\end{itshape}\footnotesize}
\lstset{tabsize=2}

{\vspace{8mm}\begin{center}{\large \bf makefile}\end{center}\vspace{3mm}}


{\vspace{8mm}\begin{center}{\large \bf circulante.cpp}\end{center}\vspace{3mm}}


{\vspace{8mm}\begin{center}{\large \bf inverse.cpp}\end{center}\vspace{3mm}}


{\vspace{8mm}\begin{center}{\large \bf isomorphes.cpp}\end{center}\vspace{3mm}}


{\vspace{8mm}\begin{center}{\large \bf pur.cpp}\end{center}\vspace{3mm}}


{\vspace{8mm}\begin{center}{\large \bf bg\_liste.hpp}\end{center}\vspace{3mm}}


{\vspace{8mm}\begin{center}{\large \bf classe\_de\_parite.hpp}\end{center}\vspace{3mm}}


{\vspace{8mm}\begin{center}{\large \bf graphe.hpp}\end{center}\vspace{3mm}}


{\vspace{8mm}\begin{center}{\large \bf keytype.hpp}\end{center}\vspace{3mm}}


{\vspace{8mm}\begin{center}{\large \bf bg\_liste.cpp}\end{center}\vspace{3mm}}


{\vspace{8mm}\begin{center}{\large \bf bigraphe.cpp}\end{center}\vspace{3mm}}


{\vspace{8mm}\begin{center}{\large \bf classe\_de\_parite.cpp}\end{center}\vspace{3mm}}


{\vspace{8mm}\begin{center}{\large \bf coloriage.cpp}\end{center}\vspace{3mm}}


{\vspace{8mm}\begin{center}{\large \bf graphe.cpp}\end{center}\vspace{3mm}}


{\vspace{8mm}\begin{center}{\large \bf keytype.cpp}\end{center}\vspace{3mm}}


{\vspace{8mm}\begin{center}{\large \bf vertex.cpp}\end{center}\vspace{3mm}}


{\vspace{8mm}\begin{center}{\large \bf vset.cpp}\end{center}\vspace{3mm}}

\renewcommand{\baselinestretch}{1.5}
\normalsize


\end{document}

%% file: sommaire.tex
Nous définissons dans cette thèse la notion de complémentation de
graphes bicoloriés et les notions de graphe pur et de graphe
inversible. L'étude des graphes purs est motivée par deux conjectures
concernant les systèmes de transitions de graphes eulériens et par la
conjecture de double recouvrement.

L'utilisation de règles de substitution nous permet de déterminer
quand deux suites de complémentation donnent le même graphe. Pour les
graphes bicoloriés, ces suites de complémentation font place à des
ensembles de complémentation.

Les graphes inversibles (les graphes bicoloriés dont l'ensemble des
sommets est un ensemble de complémentation) ont ceci de particulier
que leur inverse possède les mêmes automorphismes. L'inversibilité se
définit aussi pour les graphes non coloriés en les munissant de leur
coloriage naturel.

Il est proposé que la caractérisation des graphes purs permettrait
de valider la conjecture de double recouvrement. Nous
décrivons comment les graphes purs ont des factorisations essentielles
en graphes purs primitifs. Les quatre classes de parité connues de graphes
purs primitifs sont présentées. Les listings des programmes ayant
permis d'établir la pureté de ces graphes sont inclus.\\

\noindent {\bf mots clés} : complémentarité locale, systèmes de transitions, graphes purs, graphes inversibles, graphes eulériens, double recouvrement, orthogonalité.

%% file: summary.tex
We define pure graphs, invertible graphs, and the notion of complementation of bicoloured graphs. The study of pure graphs is motivated by two conjectures about the transition systems of eulerian graphs and by the Cycle Double Cover Conjecture.

We show how substitution rules can be used to determine when two complementation words produce the same graph. For bicoloured graphs, complementation words give way to complementation sets.

The invertible graphs (bicoloured graphs whose vertex set is a complementation set) are shown to have the property that their inverse has the same automorphisms. The property of being invertible can also be defined for non-coloured graphs by endowing them with their natural colouring.

It is proposed that a characterization of pure graphs would contribute to establish the truth of the Cycle Double Cover Conjecture. We show how pure graphs have essential factorizations into primitive pure graphs. The four primitive pure parity classes are presented. Included are the listings of the programs used to test graphs for purity.\\

\noindent{\bf keywords}: local complementation, transition systems, pure graphs, invertible graphs, eulerian graphs, double cycle cover, orthogonality.

%% file: remerciements.tex
Avant tout, je tiens à remercier mon directeur de recherche, Gert Sabidussi. À chacun de nos entretiens, j'ai apprécié sa rigueur de raisonnement et son caractère franc et jovial. Ces années de recherches ont comporté leur part de doutes et si j'ai gardé confiance en la qualité de mon travail, c'est en grande partie parce qu'il m'a traité en collègue.\\

Je remercie Marie, pour tout ce qu'elle est.\\

Comme d'autres avant moi, je remercie Jérome Fournier, qui sait mettre les gens à l'aise et contribue beaucoup à rendre l'atmosphère du département accueillante.\\

Merci à Alexandre Girouard, Louis-Sébastien Guimond, Carsten Heinz, Donald Knuth, Leslie Lamport, Linus Torvalds, et tous ceux qui ont facilité mon travail de rédaction.\\

Merci au FCAR et à l'Université de Montréal pour leurs contributions financières.

%% file: introduction.tex
Voici un graphe simple (sans boucles ni arêtes
multiples) dont chaque sommet est colorié en blanc ou
en noir :

\psfrag{u1}[][][1]{$u_1$}
\psfrag{u2}[][][1]{$u_2$}
\psfrag{u3}[][][1]{$u_3$}
\psfrag{u4}[][][1]{$u_4$}
\psfrag{u5}[][][1]{$u_5$}
\psfrag{u6}[][][1]{$u_6$}
\psfrag{u7}[][][1]{$u_7$}
\psfrag{u}[][][1]{$u$}
\psfrag{v}[][][1]{$v$}
\psfrag{w}[][][1]{$w$}
\psfrag{G}[][][1]{$G$}
\psfrag{H}[][][1]{$H$}

\[ \epsfig{file=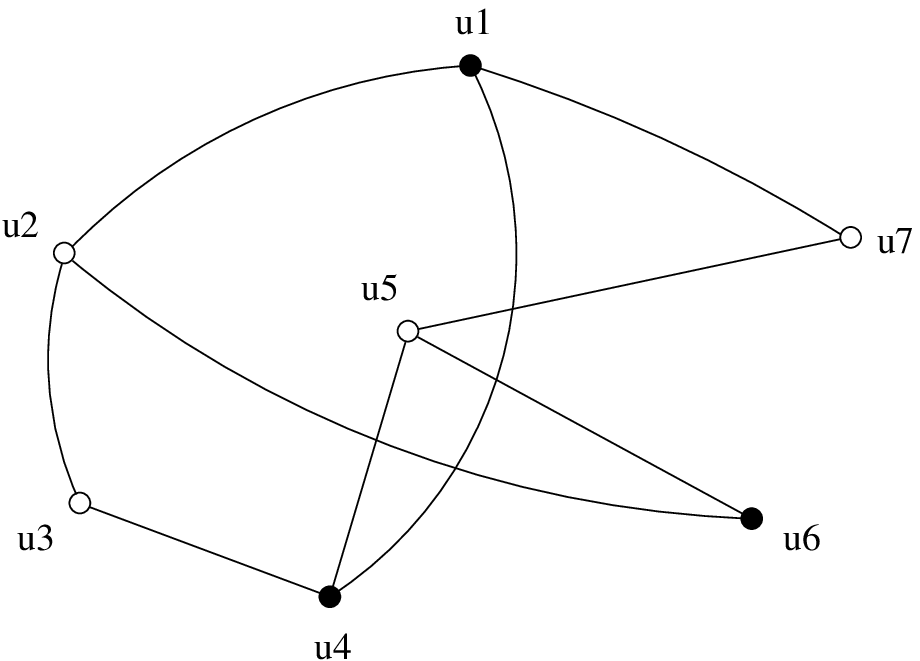,height=5cm}
\]

\noindent Peut-on trouver un sous-ensemble $A$ des sommets du graphe qui
soit un \emph{indépendant} (aucune paire de sommets ne sont adjacents) maximal
(tout autre sommet a un voisin dans $A$) et qui ne contienne que des
sommets noirs? Un tel sous-ensemble est appelé une \emph{anticlique
noire}.

Il est assez facile de se convaincre que, dans l'exemple présenté, le
graphe n'admet aucune anticlique noire. En effet, les sommets blancs
$u_3$ et $u_7$ n'ont chacun qu'un seul sommet noir comme voisin. Toute
anticlique noire devrait donc contenir à la fois les sommets $u_1$ et
$u_4$, mais ceux-ci sont adjacents.

Définissons maintenant un jeu dont le but est encore de trouver une
anticlique noire, mais dans lequel il est permis de modifier le graphe
de départ en
effectuant certaines \emph{complémentations locales}. Complémenter
localement un graphe par rapport à un sommet donné, c'est inverser les
adjacences entre ses voisins. (c.-à-d. si $v$ et $w$ sont adjacents à $u$
et qu'il y a une arête entre $v$ et $w$, la complémentation locale par
rapport à $u$ fait disparaître cette arête; s'il n'y en a pas, elle en
fait apparaître une). Par exemple, les graphes $G$ et $H$ suivants s'obtiennent l'un de l'autre en complémentant par rapport à $u$:

\begin{center}
\epsfig{file=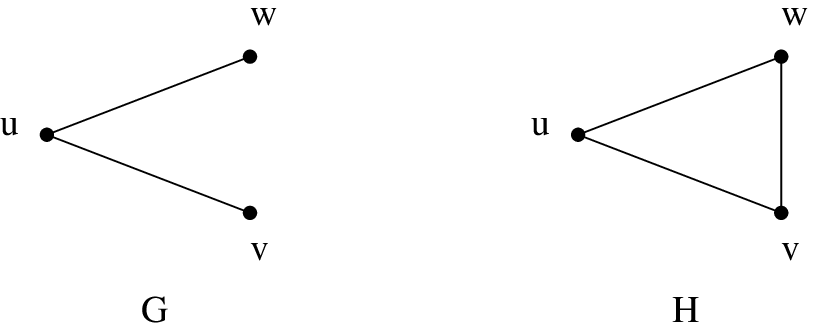,height=4cm}
\end{center}

\noindent Dans ce jeu, il est permis de jouer à un sommet blanc,
ce qui effectue une complémentation locale par rapport à ce sommet et
qui, en plus, inverse la couleur de chacun des voisins. Il est
également permis de jouer à une arête incidente avec deux sommets noirs,
disons $u$ et $v$, ce qui produit le même effet que trois
complémentations locales successives : d'abord par rapport à $u$, puis à $v$
et de nouveau à $u$ (ceci est bien défini car inverser les rôles de
$u$ et $v$ produit le même résultat). En jouant au sommet
$u_2$ dans l'exemple donné, nous obtenons :

\[ \epsfig{file=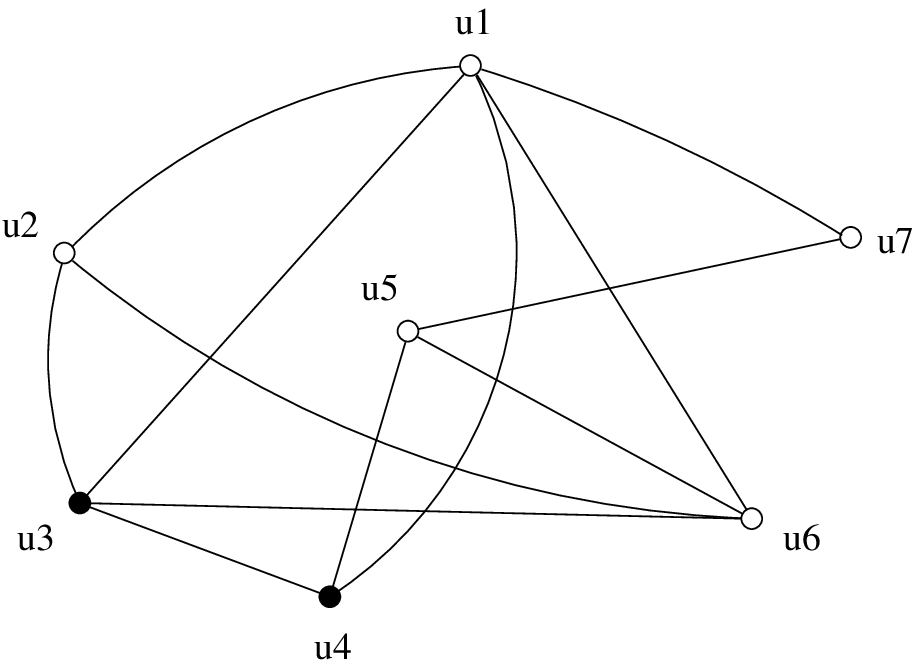,height=5cm}
\]

\noindent Il ne restera plus qu'à jouer au sommet $u_7$ pour créer une
anticlique noire $A=\{u_3,u_5\}$ dans le graphe résultant.

En permettant ces opérations, l'expérience démontre que le jeu a une
solution pour la grande majorité des graphes bicoloriés blanc et
noir. C'est-à-dire qu'après un certain nombre de coups nous avons
toutes les chances de faire surgir une anticlique noire. Comme souvent
en mathématiques, ce sont les exceptions qui vont nous intéresser :
les graphes bicoloriés pour lesquels le jeu ne fait apparaître aucune
anticlique noire sont appelés les graphes \emph{purs}.

Il est difficile de croire que ce jeu étrange soit en relation avec les graphes
eulériens. C'est pourtant le cas et, lorsqu'on se restreint à un
certain type de graphes (les graphes de cordes\footnote{Les graphes de
  cordes, aussi appelés \emph{graphes d'alternance}, ne doivent pas être
  confondus aves les graphes \emph{à cordes} (en englais, \emph{chordal graphs}).}, \emph{circle graphs} en
anglais), ce jeu est équivalent à un problème connu. Nous verrons
comment l'étude des graphes purs permet de mieux cerner ce problème et
jette un éclairage nouveau sur la conjecture de double recouvrement.

Pour mieux manipuler et décrire les graphes qui apparaissent au cours du jeu, il est naturel de se fixer un graphe de départ et d'identifier les autres graphes à l'aide des suites de coups permettant de les obtenir. Dans le premier article, nous chercherons à formaliser ces suites de coups, ce seront les \emph{suites de parité}, et nous verrons comment trouver des suites les plus simples possibles. Ce faisant, nous serons récompensés de notre diligence par un résultat étonnant permettant d'obtenir, à partir de certains graphes que nous appellerons \emph{inversibles}, d'autres graphes ayant les mêmes automorphismes. Dans le second article, nous établirons le lien entre le problème auquel nous faisions allusion plus haut, qui concerne les systèmes de transitions de graphes eulériens, et les graphes purs. Aussi, nous verrons comment les graphes purs ont des \emph{factorisations essentielles} en graphes purs \emph{primitifs}. Enfin, les quatre graphes purs primitifs connus (à transformation par le jeu près) seront présentés.

\section{Définitions}

Les définitions les plus importantes sont incluses. Là où elles
diffèrent de celles qu'on peut trouver dans la littérature, c'est
par souci de concision et de clarté dans le cadre du sujet exposé dans cette thèse. Pour les définitions de base qui
n'apparaissent pas dans cette section, comme les notions d'isomorphisme de graphes et de connexité, le lecteur peut se référer au livre de Bondy et Murty \cite{Bondy}.

\begin{defin} Un \emph{graphe} $G=(V,E)$ est constitué de deux
  ensembles disjoints $V$ et $E$ dont les éléments sont appelés
  respectivement les \emph{sommets} et les \emph{arêtes} de $G$ et
  d'une fonction d'\emph{incidence} qui associe à chaque arête $e$
  soit un sommet, auquel cas $e$ est appelée une \emph{boucle}, ou
  soit une paire non ordonnée de sommets distincts. On dit de ces
  sommets qu'il sont \emph{incidents} avec l'arête $e$, ou encore que
  ce sont les \emph{incidences} de $e$.
\end{defin}

L'\emph{ordre} de $G$ est le nombre de sommets du graphe. Si une arête
$e$ a les incidences $u$ et $v$, on dit que $e$ \emph{relie} ces
sommets et que $u$ et $v$ sont $adjacents$. Deux arêtes distinctes ayant une incidence en commun sont également dites $adjacentes$.

\begin{defin}
Des arêtes distinctes sont \emph{multiples} si elles ont les mêmes incidences. Un graphe sans boucles ni arêtes multiples est un graphe \emph{simple}.
\end{defin}

Dans le cas d'un graphe simple, nous adoptons la convention d'identifier chaque arête avec ses incidences, que nous notons entre crochets (c.-à-d. $[u,v]$ pour l'arête reliant $u$ et $v$).

\begin{defin}
Un \emph{chemin} dans un graphe est une suite $u_0e_1u_1e_2u_2...e_nu_n$
telle que $e_i$ est une arête ayant les incidences $u_{i-1}$ et $u_i$, $i=1,...,n$. Le chemin est \emph{fermé} si $u_0=u_n$. 
\end{defin}

\begin{defin}
Une \emph{chaîne} dans un graphe est une marche dont les sommets sont distincts.
\end{defin}

\begin{defin}
Un \emph{parcours} dans un graphe est une suite $u_0e_1u_1e_2u_2...e_n$ telle que $u_0e_1u_1e_2u_2...e_nu_0$ est un chemin fermé dont les arêtes sont distinctes.
\end{defin}

En général, le début ou la direction du parcours importent peu. Aussi, on dira que $e_n$ et $e_1$
sont des arêtes successives. Il existe une formalisation des parcours
permettant d'éviter de choisir un sommet de départ, ce sont les
\emph{permutations eulériennes} (Sabidussi \cite{Sabidussi}). Cependant, dans cette thèse, nous nous limiterons aux parcours tels que définis précédemment.

\begin{defin}
Un parcours est \emph{eulérien} s'il contient toutes les arêtes du graphe. Un graphe est \emph{eulérien} s'il admet un parcours eulérien.
\end{defin}

\begin{defin}
Un graphe $G_1=(V_1,E_1)$ est un \emph{sous-graphe} du graphe $G_2=(V_2,E_2)$ si $V_1\subset V_2$, $E_1\subset E_2$ et si chaque arête $e\in E_1$ a les mêmes incidences dans $G_1$ et $G_2$. C'est un sous-graphe \emph{induit} si toute arête de $G_2$ ayant ses incidences dans $V_1$ est dans $E_1$.
\end{defin}

\begin{defin}
Le \emph{degré} d'un sommet $u$ est deux fois le nombre de boucles
incidentes avec $u$ plus le nombre des autres arêtes incidentes avec
$u$. Un graphe est \emph{d-régulier} (ou simplement \emph{régulier})
si tous ses sommets ont le même degré $d$. Un sommet de degré nul est
dit \emph{isolé}.
\end{defin}

\begin{defin}
Le sous-graphe de $G=(V,E)$ \emph{induit} par $V'\subset V$ est l'unique sous-graphe induit de $G$ dont l'ensemble de sommets est $V'$. Le sous-graphe de $G$ \emph{induit} par $E'\subset E$ est le sous-graphe $G'=(V',E')$ où $V'\subset V$ est l'ensemble des incidences des arêtes dans $E'$.
\end{defin}

\begin{defin}
Un \emph{cycle} est un graphe non vide, 2-régulier et connexe. Un
\emph{cycle} d'un graphe $G$ est un cycle qui est sous-graphe de
$G$. Un \emph{m-cycle} est un cycle d'ordre $m$ (par exemple, un
1-\emph{cycle} est induit par une boucle). Une arête n'appartenant à aucun cycle de $G$ est un \emph{isthme} de $G$.
\end{defin}

\begin{defin}
Une \emph{décomposition en cycles (ou en parcours)} d'un graphe $G$ sans sommet isolé est une famille de cycles (parcours) de $G$ telle que chaque arête de $G$ appartient à exactement un de ces cycles (parcours).
\end{defin}

\section{Caractérisation des graphes eulériens}

Le théorême élémentaire suivant, dont l'objet est la caractérisation des graphes eulériens, réunit des résultats apparus au cours d'une période de près de deux cents ans, depuis un article d'Euler déposé en 1735 jusqu'à un résultat de Veblen en 1922. Il est utile de s'y attarder un moment car il restera en filigrane dans le reste de cette thèse.

\begin{theoreme}
Soit $G$ un graphe sans sommet isolé et avec un nombre fini de sommets
et d'arêtes. Les énoncés suivants sont équivalents: \\
(1) $G$ est eulérien; \\
(2) tous les degrés de $G$ sont pairs et $G$ est connexe; \\
(3) $G$ a une décomposition en cycles et il est connexe.
\end{theoreme}

\begin{proof}
(1)$\Rightarrow$(2) Soit $w$ un parcours eulérien de $G$. Comme $w$ relie chaque paire de sommet de $G$, celui-ci est connexe. Chaque passage du parcours à un sommet $u$ contribue 2 au degré de $u$ (une fois à l'entrée et une fois à la sortie). Donc tous les degrés de $G$ sont pairs.\\

(2)$\Rightarrow$(3) Démontrons l'affirmation plus générale que tout
   graphe sans sommet isolé dont les degrés sont pairs possède une
   décomposition en cycles. C'est vrai pour le graphe vide. Soit une
   chaîne $w=u_0e_1u_1...e_nu_n$ non vide de longueur maximale. Le
   sommet $u_n$ étant de degré pair, il est incident avec une arête
   $e_{n+1}\neq e_{n}$. Si $e_{n+1}$ est une boucle, elle induit un
   cycle. Sinon, puisque la chaîne est maximale, l'autre incidence de
   $e_{n+1}$ doit être $u_i$ où $i\in \{0,...,n-1\}$ et dans ce cas,
   le chemin fermé $u_{i}e_{i+1}u_{i+1}...e_{n+1}u_i$ coïncide avec un
   cycle du graphe. En retirant du graphe les arêtes d'un cycle, puis
   les sommets isolés restants, nous obtenons un graphe dont les
   degrés sont pairs. Le résultat en découle, par induction sur le nombre d'arêtes. \\

(3)$\Rightarrow$(1) Si le graphe est vide, le résultat est trivial. Nous allons construire une suite de parcours contenant successivement plus d'arêtes. Soit $C_1$ un cycle de la décomposition. Choisissons un parcours $w$ de $C_1$. Si ce parcours ne contient pas toutes les arêtes, puisque $G$ est connexe, il existe un cycle $C_2$ de la décomposition ayant un sommet en commun avec $w$. Lors d'un passage du parcours à ce sommet, on peut interrompre $w$ et parcourir $C_2$ avant de continuer. De cette fa\c{c}on, on obtient un parcours utilisant les arêtes de $C_1$ et de $C_2$. Ce procédé permet de modifier $w$ pour parcourir successivement plus de cycles, jusqu'à ce que tous les cycles aient été parcourus, auquel cas $w$ est eulérien.
\end{proof}

\section{Historique du problème}

C'est avec l'article d'Euler sur le problème des ponts de Königsberg que certains associent l'avènement de la théorie des graphes. Il y est question de la cité médiévale aujourd'hui appelée Kaliningrad, située sur la rivière Pregel en Prusse. Euler y demande s'il est possible de trouver un chemin qui traverse chacun des sept ponts une et une seule fois.

\[ \epsfig{file=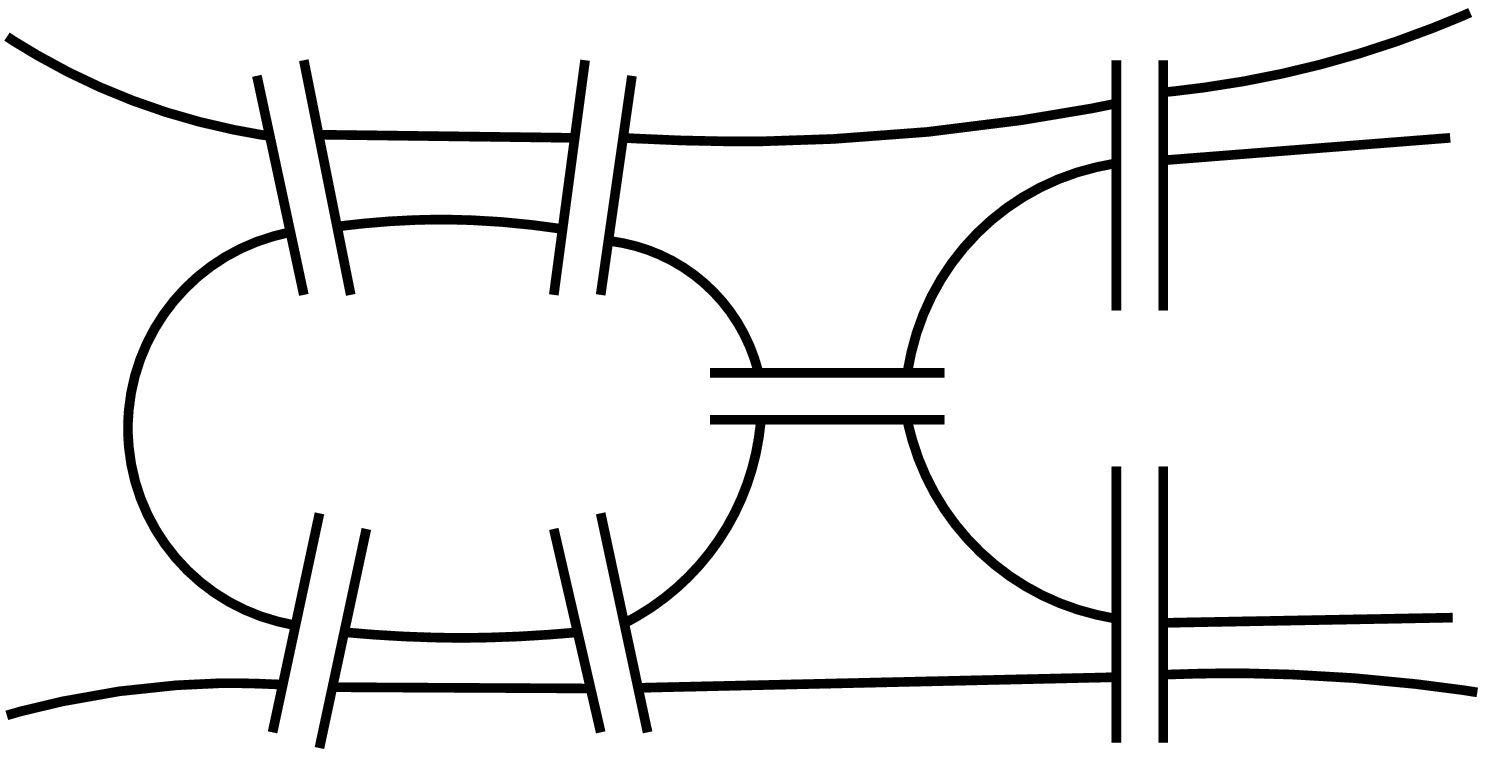,height=3cm}
\]
\begin{center}
Y a-t-il un chemin traversant les sept ponts de Königsberg?
\end{center}
\vspace{4mm}
\noindent Pour en donner la réponse, Euler fait une démonstration qui, dans le
cas où on demande de revenir au point de départ, se résume à la partie
(1)$\Rightarrow$(2) du théorême élémentaire. Pour lui, la réciproque a
peu d'intérêt et il faudra attendre en 1871 pour que Hierholzer en
fasse la preuve. Une historique intéressante du problème des ponts de
Königsberg se retrouve dans Wilson \cite{Wilson}. L'équivalence
(2)$\Leftrightarrow$\hspace{0in}(3) fut d'abord esquissée en 1912 puis établie en
1922 par Veblen \cite{Veblen1,Veblen2}.

En 1966, Kotzig \cite{Kotzig1} a introduit la notion d'\emph{orthogonalité} dans
les graphes eulériens. Étant donné une
décomposition en parcours d'un graphe connexe $4$-régulier, il s'est demandé s'il était toujours possible de
trouver un parcours eulérien orthogonal à la décomposition, en ce sens
qu'aucune paire d'arêtes ne se succèdent à la fois dans le parcours
eulérien et dans un des parcours de la décomposition (nous appelerons
des arêtes successives et le sommet entre les deux une $transition$ du parcours).

Se limitant aux graphes sans boucles, Kotzig a démontré que oui. Pour
le voir, considérons un graphe eulérien dont les degrés sont tous
$\geq 4$. Choisissons comme point de départ un parcours eulérien $w$ du
graphe et modifions ce parcours par étapes, de fa\c{c}on à réduire le nombre de
transitions fautives jusqu'à zéro. Soit
$\{u,e_1,e_2\}$ une transition fautive de $w$. Puisque $u$ est de degré $\geq$ 4, il existe au moins
une autre paire d'arêtes incidentes avec $u$, disons $\{e_3,e_4\}$,
qui se succèdent dans le parcours eulérien. À partir de ces deux
transitions, on peut modifier le
parcours de fa\c{c}on canonique : si possible, on fait se succéder $e_1$ et $e_3$ et aussi $e_2$ et $e_4$, tout en laissant le reste du parcours inchangé; sinon on fait se succéder $e_1$ et $e_4$ et également $e_2$ et
$e_3$. Il est facile de voir que dans exactement un des deux cas nous obtenons un nouveau parcours
eulérien. De plus, le nombre de transitions fautives dans ce parcours
aura diminué de 1 ou 2.

\psfrag{e1}[][][1]{$e_1$}
\psfrag{e2}[][][1]{$e_2$}
\psfrag{e3}[][][1]{$e_3$}
\psfrag{e4}[][][1]{$e_4$}
\psfrag{u}[][][1]{$u$}
\psfrag{transition fautive}[][][1]{transition fautive}
\vspace{4mm}
\begin{center}
\epsfig{file=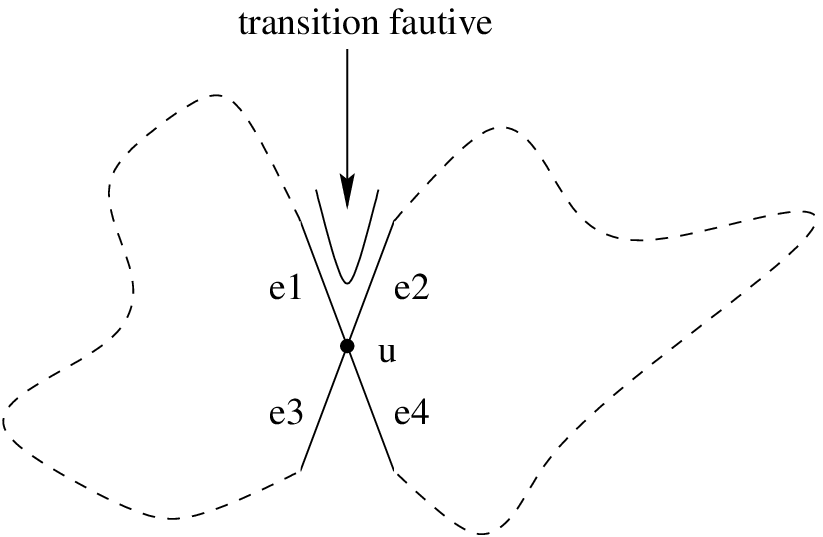,height=4cm}
\end{center}

\begin{defin}
Une \emph{transition} à un sommet $u$ est soit un couple $\{u,e\}$ où $e$ est une boucle incidente avec $u$ ou soit un ensemble $\{u,e_1,e_2\}$ où
$e_1$ et $e_2$ sont des arêtes distinctes incidentes avec $u$. Une \emph{transition}
d'un parcours $w$ (ou d'un cycle $C$) à un sommet $u$ est une transition dont les arêtes
se succèdent dans $w$ (sont adjacentes dans $C$). Un \emph{système de transitions} d'un graphe
$G$ est l'ensemble des transitions d'une décomposition en parcours ou
en cycles de
$G$.
\end{defin}

\begin{defin}
Deux systèmes de transitions sont \emph{orthogonaux} s'ils n'ont aucune
transition en commun. Un parcours eulérien ou une décomposition en
cycles sont orthogonaux à un système de transitions donné si les sytèmes
de transitions qu'ils induisent le sont.
\end{defin}

Le résultat de Kotzig peut donc s'exprimer ainsi : un graphe eulérien $G$ de
degré minimal $>2$ muni d'un système de transitions $S$ admet un
parcours eulérien orthogonal à $S$ si et seulement si pour tout sommet $u$ de degré 4 incident à une boucle $e$, la transition $\{u,e\}$ est dans $S$.

Étant donné les différentes caractérisations des graphes eulériens
données par le théorême élémentaire, il est naturel de poser la
question analogue : quand un graphe eulérien de degré
minimal $>2$ muni d'un système de transitions $S$ admet-il une
décomposition en cycles orthogonale à $S$?

Une condition nécessaire est que $S$ ne contienne aucune transition induite par une boucle (les transitions de type $\{u,e\}$). De plus, tout graphe obtenu en retirant les arêtes $e_1$ et $e_2$ d'une transition $\{u,e_1,e_2\}\in S$ doit être connexe. Si ces conditions sont
remplies, $S$ est dit $admissible$.

Cependant, il existe des cas où l'on ne peut trouver de décomposition
en cycles orthogonale à un système de transitions admissible. Par exemple, le
graphe complet sur cinq sommets ($K_5$) muni d'un système de transitions
correspondant à une décomposition en deux 5-cycles :

\vspace{2mm}
\begin{center}
\epsfig{file=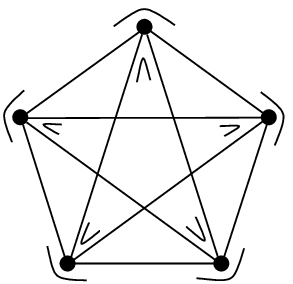,height=3.5cm}
\end{center}

\begin{center}
Un système de transitions sans décomposition en cycles orthogonale.
\end{center}
\vspace{3mm}

En 1975, Sabidussi (voir \cite{Fleischner1}) a émis la conjecture qu'une décomposition en cycles
orthogonale existe dans le
cas où $S$ correspond à un parcours eulérien (conjecture d'orthogonalité). Incidemment, ce sont les
travaux qu'il a effectués sur cette conjecture qui ont inspiré les
résultats que l'on retrouve dans cette thèse concernant la question
plus générale. Une condition suffisante dans le cas général (Fan et
Zhang \cite{Fan-Zhang}) est que $S$ soit admissible et que $G$ ne contienne pas de
sous-graphe isomorphe à une subdivision de $K_5$ (une $subdivision$
d'un graphe est obtenue en ajoutant des sommets qui \flqq
subdivisent\frqq \ les arêtes du graphe).
\vspace{3mm}
\begin{center}
\epsfig{file=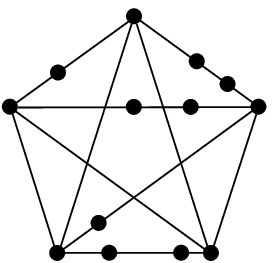,height=3.5cm}
\end{center}
\begin{center}
Une subdivision de $K_5$
\end{center}
\vspace{3mm}

\noindent En particulier, si $G$ est planaire et $S$ est admissible, on peut
trouver une décomposition en cycles orthogonale à $S$.

En ce qui concerne la conjecture d'orthogonalité, Sabidussi \cite{Sabidussi} nous dit
qu'il suffit de la démontrer pour les graphes bipartis de degrés 4 et
6 (ayant donc un nombre pair de sommets de degré 6). La conjecture est
vraie pour les graphes dont les degrés sont divisibles par 4
(Fleischner \cite{Fleischner2}, ou voir Jackson \cite{Jackson2}) et pour les graphes ayant exactement un
sommet de degré 6 et les autres de degré 4 (Fleischner \cite{Fleischner3}).

Revenant au théorême de caractérisation des graphes eulériens, on voit qu'il est futile
de chercher une décomposition en cycles d'un graphe ayant des sommets
de degré impair. Dans ce cas, si on veut recouvrir le graphe par des
cycles, il faut permettre aux arêtes d'apparaître dans plus d'un cycle. C'est peut-être cette réflexion qui est à l'origine
de la conjecture de double recouvrement.

\begin{defin}
Un \emph{double recouvrement par des cycles} (ou simplement un
\emph{double recouvrement}) d'un graphe $G$ est une famille de cycles
de $G$ telle que chaque arête de $G$ appartient à exactement deux de
ces cycles.
\end{defin}

\begin{conjec}[Conjecture de double recouvrement]
Tout graphe sans isthme possède un double recouvrement.
\end{conjec}

La paternité de cette conjecture n'est pas bien établie. Elle est
cependant très importante étant donné son lien avec la théorie des
flots à valeurs entières et avec les plongements de graphes dans des surfaces
(voir Jaeger \cite{Jaeger} et Jackson \cite{Jackson2}). De plus, cette conjecture est
intimement liée au problème d'existence d'une décomposition en cycles
orthogonale à un système de transitions. Comme nous le verrons dans le
deuxième article, tout contre-exemple à la conjecture de double
recouvrement qui serait minimal par rapport au nombre d'arêtes révélerait deux nouveaux graphes purs primitifs (à complémentation
près). Fleischner \cite{Fleischner3} a également montré que la conjecture dite \flqq de cycle
dominant\frqq \ et la conjecture d'orthogonalité impliquent, ensemble, la
conjecture de double recouvrement.\\

Abordons maintenant le lien entre les systèmes de transitions et le
jeu du départ. Étant donné un graphe simple quelconque, colorions en
blanc les sommets de degré pair et en noir les sommets impairs. En
jouant à un sommet blanc (pair) ou à une arête incidente avec des
sommets noirs (impairs), les couleurs du graphe résultant restent en
accord avec les parités des sommets, ce qu'on appelle un coloriage
$naturel$. À cause de ce rapprochement entre couleur et parité pour
une partie des graphes bicoloriés, on appele \emph{classe de parité}
de $G$ la famille des graphes bicoloriés pouvant s'obtenir par le jeu à partir
d'un graphe bicolorié $G$.

Sabidussi \cite{Sabidussi} a montré que pour chaque graphe eulérien $G$ de degrés 4
et 6 muni d'un parcours eulérien $w$, il correspond une unique classe
de parité de graphes naturellement coloriés et que $G$ possède une
décomposition en cycles orthogonale à $w$ si et seulement si la classe
de parité correspondante contient un graphe ayant une anticlique
noire. D'où l'intérêt de caractériser les graphes purs.


\section{Les graphes purs}

Il est assez facile de vérifier que le pentagone est pur :

\vspace{2mm}
\begin{center}
\epsfig{file=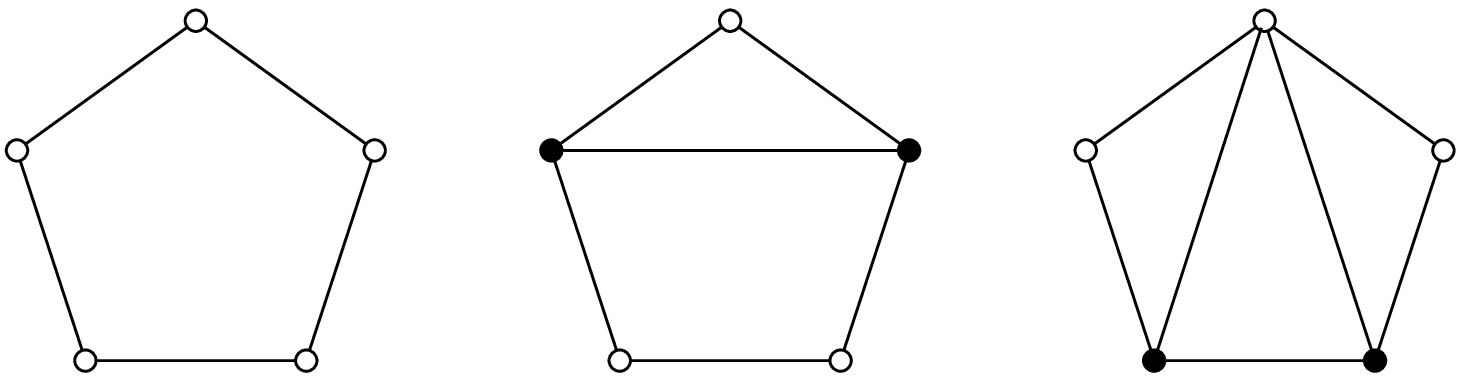, width=8.324cm}
\end{center}

\begin{center}
À isomorphisme près, les graphes de la classe de parité du pentagone.
\end{center}
\vspace{4mm}

\noindent Au moment de commencer mes travaux, les graphes de la classe
de parité du pentagone étaient les seuls graphes purs connus non
triviaux (un sommet isolé est pur). Sabidussi a alors émis l'hypothèse que, lorsque des pentagones sont identifiés en
un sommet, ils forment un graphe pur; comme exemple, voici le trèfle :

\begin{center}
\epsfig{file=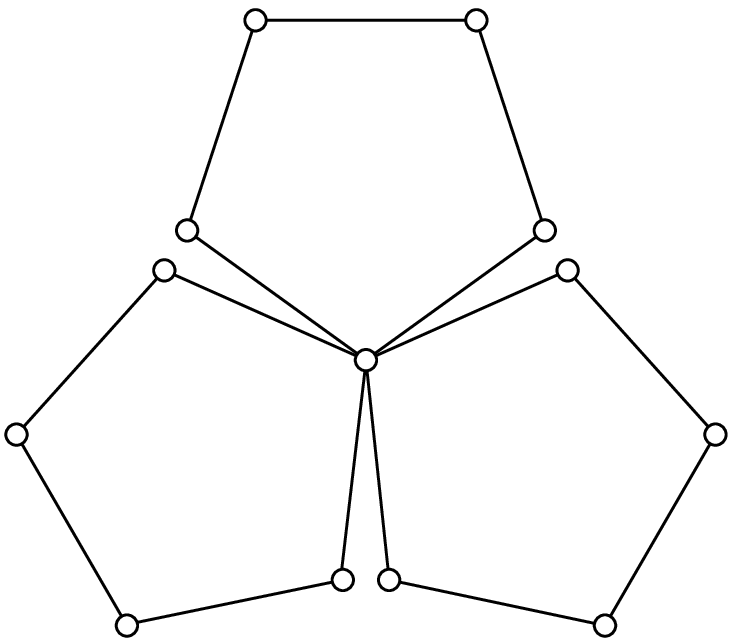, width=4.111cm}
\end{center}

\noindent Commen\c{c}ant avec le graphe des pentagones siamois :

\vspace{2mm}
\begin{center}
\epsfig{file=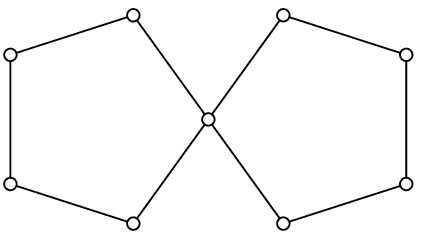, width=4cm},
\end{center}

\noindent j'ai pu vérifier qu'il était pur en faisant, à la main, la liste
exhaustive des graphes non isomorphes de sa classe de parité (cette
liste consiste en 60 graphes, voir l'annexe A). La taille de la classe
de parité croît apparemment de fa\c{c}on exponentielle au fur et à
mesure qu'on augmente le nombre de pentagones : il y a 197 graphes non
isomorphes dans la classe du trèfle, ce nombre grimpe à 571 pour quatre pentagones, et pour cinq pentagones, à 1459. De plus, ces classes sont relativement petites : à titre de comparaison, les classes des cycles d'ordre 9, 13, 17 et 21 contiennent respectivement 23, 138, 1034 et 8957 graphes non isomorphes. Déterminer ces classes à la main n'étant pas envisageable, il a fallu mécaniser le procédé et écrire un programme me permettant d'effectuer ce genre de vérifications à
l'ordinateur (le listing de ce programme est fourni en annexe C). J'ai
ainsi vérifié que le trèfle est pur, tout comme les deux graphes
suivants :
\begin{center}
\epsfig{file=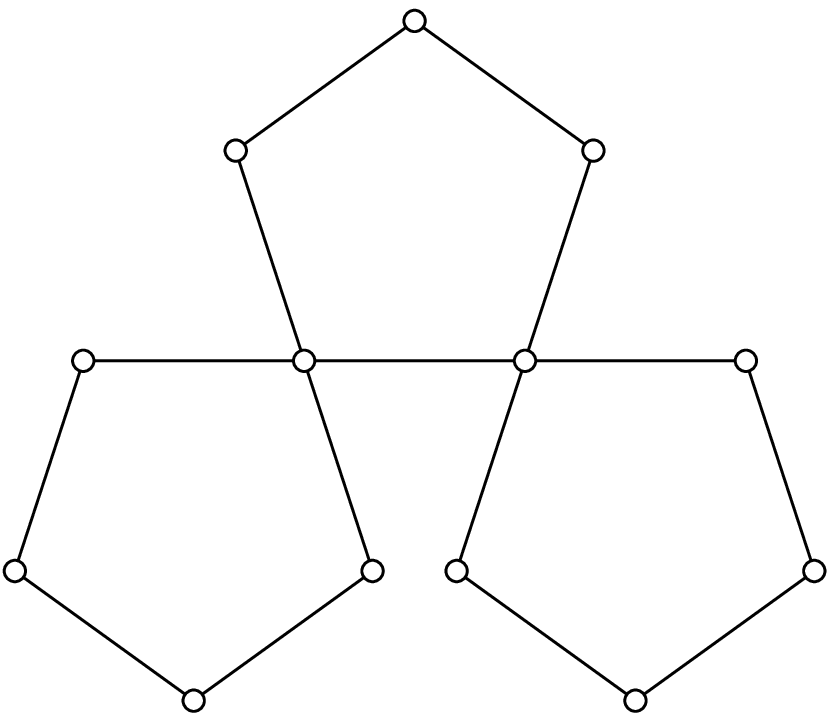, width=4.691cm}
\end{center} 
\begin{center}
\epsfig{file=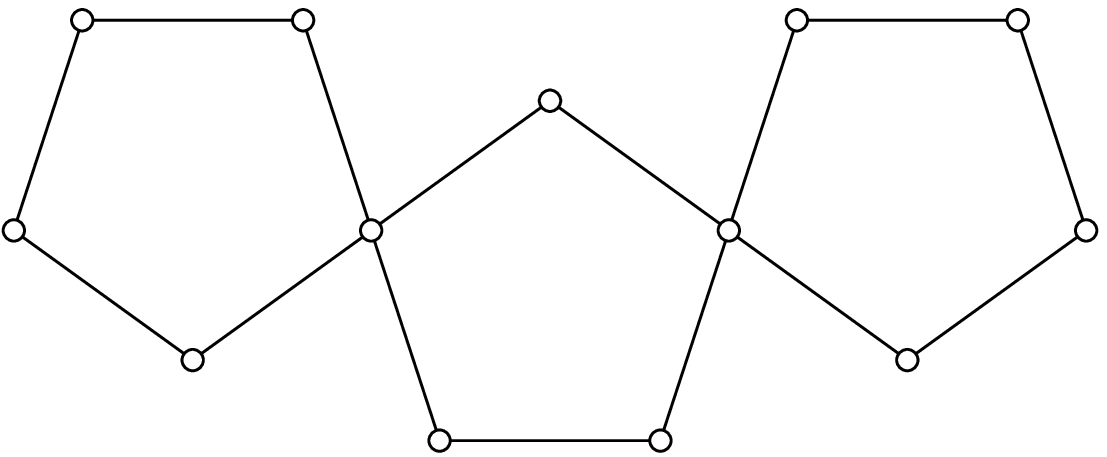, width=6.272cm}
\end{center}
\noindent En fait, même en élargissant la classe de parité en
permettant les complémentations locales aux sommets noirs (impairs),
on ne trouve, pour ces exemples, aucune anticlique noire. Appelant de
tels graphes \emph{fortement purs}, on peut montrer
que l'identification en un sommet de deux graphes fortement purs donne un graphe fortement pur (toujours en prenant le coloriage
naturel).

Il n'y a qu'un seul autre graphe eulérien dans la classe de parité des pentagones siamois, c'est

\begin{center}
\epsfig{file=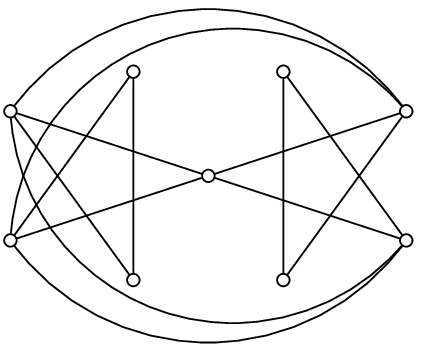, width=4cm},
\end{center} 

\noindent qui est isomorphe à

\begin{center}
\epsfig{file=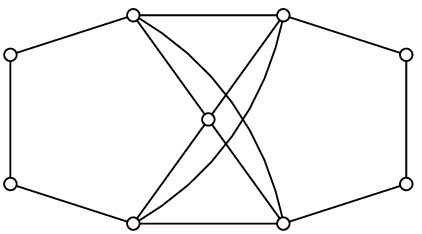, width=4cm},
\end{center} 

\noindent ce qui suggère de tester le graphe

\begin{center}
\epsfig{file=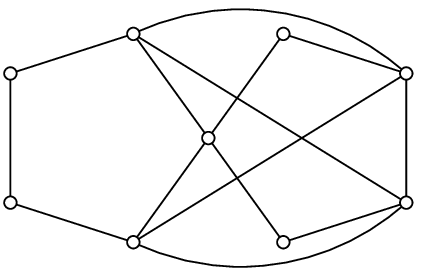, width=4cm}.
\end{center}

\noindent Ce dernier est effectivement pur. Il est intéressant, car
c'est le premier exemple que j'ai trouvé d'un graphe pur qui ne soit pas fortement pur. En testant
systématiquement tous les graphes connexes (naturellement coloriés)
d'ordre $\leq 9$, je suis tombé sur d'autres graphes purs non fortement purs, entre autres les graphes :

\psfrag{et}[][][1]{et}
\begin{center}
\epsfig{file=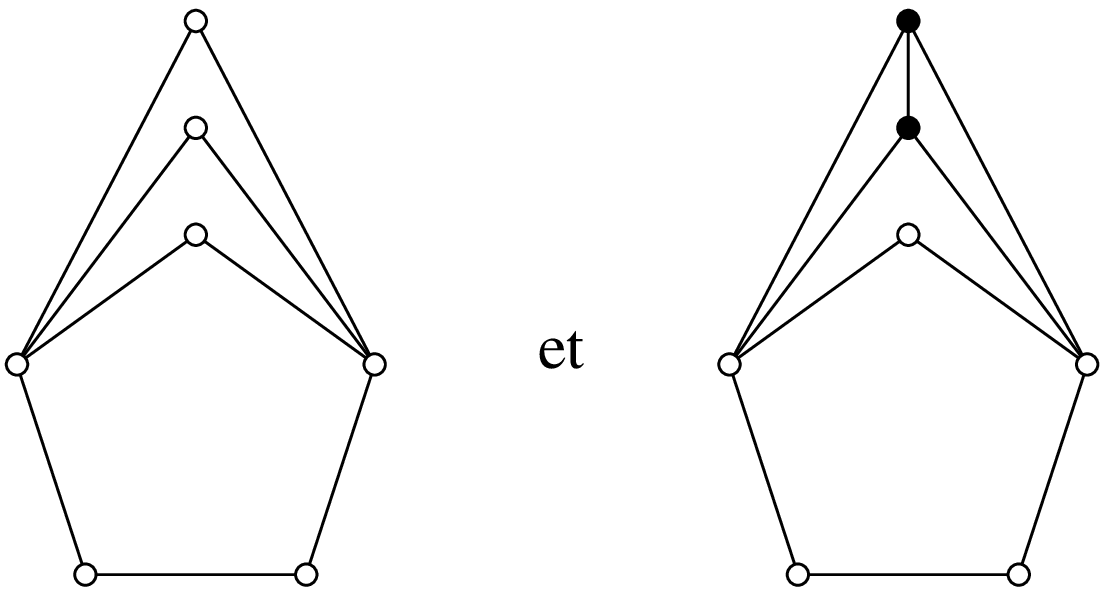, width=6.26cm}.
\end{center}

Ces graphes ont en commun d'avoir une \emph{coupe complète} (une partition $V=V_1\cup V_2$ telle que $|V_1|,|V_2|\geq 2$ et telle que pour tout $u\in V_1$ adjacent à un sommet de $V_2$ et pour tout $v\in V_2$ adjacent à un sommet de $V_1$, $[u,v]\in E$). C'est ce qui m'a mis sur la piste du théorême de
décomposition 3.14 du deuxième article. Dans une première version, ce théorême
ne s'appliquait qu'aux graphes
naturellement coloriés. Un peu plus tard est apparu le graphe pur 

\begin{center}
\epsfig{file=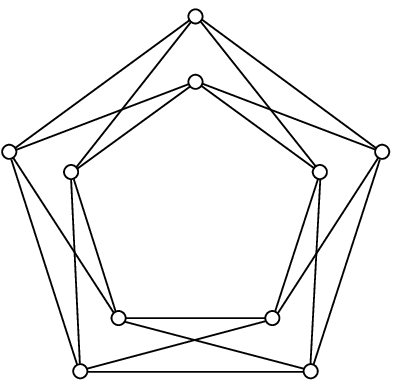, width=4cm}
\end{center}

\noindent qui, afin d'éviter de le considérer comme un graphe
primitif, m'a forcé à élargir mon étude aux bicoloriages
non naturels.

Cette généralisation m'a révélé l'existence d'une
correspondance biunivoque entre les classes de parité de graphes de
cordes arbitrairement bicoloriés
et les systèmes de transitions de graphes 4-réguliers connexes. Le
deuxième article se conclut sur les implications de ces résultats pour
la conjecture de double recouvrement et pose la question : comment
généraliser les systèmes de transitions de graphes 4-réguliers
connexes afin d'obtenir une correspondance avec les classes de parité
(de graphes arbitrairement bicoloriés)?

%% file: words.tex


\chapter[]{Word and set complementation of graphs, invertible graphs.}
\def\proofname{Proof}

\theoremstyle{theorem}
\newtheorem{theorem}{Theorem}[section]
\newtheorem{lemma}[theorem]{Lemma}
\newtheorem{corollary}[theorem]{Corollary}
\newtheorem{conjecture}[theorem]{Conjecture}
\newtheorem{proposition}[theorem]{Proposition}

\theoremstyle{definition}
\newtheorem{definition}[theorem]{Definition}
\newtheorem{notation}[theorem]{Notation}
\newtheorem{problem}{Problem}[section]




\begin{center}
{\bf\Large Fran\c cois Genest}
\end{center}

\section{Abstract}

Local complementation was first introduced as a way
  to establish relationships between Euler tours of an eulerian graph. It also appears in
  isotropic systems. We formalize the notion of substitution rules and
  introduce complementation with respect to sets of vertices in
  bicolored graphs, a concept intimately related to
  orthogonality (or compatibility) of transition systems in
  eulerian graphs and to the Cycle Double Cover Conjecture. Graph
  inversion is also introduced and we show that the inverse of a graph,
  when it exists, has the same automorphism group as the initial graph.

\section{Introduction}
Kotzig [7][8] introduced local complementation when he realized
that all Euler tours of a 4-regular graph could be transformed into
one another by a sequence of re-routings at vertices. The Euler tours
are in correspondence with the members of a \emph{complementation
class} where the re-routings become \emph{complementations} at vertices. Sabidussi
[10] came upon the notion of a \emph{parity class} while working on his
Orthogonality Conjecture (also known as
Sabidussi's Compatibility Conjecture). Sabidussi's approach involves looking at the Euler tours of a 4-regular
graph that are, in a specified way, \emph{orthogonal} to one particular tour. The resulting
subset of Euler tours is what becomes a parity class when translated
into the language of complementation. It is not within the scope of this
paper to present the Orthogonality Conjecture or its more famous relative, the Cycle Double
Cover Conjecture. The interested reader is referred to the surveys
by Jackson [5] and Jaeger [6]; see also the companion paper [4]. Our aim is to
describe complementation and parity classes, using words and sets,
respectively. A natural question will come up: when are the
complementation subsets of a graph in bijection with the graphs in its
parity class? First, we need to
introduce the concepts of local and global subtitution rules. Fon-der-Flaass [2] gave a tight bound on the
diameter of a complementation class. We give an alternative way to
obtain this bound, using substitution rules.

In this article all graphs considered will be simple with edges
represented by unordered pairs of vertices inside brackets. We will be
interested in families of graphs sharing the same vertex set
$V$ of some reference graph $G = (V,E)$. Moreover, for
$A,B\subseteq V$ with $A\cap B=\emptyset$, $K_A$ will stand for the
graph with vertex set $V$ and edge set $\{[u,v]|u,v\in A, u\neq v\}$
and $K_{A,B}$ will have vertex set $V$ and edge set $\{[u,v]|u\in A,
v\in B\}$. If $A=\{u\}$ is a singleton, we omit the
parentheses and write $K_{u,B}$. The symmetric difference of two
graphs with the same vertex set, say $G_1=(V,E_1)$ and $G_2=(V,E_2)$,
will be $G_1\triangle G_2=(V,E_1\triangle E_2)$.

We will also work
with sets of words. Following the notation of semigroup theory, the
set of words in an alphabet $V$ is denoted by $V^*$. The empty word is written
$\epsilon$. The letters $s,t$ will be used to represent words while $u,v,...$ will be vertices.

We use the standard notation $N_H(u)$ for the neighbourhood of a vertex $u$ in the graph $H$, and write $\overline{N_H(u)}$ for the closed neighbourhood of $u$.

\section{\label{section_lc}Local complementation and substitution rules}

\begin{definition} The \emph{(local) complement} at a vertex $u$ of a graph
  $G=(V,E)$ is $Gu=G\triangle K_{N_G(u)}$. In other words, the adjacency relation of $Gu$ coincides with that of $G$ except on $N_G(u)$, where it is replaced by its complement.
\end{definition}

Letting $W(G)=V(G)^*$, complementation is extended recursively to words in $W(G)$ by putting $G\epsilon = G$ and $Gsu=(Gs)u$, where $s\in W(G)$ and $u\in V(G)$. The \emph{complementation class} of $G$ is $\mathcal{C}G=\{Gs|s\in
W(G)\}$.

\begin{notation} \ \newline
(i) $V(s)$ is the set of vertices appearing in the word $s$ (the \emph{support} of $s$);

\noindent (ii) $\lambda (s) = |V(s)|$;

\noindent (iii) Given a word $s$ beginning with the letter $u$, $[s]$ is the word $su$ (for example, $[uv]=uvu$).
\end{notation}

Different words may complement a graph $G$ in the
same way. Accordingly, we define an equivalence relation $\sim_G$ on $W(G)$ by $s\sim_Gs'$ if and only if
$Gs=Gs'$ so that the resulting quotient set, denoted $\Omega (G)$, is in
bijection with $\mathcal{C}G$. Note that $s\sim_Gs' \Rightarrow
st\sim_Gs't$, for any $t\in W(G)$. We are interested in describing
$\mathcal{C}G$ using words. When $G$ is finite, so is
$\mathcal{C}G$. In that case, an obvious first goal would be to find a finite set of
words complementing $G$ to all of $\mathcal{C}G$. To that end, we
introduce the notion of a substitution rule.

Substitution rules describe when a subword can be replaced by another, while ensuring that the resulting word is equivalent to the original word. For example, it is clear that complementing with respect to one vertex twice in a row will result in a graph identical to the original graph. Hence, given a graph $G$ and a word $s=s'uus''\in W(G)$, we know that $s\sim_G s's''$. To put it differently, we can replace the subword $uu$ with the empty word. We want a substitution rule to express this possibility in general terms, without being tied to a particular graph and vertex. To achieve this, the ``symbols'' $u$ and $G$ appearing in the ``rule'' $uu\nobreak \sim_G\nobreak \epsilon$ are understood to be a vertex variable and a graph variable, respectively. Furthermore, for the rule to make any sense when considering actual values of $u$ and $G$, we assume that the variables are tied by the relationship $u\in V(G)$. The following definition attempts to formalize substitution rules just enough for our needs, without resorting to a full description of graphs in terms of logic.

\begin{definition}
Let $s$ be a word on a set of vertex variables, let $G$ be a graph variable with $V(s)\subset V(G)$, and let $P$ be a
logical formula dependent on $G$ (a property of $G$). The couple $R:(P,s)$ is a \emph{substitution rule} if
  $$P(G)\Rightarrow s\sim_G \epsilon.$$

Such a rule will often be
  written $R:P\Rightarrow s\sim_G \epsilon$. Unless otherwise
  specified, we asume rules to be non-trivial, i.e. at least
  one graph $G$ satisfies $P$.
\end{definition}

Given a substitution rule $R:P\Rightarrow s\sim_G\epsilon$,
some fixed graph $G$
and words $s_1=s'ss''$ and $s_2=s's''$ such that $P(Gs')$ is true, we
deduce that $s_1\sim_G s_2$. To emphasize that $R$ was used, we sometimes
write $s_1 \stackrel{R}{\sim}_G s_2$. The following are straightforward rules:

\begin{gather}
uu\sim_G \epsilon, \tag{R1}                         \\         
u\neq v \mbox{ and } [u,v]\notin E(G) \Rightarrow uvuv\sim_G \epsilon. \tag{R2}    
\end{gather}

\begin{proof}
\begin{align} (G\triangle K_{N_G(u)})\triangle
  K_{N_{Gu}(u)} & = G\triangle (K_{N_G(u)}\triangle K_{N_G(u)}) =  G \notag
  \\ (G\triangle K_{N_G(u)})\triangle K_{N_{Gu}(v)} &
  = G\triangle K_{N_G(u)} \triangle K_{N_G(v)} =  G\triangle
  K_{N_G(v)}\triangle K_{N_G(u)} \label{eq1}\\ & = (G\triangle
  K_{N_G(v)})\triangle K_{N_{Gv}(u)} \label{eq2}
\end{align}
Thus given $[u,v]\notin E(G)$, we have $uvuv\sim_G
vuuv \stackrel{R1}{\sim}_G vv \stackrel{R1}{\sim}_G \epsilon.$ 
\end{proof}

Equations \ref{eq1} and \ref{eq2} show that rule R2 is really about the commutativity of some local complementations. While we could define substitution rules to allow the form $[u,v]\notin E(G)\Rightarrow uv\sim_G vu$, the definition we chose seems to be more manageable in the handling of proofs.

Defining the inverse of a word $s$, written $s^{-1}$, as the word
obtained from $s$ by reversing the order of its letters, some direct
consequences of R1 are that for any $s,s',s'',t\in V(G)^*$, we have $tt^{-1}\sim_G \epsilon$ and $s'\sim_{Gs}
s'' \iff ss't\sim_G ss''t$. Unfortunately, $s \sim_G s'$ and $t
\sim_G t'$ do not guarantee that $st \sim_G s't'$. This means that
$\Omega (G)$ cannot be endowed with a group structure using the
concatenation operation. It does have an algebraic structure,
that of an automaton. However this does not seem to be of any help
regarding complementation. The following rule is known [1][10]:

\begin{equation}
  [u,v]\in E \Rightarrow [uv][vu]\sim_G \epsilon \tag{R3}
\end{equation}

This follows from the following lemma, by symmetry between $u$ and $v$:

\begin{lemma} \label{lemme_arete}
If $[u,v]\in E(G)$, then
$$Guvu = G\triangle K_{\{u,v\},V_u\cup V_v} \triangle K_{V_u,V_v} \triangle K_{V_u,V_{uv}}\triangle K_{V_v,V_{uv}},$$

\noindent where $V_u= N(u)\setminus
\overline{N(v)}$, $V_v= N(v)\setminus
\overline{N(u)}$ and $V_{uv} = N(u)\cap N(v)$.
\end{lemma}


\begin{proof}
 Partition $E(G)$ into $E(G)=\{[u,v]\}\cup E(K_{u,V_u} \triangle K_{u,V_{uv}} \triangle K_{v,V_v} \triangle K_{v,V_{uv}}) \cup E'$.
Since $N(u)=\{v\}\cup V_u \cup V_{uv}$ we have
\begin{eqnarray*} E(Gu)& = & E(G) \triangle E(K_{v,V_u} \triangle K_{v,V_{uv}}
  \triangle K_{V_u,V_{uv}}\triangle K_{V_u} \triangle K_{V_{uv}})\\ & = & \{[u,v]\}
  \cup E(K_{u,V_u} \triangle K_{u,V_{uv}} \triangle K_{v,V_u} \triangle K_{v,V_v}) \cup
  E' \\ & & \triangle E(K_{V_u,V_{uv}} \triangle K_{V_u} \triangle K_{V_{uv}} ).
\end{eqnarray*}
From this we see that $N_{Gu}(v)=\{u\}\cup V_u \cup V_v$ and
\begin{eqnarray*}
E(Guv) &=& E(Gu)\triangle E(K_{u,V_u} \triangle K_{u,V_v} \triangle K_{V_u,V_v} \triangle K_{V_u} \triangle K_{V_v}) \\ &=& \{[u,v]\} \cup E(K_{u,V_v} \triangle K_{u,V_{uv}} \triangle K_{v,V_u} \triangle K_{v,V_v}) \cup E' \\ & & \triangle E(K_{V_u,V_v} \triangle K_{V_u,V_{uv}} \triangle K_{V_v} \triangle K_{V_{uv}}).
\end{eqnarray*}
Finally, $N_{Guv}(u) = \{v\}\cup V_v \cup V_{uv}$ and
\begin{eqnarray*}
E(Guvu) &=& E(Guv)\triangle E(K_{v,V_v} \triangle K_{v,V_{uv}} \triangle K_{V_v,V_{uv}}) \\ &=& \{[u,v]\} \cup E(K_{u,V_v} \triangle K_{u,V_{uv}} \triangle K_{v,V_u} \triangle K_{v,V_{uv}}) \cup E' \\ & & \triangle E(K_{V_u,V_v} \triangle K_{V_u,V_{uv}} \triangle K_{V_v,V_{uv}})\\ &=& E(G) \triangle E(K_{\{u,v\},V_u\cup V_v} \triangle K_{V_u,V_v} \triangle K_{V_u,V_{uv}} \triangle K_{V_v,V_{uv}}).
\end{eqnarray*} 
\end{proof}

\begin{proposition}
\begin{equation}[u,v],[v,w],[u,w]\in E\Rightarrow
  [uv][vw][uw]\sim_G \epsilon.  \tag{R4}
\end{equation}
\end{proposition}

\begin {figure} [htb!]
\begin {center}
\psfrag{u}[][][1]{$u$}
\psfrag{v}[][][1]{$v$}
\psfrag{w}[][][1]{$w$}
\psfrag{setu}[][][1]{$\{u\}$}
\psfrag{setv}[][][1]{$\{v\}$}
\psfrag{setw}[][][1]{$\{w\}$}
\psfrag{uv}[][][1]{$\{u,v\}$}
\psfrag{uw}[][][1]{$\{u,w\}$}
\psfrag{vw}[][][1]{$\{v,w\}$}
\psfrag{uvw}[][][1]{$\{u,v,w\}$}
\psfrag{vide}[][][1]{$\emptyset$}
\psfrag{G'}[][][1]{$G'$}
\psfrag{G'[uv]}[][][1]{$G'[uv]$}
\psfrag{G'[uv][vw]}[][][1]{$G'[uv][vw]$}
\psfrag{G'[uv][vw][uw]}[][][1]{$G'[uv][vw][uw]$}
\epsfig {file=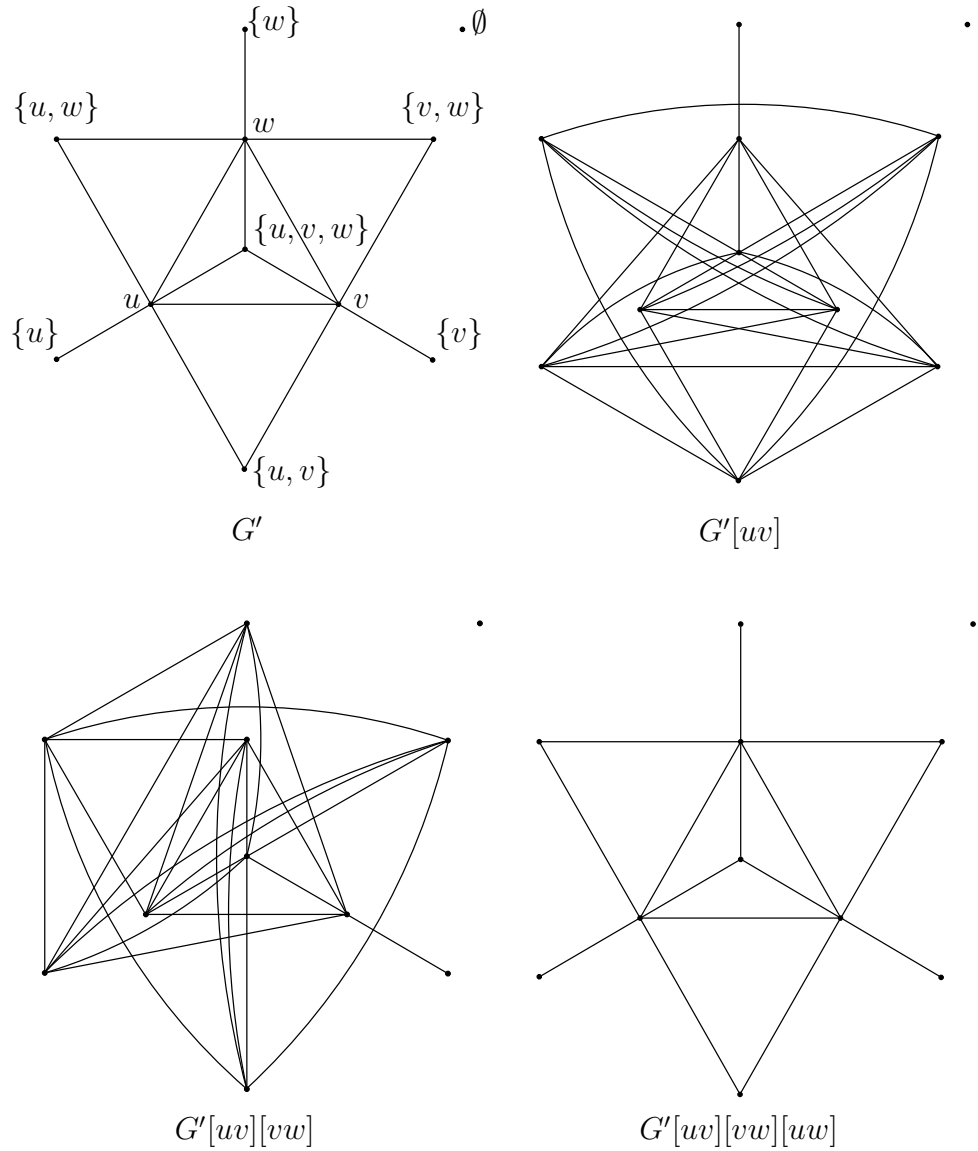,height=15cm}
\end {center}
\caption{\label{preuve_par_dessin}Successive complementations show that $G'=G'[uv][vw][uw]$.}
\end {figure}

\begin{proof}
This follows from Figure~\ref{preuve_par_dessin} and Lemma~\ref{technique}.
\end{proof} 

\begin{lemma} \label{technique}
Let $u,v,w$ be vertices inducing a triangle in $G$. Consider the graph $G'=(V',E')$, where $V'=\{u,v,w\}\cup \mathcal{P}(\{u,v,w\})$ and
        $$E' = \{[u,v],[v,w],[u,w]\} \cup  \{[x,y]|x\in \{u,v,w\}, y\in \mathcal{P}(\{u,v,w\}), x\in y\}$$

and define
\[ \begin{array}{lll} \Phi &:& V \longrightarrow V'\\
        & & x\longmapsto \left\{ \begin{array}{l} x \quad \mbox {if } x\in \{u,v,w\}\\N(x)\cap \{u,v,w\} \quad\mbox{if } x\notin \{u,v,w\}. \end{array} \right. \end{array} \]
Then for any word $s$ in the alphabet $\{[uv],[vw],[uw]\}$
        $$ [x,y]\in E(G\triangle Gs) \iff [\Phi (x),\Phi(y)]\in E(G'\triangle G's)$$
\end{lemma} 

\begin{proof} An easy induction on the length of a word $t$ in the alphabet $\{u,v,w\}$ shows that
\begin{equation} \label{allo} N_{Gt}(x)\cap\{u,v,w\}=N_{G't}(\Phi(x))\cap\{u,v,w\}\quad \mbox{ for all } x\in V.
\end{equation}
As a special case, this is also true for $t$ in the alphabet $\{[uv],[vw],[uw]\}$.
Now use induction on the length of $s$. The assertion is true for the
empty word. Without loss of generality a non-empty $s$ decomposes into $s'[uv]$ with
\begin{equation} \label{allo2} [x,y]\in E(G\triangle Gs') \iff [\Phi(x),\Phi(y)]\in E(G'\triangle G's').
\end{equation}
From (\ref{allo}), $$N_{Gs'}(x)\cap \{u,v\}=N_{G's'}(\Phi(x))\cap
\{u,v\}$$ and $$N_{Gs'}(y)\cap \{u,v\}=N_{G's'}(\Phi(y))\cap
\{u,v\}$$ so, using Lemma~\ref{lemme_arete}, we deduce that
\begin{equation} \label{allo3} [x,y]\in E(Gs'\triangle Gs) \iff [\Phi(x),\Phi(y)]\in E(G's'\triangle G's).
\end{equation}
The result follows from (\ref{allo2}) and (\ref{allo3}).
\end{proof}

\begin{definition} Given a word $s\in W(G)$ and a set of substitution
  rules $\mathcal{R}$, $\mathcal{R}_G(s)$ is the set of words which can be deduced to be equivalent to $s$ using the rules in
  $\mathcal{R}$, i.e., $s'\in \mathcal{R}_G(s) \iff \exists R_1,...,R_n\in \mathcal{R}$ and $s_0,...,s_n \in
W(G)$, $n\geq 0,$ such that $s_0=s,$ $s_n=s',$ and $s_i \stackrel{R_{i+1}}{\sim}_G s_{i+1}$, $i=0,...,n-1$.
\end{definition}

Once local substitution rules will have been introduced in Section~\ref{section_rules},
we will see that R1,...,R4 determine all local rules. Accordingly, we
write $Loc_G(s)$ instead of $\{R1,R2,R3,R4\}_G(s)$.

\begin{definition} Let $\mathcal{R}$ be a set of substitution
  rules. The set of substitution rules \emph{generated} by $\mathcal{R}$, written $\langle\mathcal{R}\rangle $, is the set of rules $(P,s)$ such that for every graph $G$ satisfying $P$, we have $s\in$\ $\mathcal{R}_G(\epsilon)$. $\mathcal{R}$ is \emph{independent} if no
  subset of $\mathcal{R}$ generates $\langle\mathcal{R}\rangle $.
\end{definition} 

\begin{definition} A word $s\in W(G)$ is \emph{bracket-writable} if $s=s_1s_2...s_n$, where each $s_i$ contains at most one repeated letter, in which case it is of the form $s_i=utu=[ut]$, and no letter appears in different $s_i$'s. Such a word is said to be \emph{reduced} if each $s_i$ either consists of a single vertex or can be expressed as $s_i=[u_{i1}u_{i2}]$, where $[u_{i1},u_{i2}]\in E(Gs_1s_2...s_{i-1})$.
\end{definition}

\begin{theorem} \label{reduction}
  For any $s\in W(G)$ and $u\in V$ there exists a reduced
  $s'\in Loc_G(s)$ such that $V(s')\subset V(s)$, and $u$ appears in
  position $1$ or $2$ of $s'$, if at all.
\end{theorem}

\begin{proof} By way of contradiction, suppose that $G$ and $s$ constitute a
  counter-example with $\lambda(s)$ minimal. Clearly $s$ is non-empty. Choose $u_0\in V(s)$, with the restriction that $u_0=u$ if $u\in V(s)$. Writing $\rho (v,t)$ for the position of the last occurence of
  $v$ in $t$, we can suppose without loss of generality that $\rho (u_0,s)\leq \rho (u_0,s')$
  for all $s'\in Loc_G(s)$ such that $V(s')\subset V(s)$. Suppose that $\rho
  (u_0,s)>2$. Consider the subwords of $s$ of length 2 and 3 ending
  with the last occurrence of $u_0$. It is easy to verify that at least
  one of the following sequences of substitutions can be performed (in
  each case, $E$ is meant to be the edge set just prior to complementation
  with respect to the subword considered):

$\\
\mbox{case 1 } \quad u_0u_0\sim \epsilon, \\
\mbox{case 2 } \quad
uu_0\stackrel{R2}{\sim}(u_0uu_0u)uu_0\stackrel{R1}{\sim}u_0u, \mbox
{ if } [u,u_0]\notin E, \\
\mbox{case 3 } \quad
u_0uu_0=[u_0u]\stackrel{R3}{\sim}[uu_0][u_0u][u_0u]\stackrel{R1}{\sim}[uu_0]=uu_0u,
\mbox { if } [u,u_0]\in E, \\
\mbox{case 4 } \quad uuu_0 \stackrel{R1}{\sim}u_0, \\
\mbox{case 5 } \quad vuu_0 \stackrel{case 2}{\sim} uvu_0 \stackrel{case
  2}{\sim} uu_0v, \mbox{ if } [u,v],[v,u_0]\notin E, [u,u_0]\in E, \\
\mbox{case 6 } \quad vuu_0 \stackrel{R4}{\sim}[uu_0][vu_0][vu]vuu_0
\stackrel {R1}{\sim} uu_0uv, \mbox { if } [u,u_0],[u,v],[v,u_0]\in E, \\
\mbox{case 7 } \quad vuu_0 \stackrel{R1}{\sim}uuvuu_0 \stackrel{case
  6}{\sim} uu_0uv, \mbox{ if } [u,u_0],[u,v]\in E, [v,u_0]\notin E, \\
\mbox{case 8 } \quad vuu_0 \stackrel{R1}{\sim}u_0u_0vuu_0 \stackrel {case 6}
{\sim} u_0uvu, \mbox { if } [u,u_0],[v,u_0]\in E, [u,v]\notin E. \\
$

Since this would contradict the minimality of $\rho (u_0,s)$, we must
have $\rho (u_0,s)\linebreak \leq 2$. By the minimality of $\lambda (s)$, $s$
cannot have the prefix $u_0u_0$ ( if $s=u_0u_0s''$ then $s\sim_G s''$
with $\lambda(s'')<\lambda(s)$). If $s=u_0s''$ with $u_0\notin V(s'')$
then by the minimality of $\lambda (s)$ we know that $s''$ can be
replaced by a reduced word not containing $u_0$, resulting in a
reduced word equivalent to $s$, a contradiction. If $s=vu_0s''$
with $[v,u_0]\notin E(G)$ then permuting $v$ and $u_0$ yields the
previous case. Thus $s$ is of the form $s=vu_0s''$ with $[v,u_0]\in
E(G)$ and $u_0\notin V(s'')$. Now consider the graph $G'=Gvu_0$. We
can find a reduced word $t\in Loc_{G'}(s'')$ with $V(t)\subset V(s'')$ and $v$ absent from $t$
or in position 1 or 2. However, $v$ cannot be absent from $t$ or else
$s\sim_G vu_0t$, a reduced word. If $t=vt'$ then $s\sim_G[vu_0]t'$, a
reduced word. The only remaining possibility is $t=wvt'$ with
$[v,w]\in E(G')$ (as before, $[v,w]\notin E(G')$ reduces to an earlier
case). Knowing that $[v,u_0]\in E(G)$ and $[v,w]\in E(Gvu_0)$, we must
have $[w,u_0]\in E(G)$. Using substitutions as in cases 6 or 7,
according to whether $[v,w]\in E(G)$ or not, we have $s\sim_G vu_0wvt' \sim_G
wu_0wvvt' \sim_G [wu_0]t'$. In order to avoid $[wu_0]t'$ being reduced, we
must have $w\in V(t')$. But $t$ is reduced, so $t=wvwt''$ and thus
$s\sim_G [wu_0]wt''\sim_G wu_0t''$, which is reduced. Therefore, no
counter-example exists. 

\end{proof}

Thus to obtain all of $\mathcal{C}G$, we need only look at the reduced
words of $G$, which are finite in number if $G$ is finite. The proof
of the following involves a case analysis as in Theorem \ref{reduction}, and is
left to the reader:

\begin{lemma} \label{reduc_aretes}
  Let $s=[uv][wx]\in W(G)$ be reduced (i.e. $[u,v]\in
  E(G)$ and $[w,x]\in E(G[uv])$), then at least one of
  $[wx][uv],[wu][vx]\mbox { or } [wv][ux]$ is a reduced word in
  $Loc_G(s)$ (according to whether $[w,x],[w,u]$ or $[w,v]\in E(G)$, respectively).
\end{lemma}

\section{The diameter of a complementation graph}

The structure of $\mathcal{C}G$ can be studied in the
\emph{complementation graph} of $G$. This graph, say $H$, is defined by

$$V(H)=\mathcal{C}G,$$

$$[G_1,G_2]\in E(H) \iff G_2=G_1u\mbox{ for some }u\in V(G)$$

One may ask what is the diameter of $H$. \nobreak{Fon-der-Flaass} [2] found a tight bound of $\max \{|V(G)|+1, 10|V(G)|/9\}$ for this diameter when $G$ is connected. This yields an upper bound of $7|V(G)|/6$
in general, which is attained, for example, when all the components of
$G$ are of the form:

\vspace{5mm}

\begin{center}
\epsfig{file=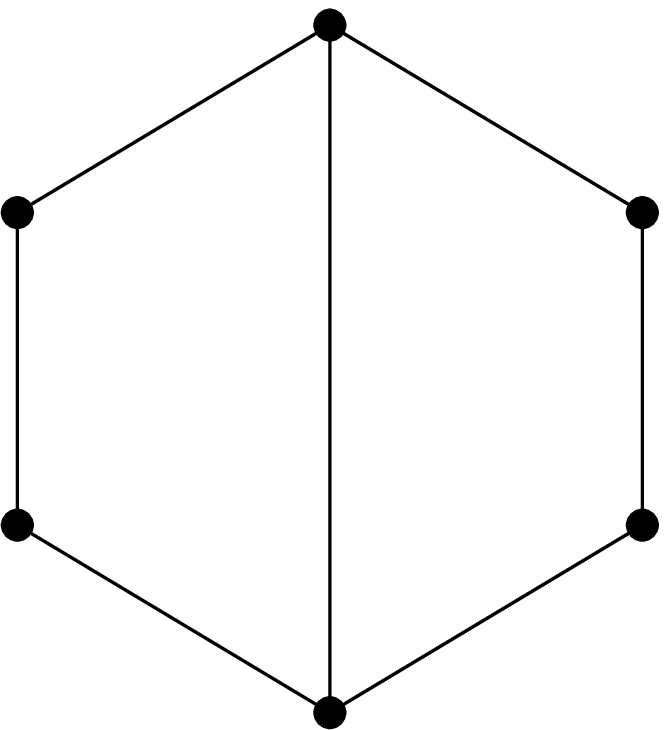,height=2cm,bbllx=0,bblly=0,bburx=190,bbury=208}

\end{center}

In this section, we give an alternative proof of Fon-der-Flaass's bound using substitution rules.

\begin{lemma} \label{mot_minimal}
Let $s=[u_{11}u_{12}][u_{21}u_{22}]...[u_{r1}u_{r2}]$ be a non-empty
reduced word in $W(G)$ such that $V(s)$ induces a connected subgraph
of $G$. For any $u_0\in V(s)$, there exists $s'\in Loc_G(s)$ such that
$l(s')$ (= the length of $s'$) is $\lambda (s)+1$, and the letter $u_0$ occurs both in the first and last
position.
\end{lemma}

\begin{proof}
Choose a bracket-writable $s'=u_0u_1u_2...u_ru_0[v_{11}v_{12}]...[v_{m1}v_{m2}]\in
 Loc_G(s)$ with $V(s')=V(s)$ and $r$
maximal (this exists by Lemma~\ref{reduc_aretes}). Again by Lemma~\ref{reduc_aretes} and R3, since $G|_{V(s')}$ is connected, we can suppose
that $v_{11}$ is adjacent to one of $u_0,...,u_r$. Substituting
repeatedly $u_i[v_{11}v_{12}]$ with $[v_{11}v_{12}]u_i$ (and thus
 moving $v_{11}$ towards the beginning of the word) as long as $[u_i,v_{11}]$, $[u_i,v_{12}]\nobreak \notin E(Gu_0u_1...u_{i\mbox{-}1})$, we
eventually get to make one of the substitutions
 $$u_i[v_{11}v_{12}]\nobreak \sim \nobreak v_{11}v_{12}u_i\quad \mbox{or}\quad u_i[v_{11}v_{12}]\sim v_{12}v_{11}u_i$$ depending on whether or not $[u_i,v_{11}]\in E(Gu_0...u_{i\mbox{-}1})$. This contradicts the maximality of $r$. Thus $s'$
is of the form $u_0u_1...u_ru_0$ with $u_0,...,u_r$ distinct and $l(s')=\lambda(s)+1$.
\end{proof}
  
\begin{lemma} \label{mot_canonique}
For any $s\in W(G)$, there exists
$s'=u_1u_2...u_r[v_{11}v_{12}]...[v_{m1}v_{m2}]\linebreak \in Loc_G(s)$, a
reduced word such
that $V(s')\subset V(s)$ and with no edge in $G$ between
$\{u_1,...,u_r\}$ and $\{v_{11},v_{12},...,v_{m1},v_{m2}\}$.
\end{lemma}

\begin{proof}
By Theorem~\ref{reduction}, we can choose a reduced $s'\in Loc_G(s)$ such that
 $V(s')\subset V(s)$. Also, we can impose the condition that the position of
 the first occurrence of a double occurrence letter is
 maximal, i.e. $s'$ is bracket-writable as
 $u_1u_2...u_r[v_{11}v_{12}]s_2...s_m$, where $u_1,...,u_r$ are distinct, each $s_i$ stands for $v_i$ or $[v_{i1}v_{i2}]$, and $r$ is maximal. Suppose $s_{i+1}=v\in V(G)$ for some $i$. Without loss of generality $i$ is
 minimal, but then one of the following substitutions can be
 performed:

$$[v_{i1}v_{i2}]v \sim v[v_{i1}v_{i2}],\quad [v_{i1}v_{i2}]v
\sim vv_{i1}v_{i2}, \quad \mbox{or}\quad [v_{i1}v_{i2}]v \sim vv_{i2}v_{i1}.$$

\noindent This contradicts the minimality of $i$ (or the maximality of $r$, if
 $i=1$). Thus $s'=u_1...u_r[v_{11}v_{12}]...[v_{m1}v_{m2}]$. Suppose
 that $s'$ does not satisfy the desired conditions. Then some $u_i$ and some $v_{jk}$ must be adjacent in $Gu_1...u_r$. By Lemma~\ref{reduc_aretes} and R3,
 we can suppose that $v_{jk}=v_{11}$. Just as in the proof of Lemma~\ref{mot_minimal}, transform $s'$, in successive steps, by replacing $u_j[v_{11}v_{12}]$ with $[v_{11}v_{12}]u_j$ as long as there is a subword of the form $u_j[v_{11}v_{12}]$, where $[u_j,v_{11}],[u_j,v_{12}]\notin E(Gu_1...u_{l-1})$. Eventually, this produces a subword of the form $u_j[v_{11}v_{12}]$ which can be replaced by $v_{11}v_{12}u_j$ or $v_{12}v_{11}u_j$. The existence of the resulting reduced word contradicts the maximality of $r$.
\end{proof}

\begin{theorem} \label{diametre}
If $H$ is the complementation graph of a connected
graph $G$ of finite order $\neq 6$, and $d$ is the diameter of $H$, then $d\leq 10|V(G)|/9$.
\end{theorem}

\begin{proof}
Without loss of generality, we can choose $G$ and $s\in W(G)$ so that
$Gs$ is at distance $d$ from $G$, $s$ is reduced and $\lambda(s)$ is
minimal. Because of R2, we can suppose that the components of the subgraph $G'$ of $G$ induced by
$V(s)$ have the vertex sets $V(s_1),...,V(s_k)$, respectively, where
$s=s_1s_2...s_k$. Note that any permutation of the $s_i$'s would yield
an equivalent word. By Lemma~\ref{mot_canonique}, we can also suppose
that each
$s_i$ either consists of non-repeating letters or satisfies the
conditions of Lemma~\ref{mot_minimal} (in which case $\lambda(s_i)$ is
even). In the following, $N(t)$ stands for $N_G(V(t))$. We now describe an algorithm to modify $s$ in steps that preserve
equivalence and result in a word $s'$ satisfying $l(s')\leq
10|V(G)|/9$:\\

Step 1: Set $I=\{1,...,k\}$ and $S=V(G)\setminus V(s)$.\\

Step 2: For each $i\in I$ such that
$l(s_i)=\lambda(s_i)$, modify $I$ by removing $i$.\\

Step 3: While there exists an $i\in I$ and $u\in N(s_i)\cap S\setminus
\underset {i\neq j\in I}{\bigcup} N(s_j)$, replace $s_i$ with an equivalent word
given by Lemma~\ref{mot_minimal} and remove $i$ from $I$ and $u$ from
$S$.\\

At this point, for $i\in I$ and $u\in N(s_i)\cap S$, there is
necessarily a $j\in I\setminus \{i\}$ such that $u\in N(s_j)$. Note that
$V(s_i)\cup \{u\}$ induces a connected subgraph of $Gu$, and by an
argument along the lines of the proof of Lemma~\ref{mot_canonique} there exists a word, which we will denote by
$s_i'$, such that $V(s_i')=V(s_i)$, $l(s_i')=\lambda (s_i')$ and
$uus_i\sim_G us_i'u$. Given $u\in N(s_i)\cap S$ with $i\in I$, let
$I_u=\{j\in I|u\in N(s_j)\}$.\\

Step 4: If there exists some $u\in N(s_i)\cap S$ with $i\in I$ such
that $\underset{j\in I_u}{\sum} \lambda(s_j)\geq 8$, bring the corresponding
$s_i$'s to the beginning of $s$ using R2, and relabel so that
$I_u=\{1,...,l\}$; then make the substitution indicated by
\[s_1s_2...s_l\sim_G us_1'uus_2'u...us_l'u\sim_G us_1's_2'...s_l'u,\]
remove the indices in $I_u$ from $I$, $u$ from $S$ and go back to step
3.\\

Given $i\in I$ such that $\lambda(s_i)=2$ (i.e. $s_i=uvu$ for some $[u,v]\in
E(G)$), we must have that $N(s_i)\cap S$ contains at least two
vertices. This follows from the minimality of $\lambda(s)$ and by the
substitution rules:

\begin{gather}
\mbox{degree}(u)\leq 1 \Rightarrow u\sim_G \epsilon, \tag{R5} \\
N(u)=N(v) \Rightarrow uv\sim_G \epsilon. \tag{R6} 
\end{gather}

Step 5: If there are $i,j\in I$, $i\neq j$, and distinct vertices $u$ and $v$
such that \[N(s_i)\cap
N(s_j)\cap S \setminus \underset{\substack{l\in I\\ i\neq l \neq
    j}}{\bigcup} N(s_l)=\{u,v\},\] replace $s_i$ and $s_j$ with
equivalent words as given by Lemma~\ref{mot_minimal}, remove $i$ and $j$
from $I$, $u$ and $v$ from $S$, and go back to step 3.\\

\psfrag{u}[][][.7]{$u$}
\psfrag{v}[][][.7]{$v$}
\psfrag{w}[][][.7]{$w$}
\psfrag{x}[][][.7]{$x$}
\psfrag{V(s1)}[][][.7]{$V(s_1)$}
\psfrag{V(s2)}[][][.7]{$V(s_2)$}
\psfrag{V(s3)}[][][.7]{$V(s_3)$}
\psfrag{V(s4)}[][][.7]{$V(s_4)$}
\psfrag{V(s5)}[][][.7]{$V(s_5)$}
\psfrag{V(s6)}[][][.7]{$V(s_6)$}
\psfrag{V(s7)}[][][.7]{$V(s_7)$}
\psfrag{V(s8)}[][][.7]{$V(s_8)$}
\psfrag{V(s9)}[][][.7]{$V(s_9)$}
\psfrag{Ix1234567}[][][.7]{$I_x\subset \{1,2,3,4,5,6,7\}$}
\psfrag{Iwx1234567}[][][.7]{$I_w,I_x\subset \{1,2,3,4,5,6,7\}$}
\psfrag{Ix12345678}[][][.7]{$I_x\subset \{1,2,3,4,5,6,7,8\}$}

\begin{figure}[htb!]
\begin{center}
\epsfig{file=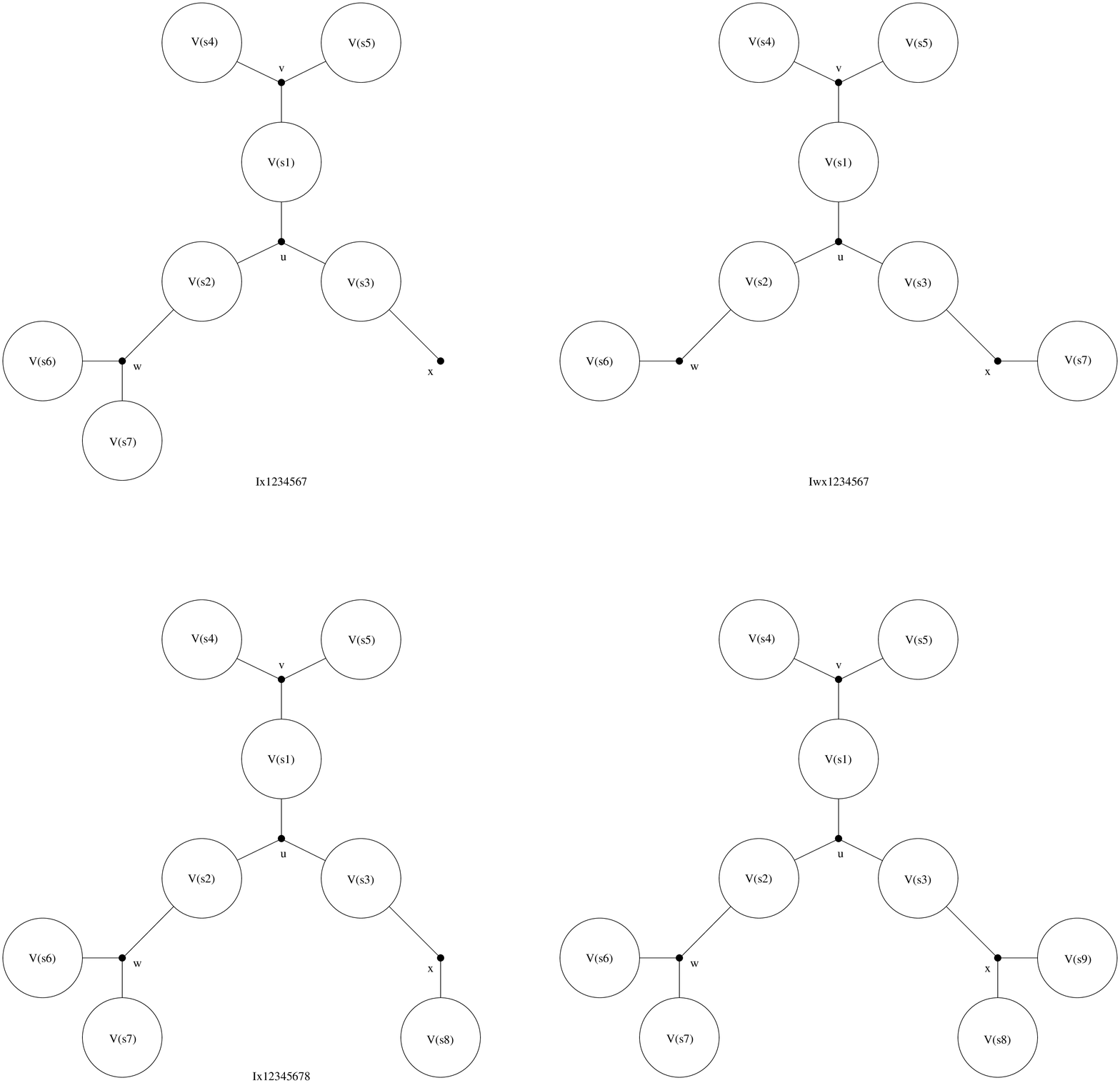,height=13cm}

\caption{\label{arbre1}}
\end{center}
\end{figure}



Step 6: If there exists $i\in I$ with $\lambda(s_i)=4$, let $u\in
N(s_i)\cap S$. Because of step 4, $|I_u|=2$ and, using R2
and relabeling, we can suppose that $I_u=\{1,2\}$ with
$\lambda(s_1)=4$ and $\lambda(s_2)=2$. Because of step 5, we can
assume that there is $v\in N(s_2)\cap S$ with $v\neq u$ and $v\notin
N(s_1)$. Let $I_v=\{2,3,...,l\}$. We replace $s_1s_2...s_l$ with an equivalent
word $s_1vs_2'...s_l'v$ in the same manner as noted earlier, remove $1,...,l$ from
$I$, $u$ and $v$ from
$S$, and go back to step 3.\\

At this point, for each $i\in I$, we have $\lambda(s_i)=2$ and, for
each $u\in N(s_i)$, $|I_u|=2$ or $3$. If for each $i\in I$ and each $u\in N(s_i)\cap S$ we have
$|I_u|=2$, we can stop, since then $|S|\geq |I|$.\\

Step 7: If for some $i\in I$ and $u\in N(s_i)\cap S$ we have
$|I_u|=2$, we can suppose (with the appropriate use of R2 and
relabeling) that $I_u=\{1,2\}$ and $I_v=\{2,3,4\}$, for some $v$. Now replace $s_1s_2s_3s_4$ with $s_1vs_2's_3's_4'v$,
remove $1,2,3,4$ from $I$, $u$ and $v$ from $S$, and go back to step 3.\\

Step 8: If there is a vertex $u\in S$ and $i\in I_u$ for which some $v\in N(s_i)\cap
S\setminus\{u\}$ has $I_v=I_u$, with say $I_u=\{1,2,3\}$ after an appropriate use of R2 and relabeling, replace $s_1s_2s_3$ with $us_1's_2's_3'u$,
remove $1,2,3$ from $I$, $u$ and $v$ from $S$ and go back to step 3.\\

Step 9: If there is a vertex $u\in S$ and $i\in I_u$ for which some $v\in N(s_i)\cap
S\setminus\{u\}$ has $|I_v\cap I_u|=2$, say, without loss of generality, $I_u=\{1,2,3\}$ and $I_v=\{2,3,4\}$, replace
$s_1s_2s_3s_4$ with $us_1's_2's_3'us_4$, remove $1,2,3,4$ from $I$,
$u$ and $v$ from $S$, and go back to step 3.\\

Now we can suppose that for some $u,v\in S$, we have $I_u=\{1,2,3\}$, $I_v=\{1,4,5\}$ and, because of steps 8 and 9, there are distinct vertices $w,x$ such that $w\in N(s_2)\cap (S\setminus\{u\})$ and $x\in N(s_3)\cap (S\setminus\{u\})$.\\

Step 10: If $I_w\subset \{1,2,3,4,5\}$ and $I_x\subset \{1,2,3,4,5\}$,
replace $s_1s_2s_3s_4s_5$ with $s_2s_3vs_1's_4's_5'v$, remove
$1,2,3,4,5$ from $I$, $u,v,w$ and $x$ from $S$, and go back to step
3.\\

We can now suppose that $6\in I_w$.\\

Step 11: If $I_w,I_x\subset \{1,2,3,4,5,6\}$, replace
$s_1s_2s_3s_4s_5s_6$ with $vs_1's_4's_5'vs_2s_3s_6$, remove $1$
through $6$ from $I$, $u,v,w$ and $x$ from $S$, and go back to step
3.\\

From this point on, we are considering subwords with at least eighteen distinct
vertices, so that the length of a replacement word can exceed by two the number of distinct
letters and still avoid any possible violation of
$\lambda(s')\leq 10|V(G)|/9$.\\

Step 12: It should be clear by now how to modify our word $s$ in each
of the situations depicted in Figure~\ref{arbre1}. In each case,
modify $I$ and $S$ appropriately and go back to step 3.

\end{proof}

\begin{figure} [htb!]
\begin{center}
\epsfig{file=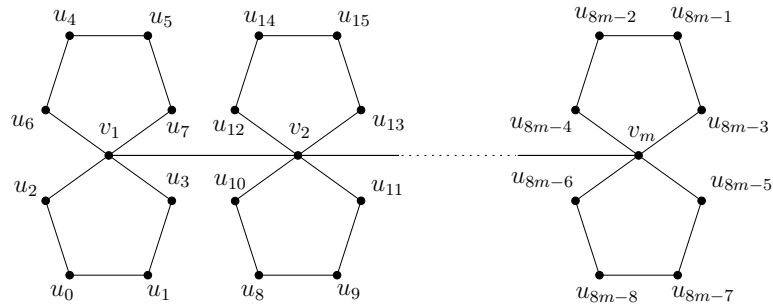,height=4cm}

\end{center}
\caption{\label{gros_diam}A connected graph $G$ whose complementation graph has diameter $10|V(G)|/9$ (also found in \emph{[2]}).}
\end{figure}

\begin{proposition}
Let $H$ be a pair of pentagons sharing a vertex. Let $G$ consist of
$m$ copies of $H$ together with a path going through all the
cut-vertices, as shown in Figure~\ref{gros_diam}. Then the diameter of the complementation graph of $G$ is $10|V(G)|/9$.
\end{proposition}

\psfrag{u0}[][][.8]{$u_0$}
\psfrag{u1}[][][.8]{$u_1$}
\psfrag{u2}[][][.8]{$u_2$}
\psfrag{u3}[][][.8]{$u_3$}
\psfrag{u4}[][][.8]{$u_4$}
\psfrag{u5}[][][.8]{$u_5$}
\psfrag{u6}[][][.8]{$u_6$}
\psfrag{u7}[][][.8]{$u_7$}
\psfrag{u8}[][][.8]{$u_8$}
\psfrag{u9}[][][.8]{$u_9$}
\psfrag{u10}[][][.8]{$u_{10}$}
\psfrag{u11}[][][.8]{$u_{11}$}
\psfrag{u12}[][][.8]{$u_{12}$}
\psfrag{u13}[][][.8]{$u_{13}$}
\psfrag{u14}[][][.8]{$u_{14}$}
\psfrag{u15}[][][.8]{$u_{15}$}
\psfrag{v1}[][][.8]{$v_1$}
\psfrag{v2}[][][.8]{$v_2$}
\psfrag{v3}[][][.8]{$v_3$}
\psfrag{v4}[][][.8]{$v_4$}
\psfrag{v5}[][][.8]{$v_5$}
\psfrag{v6}[][][.8]{$v_6$}
\psfrag{v7}[][][.8]{$v_7$}
\psfrag{v8}[][][.8]{$v_8$}
\psfrag{vm}[][][.8]{$v_m$}
\psfrag{w}[][][.8]{$w$}
\psfrag{u8m-8}[][][.8]{$u_{8m-8}$}
\psfrag{u8m-7}[][][.8]{$u_{8m-7}$}
\psfrag{u8m-6}[][][.8]{$u_{8m-6}$}
\psfrag{u8m-5}[][][.8]{$u_{8m-5}$}
\psfrag{u8m-4}[][][.8]{$u_{8m-4}$}
\psfrag{u8m-3}[][][.8]{$u_{8m-3}$}
\psfrag{u8m-2}[][][.8]{$u_{8m-2}$}
\psfrag{u8m-1}[][][.8]{$u_{8m-1}$}

\begin{proof}
Let $s=[u_0u_1][u_2u_3][u_4u_5][u_6u_7]...[u_{8m-4}u_{8m-3}][u_{8m-2}u_{8m-1}]$. Let $s'$ be a shortest word in $W(G)$ such that $s'\sim_G s$. Choose a pentagon of $G$: for example, the pentagon containing $u_0$. Identify all the vertices not on the pentagon with a new vertex $w$. Remove all loops, and identify all multiple edges. We obtain the graph $G'$:

\vspace{5mm}
\begin{center}
\epsfig{file=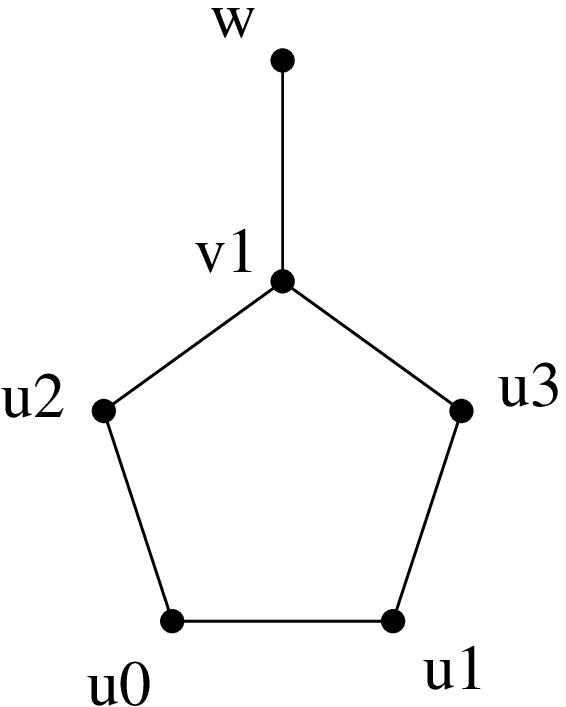,height=3.3cm}

\end{center}
\vspace{5mm}

A word in $W(G)$ induces a word in $W(G')$ in a unique way: send $\epsilon$ to $\epsilon$, and given $tx\in W(G)$ such that $t$ is sent to $t'$, send $tx$ to\\

(i) $t'x$ if $x$ is on the pentagon; 

(ii) $t'w$ if $x$ is not on the pentagon and $x\in N_{Gt}(\{u_0,u_1,u_2,u_3,v_1\})$;

(iii) $t'$ if $x$ is not on the pentagon and $x\notin N_{Gt}(\{u_0,u_1,u_2,u_3,v_1\})$.\\

In this way, $s'$ is sent to a word $t\in W(G')$ such that $G't$ is:

\vspace{5mm}
\begin{center}
\epsfig{file=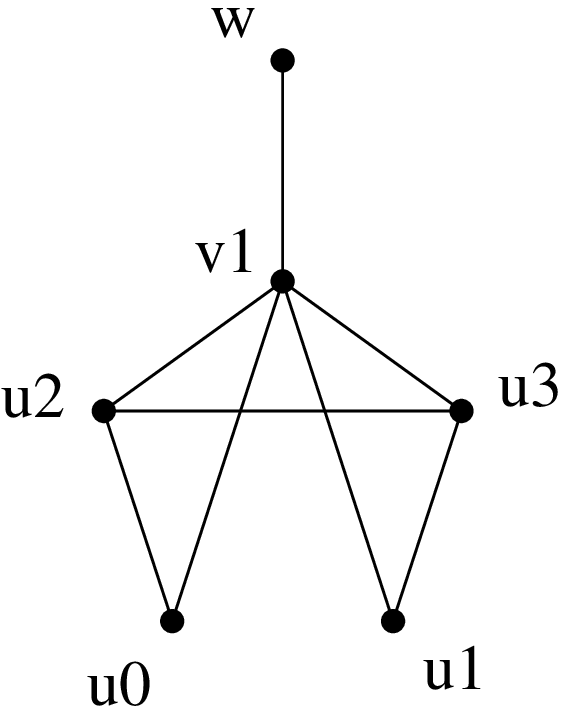, height=3.3cm}

\end{center}
\vspace{5mm}

It can be verified that for any reduced $t'\in W(G')$ such that $t'\sim_{G'}t$, we have $V(t')=\{u_0,u_1,u_2,u_3\}$ or $V(t')=\{u_0,u_1,u_2,u_3,w\}$. This implies that neither $t$ nor $s'$ contain exactly one occurrence of $v_1$. If $v_1\notin V(s')$, then $u_0,u_1,u_2,u_3$ contribute 5 to the length of $s'$ (i.e., $t$ is one of the words $u_0u_1u_3u_2u_0$, $u_1u_0u_2u_3u_1$, $u_2u_0u_1u_3u_2$, or $u_3u_1u_0u_2u_3$). If $v_1$ appears at least twice in $s'$, then the vertices of the two pentagons joined at $v_1$ contribute at least 10 to the length of $s'$. Therefore, $l(s')\geq 10m$, so the diameter of the complementation graph is precisely $10m=10|V(G)|/9$.
\end{proof} 

\section{\label{section_rules}Local and global substitution rules}

Going back to the rules R5 and R6 presented in the proof of Theorem~\ref{diametre}, note that each logical formula takes into account
adjacencies involving every vertex of the graph, not just the adjacencies
within $V(s)$.

\begin{definition} A substitution rule $P\Rightarrow s\sim_G \epsilon$ is
  \emph{local} if for any two graphs $G$ and $H$ such that $G$ is an
  induced subgraph of $H$, $$P(G)\Rightarrow s\sim_H \epsilon.$$
A non-local rule is \emph{global}.
\end{definition}

Thus R1 to R4 are local rules and R5 and R6 are global rules. From the
definition follows that:

\begin{proposition} \label{local}
Local rules generate local rules.
\end{proposition}

\begin{proposition} \label{global}
A substitution rule $P$ $\Rightarrow s\sim_G \epsilon$, where $s$ is a
non-empty reduced word, is global.
\end{proposition}

\begin{proof}
Let $G=(V,E)$ be a graph satisfying $P$. Let
$u,v\notin V$. If $s$ ends with a single occurrence letter $w$, let
$E'=E\cup\{[u,w],[v,w]\}$. If not, we have $s=s'wxw$ with $w,x\notin
V(s')$, in which case let $E'=E\cup\{[u,w],[v,x]\}$. By construction, G is an induced
subgraph of $H=(V\cup\{u,v\},E')$ but $[u,v]\notin H$ while $[u,v]\in
Hs$. Thus $Hs\neq H$.
\end{proof}

\begin{proposition} The rules 
\begin{gather}
uu\sim_G \epsilon, \tag{R1}                         \\         
u\neq v \mbox{ and } [u,v]\notin E(G) \Rightarrow uvuv\sim_G \epsilon, \tag{R2} \\   
[u,v]\in E \Rightarrow [uv][vu]\sim_G \epsilon, \tag{R3} \\
[u,v],[v,w],[u,w]\in E\Rightarrow [uv][vw][uw]\sim_G \epsilon,  \tag{R4}
\end{gather}

\noindent form an independent generating set of the local rules.
\end{proposition}

\begin{proof}
We first show independence. Consider
$G=(\{u,v\},\{[u,v]\})$. Since $\mathcal{C}G=\{G\}$,
$\{$R1,R2,R4$\}_G(\epsilon)= \{$R1$\}_G(\epsilon )$ and any word in
$\{$R1$\}_G(\epsilon)$ will contain an even number of occurrences of the letter
$u$. Thus any independent generating subset of the four rules must
contain R3.

Now let $G=(\{u,v,w\},\{[u,v],[u,w],[v,w]\})$. Any word in
$\{$R1,R2,R3$\}_G(\epsilon)$ has an even number of letters, so that
R4 is also essential. If we let $G=(\{u\},\emptyset)$, then
$\{$R2,R3,R4$\}_G(\epsilon)= \{\epsilon \}$, thus R1 is
essential. Finally, let $G=(\{u,v\},\emptyset)$. Defining the total order $u<v$ on $V(G)$, let the sign of a word $s=u_1u_2...u_n$ in $W(G)$ be $\sigma(s)= (-1)^{\alpha}$ where
$\alpha=\operatorname{card} \{(i,j)|u_i<u_j,1\leq i<j\leq n\}$. By induction on the length of words, we can show that for any $s\in \{$R1,R3,R4$\}_G(\epsilon)=\{$R1$\}_G(\epsilon)$, we have $\sigma(s)=1$. Since $\sigma(uvuv)=-1$, this completes the proof of independence.

We know from Proposition~\ref{local} that $\langle\{$R1,R2,R3,R4$\}\rangle $ is a set of local
rules. Let $P\Rightarrow s\sim_G \epsilon$ be a local rule. Consider
$V=V(s)$ and $\{G_i\}_{i\in I}$ the
  family of graphs on the vertex set $V$ satisfying $P$. By Theorem~\ref{reduction}, for any rule
  $\mathcal{R}_i:G|_V=G_i\Rightarrow s\sim_G \epsilon$ there exists a rule
  $\mathcal{R}'_i:G|_V=G_i\Rightarrow s'\sim_G \epsilon$ in
  $\langle\{\mathcal{R}_i,$R1,R2,R3,R4$\}\rangle $ where $s'$ is reduced. Since
  $\mathcal{R}'_i$ is local, Proposition~\ref{global} forces $s'=\epsilon$. Thus
  $\mathcal{R}'_i\in$\ $\langle\{$R1,R2,R3,R4$\}\rangle $, and since every substitution
  in the proof of Theorem~\ref{reduction} is reversible, we have
  $\mathcal{R}_i\in$\ $\langle\{$R1,R2,R3,R4$\}\rangle $. Since $P\Rightarrow
  s\sim_G \epsilon$ is generated by the $\mathcal{R}_i$'s, we conclude
  that it is in $\langle\{$R1,R2,R3,R4$\}\rangle $.
\end{proof}

\begin{proposition}
If $s,s'$ are reduced words such that $s'\in Loc_G(s)$, then $V(s)=V(s')$.
\end{proposition}

\begin{proof}
Suppose, by way of contradiction, that there is a $u\in V(s)\setminus V(s')$. Then
$s's^{-1}\sim_G \epsilon$ and a reduced word $t\in Loc_G(s's^{-1})$
given by Theorem~\ref{reduction} will contain $u$. But this would mean that local
rules generate a global rule of the form $P\Rightarrow t\sim_G
\epsilon$, contradicting Proposition \ref{local}.
\end{proof}

Looking for a new independent (global) rule $P\Rightarrow s\sim_G
\epsilon$, we can suppose, by Lemma~\ref{mot_canonique}, that $s=s_1...s_k$ is reduced, where each $V(s_i)$ induces a component of
$G|_{V(s)}$, and each $s_i$ is of the form $u_1...u_m$ (non-repeating
letters) or $[u_{11}u_{12}]...[u_{m1}u_{m2}]$ (reduced). Given a component induced by $s$, say $G'=G|_{V(s_i)}$, we have $s_i\sim_{G'}\epsilon$. Therefore, we ask:
 
\begin{problem} \label{prob1}
For which connected graphs $G=(V,E)$ of minimal degree $>1$
without twins (vertices $u,v$ such that $N(u)=N(v)$ or $\overline{N(u)}=\overline{N(v)}$), together with a reduced word $s$ such that $V(s)=V$, do
we have $s\sim_G \epsilon$?
\end{problem}

We will see in section 6 that there is a family of
graphs satisfying the conditions of Problem~\ref{prob1}. 
\section{Complementation sets}

Notice that if $u$ is of odd degree, then the degree of any vertex of $G$ has the
same parity in $G$ and $Gu$. If $u$ is of even degree, then the
parity of the degrees in $G$ and $Gu$ differs precisely over the
neighbours of $u$. Hence, if we colour the vertices of even degree white and the others black, and if we change the colours of
the neighbours of $u$ whenever we complement at a white vertex $u$, then
the colours agree with the parity of the degrees for each graph in
$\mathcal C G$. This will be called the \emph{natural colouring} of $G$. In
the following, a \emph{bicolouring} will always mean a $\{$black,
white$\}$-colouring.

\begin{definition}
The \emph{(local) complement} of a bicoloured graph $G$ with respect to a
vertex $u$ is a bicoloured graph $Gu$ such that
 $$V(Gu)=V(G),$$
$$E(Gu)=E(G)\triangle E(K_{N_G(u)}),$$
with its bicolouring defined to be the same as that of $G$
if $u$ is black in $G$; if $u$ is white in $G$, then it is obtained
from the bicolouring of $G$ by reversing the colours of the vertices
in $N_G(u)$.
\end{definition}

Complementation with respect to
words in the alphabet $V(G)$ is extended in the natural manner.\\

For the purposes of the next definition, call a set of words $W\subset
V(G)^*$ $parity$ $closed$ if\\

\noindent (i) $W$ contains the empty word;

\noindent (ii) if
$s\in W$ and $u$ is a white vertex of $Gs$, then $su\in W$;

\noindent (iii) if
$s\in W$ and $u,v$ are adjacent black vertices of $Gs$, then
$s[uv],s[vu]\in W$.\\

Clearly the intersection of parity closed sets is
parity closed, hence there is a smallest parity closed set, denoted by
$W^{\circ}(G)$. The words in $W^{\circ}(G)$ will be referred to as \emph{parity words}.

\begin{definition}
The \emph{parity class} of a bicoloured graph $G$ is
$$[G]=\{Gs|s\in W^\circ(G)\}.$$
\end{definition}

Following the idea of section~\ref{section_lc}, two parity words $s$ and $s'$ will be
\emph{equivalent} ($s\sim_G s'$) if $Gs=Gs'$. By verifying it for R1 to R4, we can show that local
substitution rules are valid for bicoloured graphs (i.e. if $G$ is a bicoloured graph with underlying (uncoloured) graph $H$ and $s'\in Loc_H(s)$, then $Gs=Gs'$).

\begin{theorem}
Two reduced parity words $s$ and $t$ with $V(s)=V(t)$ are equivalent.
\end{theorem}

\begin{proof}
Use induction on $\lambda(s)$. If $s=us'$ with $u\notin V(s')$, then $u$ is
white in $G$, and applying Theorem~\ref{reduction} to $t$, $t\sim_G
ut'$. If
$s$ is of the form $[uv]s'$ then $u$ is black in $G$, and applying
Theorem~\ref{reduction} to $t$, $t\sim_G [uw]t'$. If $w\neq v$, apply Theorem~\ref{reduction} again to
get $t\sim_G [uw][vx]t''$ and finally, from Lemma~\ref{reduc_aretes}, $t\sim_G [uv][wx]t''$. By changing
the reference graph to $Gu$ or $G[uv]$ accordingly, the problem
reduces to words for which $\lambda$ is smaller.
\end{proof}

Given a bicoloured graph $G$, we are now justified to speak about
complementation with respect to subsets of $V$.

\begin{definition}
A set $S\subset V(G)$ is a \emph{complementation set} of a bicoloured graph $G$ if there
exists a reduced word $s\in W^\circ(G)$ such that $V(s)=S$. In that case, the
\emph{complement of} $G$ \emph{with respect to} $S$ is $GS:=Gs$.
\end{definition}

\section{Complementation and symmetry}

\begin{definition}
A (bicoloured) graph $G$ is \emph{invertible} if $V(G)$ is one of its
complementation sets. When no bicolouring is specified, $G$ is assumed to have its
natural bicolouring. The \emph{inverse} of an invertible graph $G$,
written $G^{-1}$, is $GS$, where $S=V(G)$.
\end{definition}

\begin{theorem}\label{automo}
Let $\Phi \in Aut(G)$ (the set of all automorphisms of the bicolored
graph $G$). If $\Phi$ stabilizes $S\subset V(G)$ (i.e. $\Phi(S)=S$),
then $\Phi \in Aut(GS)$.
\end{theorem}

\begin{proof}
By definition, $\Phi$ preserves colour and $[u,v]\in E(G)\iff [\Phi u, \Phi v]\in
E(G)$. By symmetry, $u$ changes colour from $G$ to $GS$ if and only if
$\Phi(u)$
changes colour from $G$ to $G\Phi (S)=GS$. Also, we have that $[u,v]\in E(G\triangle GS)\iff [\Phi u, \Phi v]\in
E(G\triangle G\Phi (S))=E(G\triangle GS)$. Thus $\Phi$ preserves
colour over $GS$ and $[u,v]\in E(GS)\iff
[\Phi u, \Phi v]\in E(GS)$.
\end{proof}

\begin{corollary} \label{cor_aut}
$Aut(G^{-1})=Aut(G).$
\end{corollary}

Here is a point to watch out for: the two groups $Aut(G)$ and
$Aut(G^{-1})$ are in fact identical, not just isomorphic.

\begin{corollary}
The inverse of an invertible vertex-transitive (bicoloured) graph is vertex-transitive.
\end{corollary}

We present an explicit formula for the inverse of a cycle:

\begin{proposition}
A cycle $C_n=\operatorname{Cay}(\mathbb{Z}_n,\{1,-1\})$ of length $n\geq 3$ is invertible
if and only if $n\not\equiv 0$(mod $3$). Furthermore, for $m\geq 1$,
$$C_{3m+1}^{-1}=\operatorname{Cay}(\mathbb{Z}_{3m+1},\{1,3,4,6,...,
3i-2,3i,...,3m-2,3m\}),$$
$$ C_{3m+2}^{-1}=\operatorname{Cay}(\mathbb{Z}_{3m+2},\{2,3,...,3i-1,3i,...,3m-1,3m\}).$$
\end{proposition}

\begin{figure}[htb!]
\begin{center}
\epsfig{file=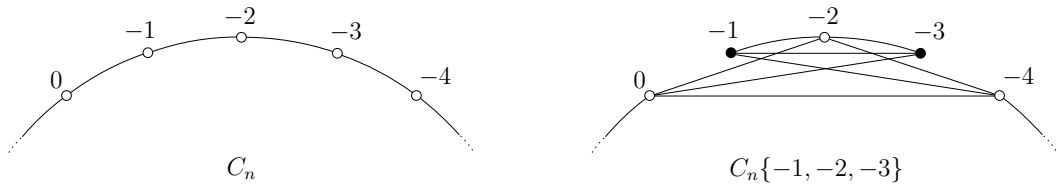,width=14cm}

\caption{\label{cycle}The complementation of $C_n$ with respect to $\{-1,-2,-3\}$.}
\end{center}
\end{figure}

\begin{proof}
$C_3$ is not invertible, while $C_4^{-1}=C_4=\operatorname{Cay}(\mathbb{Z}_4,\{1,3\})$ and
$C_5^{-1}=\operatorname{Cay}(\mathbb{Z}_5,\{2,3\})$. By Corollary~\ref{cor_aut}, the
inverse of an invertible circulant (a Cayley graph on $\mathbb{Z}_n$) is also a
circulant. For $n\geq 6$, removing vertices $-1$,$-2$ and $-3$ from $C_n\{-1,-2,-3\}$ yields $C_{n-3}$ (see Figure~\ref{cycle}). Thus, $C_n$ is
  invertible if and only if $C_{n-3}$ is, and given
  $C_{n-3}^{-1}=\operatorname{Cay}(\mathbb{Z}_{n-3},S)$, we must have $C_n^{-1}=\operatorname{Cay}(\mathbb{Z}_n,S')$
  with $S\subset S'$. The result follows by induction.
\end{proof}

\psfrag{0}[][][.8]{$0$}
\psfrag{n-1}[][][.8]{$-1$}
\psfrag{n-2}[][][.8]{$-2$}
\psfrag{n-3}[][][.8]{$-3$}
\psfrag{n-4}[][][.8]{$-4$}
\psfrag{Cn}[][][.8]{$C_n$}
\psfrag{Cn{n-1,n-2,n-3}}[][][.8]{$C_n\{-1,-2,-3\}$}

The remainder of this section is concerned with self-complementary symmetric (i.e. vertex- and edge-transitive) graphs. They were characterized in [11] and classified in [9].

From Theorem~\ref{automo}, we obtain that:

\begin{corollary}
An invertible self-complementary symmetric graph has either itself or its complement as an inverse.
\end{corollary}

\begin{definition}
Let $q=p^r$ for some odd prime $p$, with $q\equiv 1$(mod $4$). Let $\Gamma$ be the additive group of the finite field $F_q$ with $q$ elements.  Let $\omega$ be a primitive root in $F_q$ and $S=\{\omega ^2,\omega ^4,..., \omega ^{q-1}\}$, the set of non-zero squares. The \emph{Paley graph} of order $q$ is $\operatorname{Cay} (\Gamma, S)$.
\end{definition}

\begin{definition}
Let $q=p^r$ for some prime $p\equiv 3$(mod $4$) and some even $r$. Let $\Gamma$ and $\omega$ be as in the previous definition. Let $S=\{\omega^k|k\equiv 0, 1$(mod $4)\}$. The $\mathcal{P}^*$-\emph{graph} of order $q$ is $\operatorname{Cay} (\Gamma, S)$.
\end{definition}

Peisert [9] showed that, up to isomorphism, the
self-complementary symmetric graphs are the Paley and
$\mathcal{P}^*$-graphs and one additional graph $G(23^2)$ on $23^2$
vertices. Peisert gives a construction of $G(23^2)$ which allows us to interpret it as a Cayley graph $\operatorname{Cay} (\Gamma, S)$
with $\Gamma$ being the additive group of $F_{23^2}$ and such that
$\mathbb{Z}_{23}^*\subset S$.  

The following definition and theorem are taken from [3].

\begin{definition}
A graph is \emph{strongly regular modulo} $s$ \emph{with parameters} $(v,k,\lambda, \mu)$ if, modulo $s$, the number of vertices is congruent to $v$; the degree of each vertex, to $k$; the number of common neighbours of any two adjacent vertices, to $\lambda$; and the number of common neighbours of any two non-adjacent vertices, to $\mu$.
\end{definition}

\begin{theorem}[Fon-der-Flaass {[3]} Theorem 3.6]\label{thm_flaass}
$\overline{G}\in \mathcal{C}G$ if and only if $G$ is strongly regular modulo $2$ with parameters $(1,0,0,1)$.
\end{theorem}

\begin{proposition}\label{compl}
Let $G=\operatorname{Cay} (\Gamma, S)$ of order $n$ be a Paley graph, a $\mathcal{P}^*$-graph or the special graph $G(23^2)$. $\overline{G}\in \mathcal{C}G$ if and only if $2\notin S$.
\end{proposition}

\begin{proof}
The case $n=1$ is trivial, so let $n>1$. Since a self-complementary symmetric graph has order $n\equiv 1$(mod $4$) and is regular of degree $(n-1)/2$, using a simple counting argument and Theorem~\ref{thm_flaass}, the following are equivalent:\\

\noindent (1) $\overline{G}\in \mathcal{C}G$;\\
\noindent (2) the number of common neighbours of two adjacent vertices is even;\\
\noindent (3) the number of common neighbours of two non-adjacent vertices is odd.\\

We have $1,-1\in S$. Therefore $G$ admits the automorphism $\Phi$ defined by $\Phi (u)=-u$. The vertices $-1$ and $1$ are adjacent to $0$ and, because of $\Phi$, have an odd number of common neighbours. Since $1$ and $-1$ are adjacent if and only if $2=1-(-1)\in S$, the result follows.
\end{proof}

\begin{theorem}\label{compl2}
Given a self-complementary symmetric graph $G$ of order $n=p^r$, $\overline{G}\in \mathcal{C}G$ if and only if $n=1$ or $n\equiv 5$(mod $8$).
\end{theorem}

\begin{proof}
The idea is to check whether $2\in S$ and apply
Proposition~\ref{compl}. By construction, the statement is true for
$G(23^2)$. The remaining cases are nice exercices in number theory:

If $p\equiv 1$(mod $8$) (and therefore $G$ is Paley graph), we know
that $2$ is a quadratic residue modulo $p^r$ and so $2\in S$.

If $p\equiv 3$(mod $4$), we must have $r$ even and thus
$n-1=p^r-1\equiv 0$(mod $8$). Since $(p-1,8)=2$ (where $(a,b)$ is the
gcd of $a$ and $b$),

      $$n-1=(p-1)(p^{r-1}+p^{r-2}+...+p+1)=4(p-1)c$$
and 
     $$2^{\frac{n-1}{4}}=(2^{p-1})^c\equiv 1(\mbox{mod }p).$$ On the other
  hand, if $s=\omega^t$ for a primitive root $\omega$ in $F_{p^r}$, we
  have $$2^{\frac{n-1}{4}}=(\omega^{\frac{n-1}{4}})^t.$$ Therefore
  $t\equiv 0$(mod $4$) and $2\in S$.

We leave the cases $p\equiv 5$(mod $8$), $r$ even, and $p\equiv 5$(mod
$8$), $r$ odd, to the reader.
\end{proof}

As a consequence of Theorem~\ref{compl2}, we would have $G^{-1}=G$ for any invertible self-complementary graph $G$ of order $n\equiv 1$(mod $8$).

The following conjectures are suggested by computer testing:

\begin{conjecture} \label{conjA}
A self-complementary symmetric graph of order $n$ is invertible if and only if $n=1$ or $n\equiv 5$(mod $8$), in which case $G^{-1}=\overline{G}$.
\end{conjecture}

\begin{conjecture} \label{conjB}
Let $G$ be the Paley graph of prime order $p$ (i.e., $G=\operatorname{Cay}(\mathbb{Z}_p,S)$, where $S$ is the
set of quadratic residues mod $p$). Let $a,b\in \mathbb{Z}_p^*$
of orders $4$ and $k$, respectively, with $a-1$ being a quadratic
residue. Let $S_1=a\langle b\rangle \cup \langle b\rangle $ and $S_2=-a\langle b\rangle \cup \langle b\rangle $. Then at most
one of $S_1$ or $S_2$ is a complementation set of $G$. Furthermore, if
$p\equiv 5$(mod $16$), it cannot be $S_1$ and if $p\equiv 13$(mod $16$), it
cannot be $S_2$.
\end{conjecture}

\begin{conjecture}
Let $G$ be the Paley graph on $p=4q+1$ vertices
with $p,q$ prime. Then
$(-2)^q \langle 2^4\rangle \cup \langle 2^4\rangle $ is a complementation set of $G$.
\end{conjecture}

Conjecture~\ref{conjB} cannot be strengthened by saying that exactly
one of $S_1$ or $S_2$ is a complementation set. The first instances
of graphs where neither are complementation sets occur, when $p\equiv
5$(mod 16), at $p=37,$ 421, 661, 1381, 1621, 2789, 2917, 3061,
and when $p\equiv 13$(mod 16), at $p=2381,$ 3181, 5437, 5821.

\begin{proposition} Let $G=\operatorname{Cay}(\Gamma, S)$ be the Paley graph on $p^r\equiv 5$(mod $8$) vertices. Let $a\in F_{p^r}^*$ be of order $4$. Then $\{0\}\cup
  \langle a\rangle $ is a complementation set of $G$.
\end{proposition}

\begin{proof}
Since $a\notin S$, we find that $\{0\}\cup \langle a\rangle $ induces a
white pentagon in $G$, which is invertible.
\end{proof}

\begin{figure}[htb!]
\begin{center}
\epsfig{file=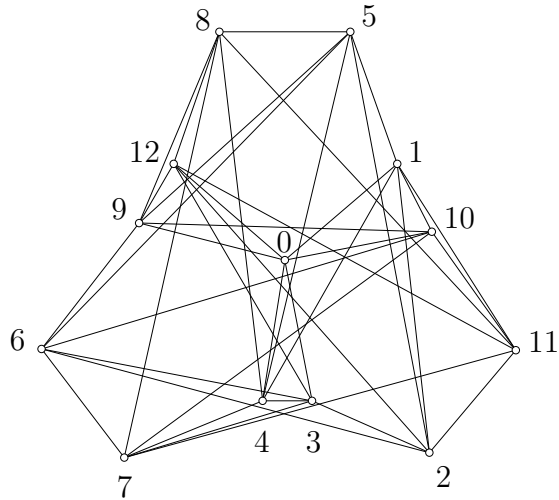,height=6.5cm}
\vspace{2mm}

\caption{\label{z13}The self-complementary symmetric graph on 13 vertices.} 
\end{center}
\end{figure}

In Figure~\ref{z13}, we can see that $S=\{0,1,-1,8,-8\}$, as well as 3$S$
and $-4S$ all induce pentagons. It can be verified that $G$ is
invertible by finding a reduced word $s$ satisfying $Gs=\overline{G}$ and
$\lambda(s)=13$. Can we find a method for constructing $s$ other than
the greedy algorithm (i.e. successively complement at any vertex you are
still allowed to) that can be generalized to other self-complementary
symmetric graphs? This seems to be a hard problem.

\psfrag{1}[][][1]{$1$}
\psfrag{2}[][][1]{$2$}
\psfrag{3}[][][1]{$3$}
\psfrag{4}[][][1]{$4$}
\psfrag{5}[][][1]{$5$}
\psfrag{6}[][][1]{$6$}
\psfrag{7}[][][1]{$7$}
\psfrag{8}[][][1]{$8$}
\psfrag{9}[][][1]{$9$}
\psfrag{10}[][][1]{$10$}
\psfrag{11}[][][1]{$11$}
\psfrag{12}[][][1]{$12$}
\psfrag{0}[][][1]{$0$}

%
%

We can construct many graphs
satisfying the conditions of Problem~\ref{prob1}. Given a
self-complementary symmetric graph $H$ such that $H^{-1}=\overline{H}=Ht$, with $t$ reduced, form $G$ by adding a new vertex $u$ to $H$, and joining $u$ to all vertices of $H$. Choosing $s=tu$ we have $Gs=G$ (by
Corollary~\ref{cor_aut}, $N_{Ht}(u)=N_H(u)$). Since there are an infinite
number of primes congruent to 5(mod 8) and assuming
Conjecture~\ref{conjA} holds, this family of graphs
would be infinite. It can be verified that if this
same $G$ is
provided with its natural bicolouring, then $G^{-1}=Gut^{-1}u\neq G$.

\begin{conjecture}
Given a bicoloured graph $G$ without isolated vertices or twins, $[G]$
is in bijection with the complementation sets of $G$.
\end{conjecture}

\section{References}

\setlength{\parindent}{0pt}

[1] Bouchet, A., Graphic presentations of isotropic systems,
\emph{J. Combin. Theory Ser. B} {\bf{45}} (1988), 58-76.\\

[2] Fon-der-Flaass, D. G., Distance between locally equivalent graphs (Russian),
\emph{Metody Diskret. Analiz.} {\bf{48}} (1989), 85-94, 106-107.\\ 

[3] Fon-der-Flaass, D. G., Local complementations of simple and oriented graphs (Russian),
\emph{Sibirsk. Zh. Issled. Oper.} {\bf{1}} (1994), 43-62, 87.\\

[4] Genest, F., Transition systems, orthogonality and local
complementation, submitted.\\

[5] Jackson, B., On circuit covers, circuit decompositions and Euler
tours of graphs, \emph{Surveys in combinatorics (Keele, 1993)}, London
Math. Soc. Lecture Note Ser. {\bf{187}}, Cambridge Univ. Press (1993),
191-210.\\

[6] Jaeger, F., A survey of the cycle double cover conjecture, \emph{Cycles in graphs (Burnaby, 1982)}, North-Holland Math. Stud. {\bf{115}}, North-Holland, Amsterdam (1985), 1-12.\\

[7] Kotzig, A., Moves without forbidden transitions in a graph,
\emph{Mat. \v Casopis Sloven. Akad. Vied} {\bf{18}} (1968), 76-80.\\

[8] Kotzig, A., \emph{Quelques remarques sur les transformations $\kappa$},
séminaire Paris (1977).\\

[9] Peisert, W., All self-complementary symmetric graphs,
\emph{J. Algebra} {\bf{240}} (2001), 209-229.\\

[10] Sabidussi, G., \emph{Eulerian walks and local complementation}, D.M.S. 84-21, Dép. de math. et stat., Universit\'e de Montr\'eal (1984).\\

[11] Zhang, H., Self-complementary symmetric graphs, \emph{J. Graph
  Theory} {\bf{16}} (1992), 1-5.\\


%% file: transition.tex







\chapter[]{Transition systems, orthogonality and local complementation.}

\begin{center}
{\bf\Large Fran\c cois Genest}
\end{center}

\section{Abstract}
Intimately related to the Cycle Double Cover
  Conjecture is the problem of finding a cycle decomposition
  orthogonal to a given transition system  in an eulerian graph. One
  approach consists in finding a black anticlique in the corresponding parity class
  of bicoloured graphs.

\section{Introduction}
This article discusses problems that seem foreign to each other yet prove to be inextricably linked.

Much effort has been spent on the Cycle Double Cover Conjecture; understandably so, given its many implications (see the surveys by Jackson [10] and Jaeger [11]).

Typical questions about transition systems of eulerian graphs involve
finding transition systems that are orthogonal to one
another. Fleischner [5][6] showed how the Dominating Circuit
Conjecture and Sabidussi's Orthogonality Conjecture imply the Cycle Double Cover Conjecture.

Kotzig [13] demonstrated how $\kappa$-transformations of Euler tours in $4$-regular connected graphs are essentially the same as local complementations. In a similar way, Sabidussi [14] made the connection between eulerian graphs with one specified Euler tour and parity classes of simple graphs.

We will go from transition systems to parity classes and back again, showing how the different problems are meshed together. The notion of pure graphs is explored; it is proposed that a characterization of primitive pure graphs would help to settle the Cycle Double Cover Conjecture.

\section{Preliminaries}
The terminology in this article
mostly follows Fleischner [6] and Jackson [10]. A \emph{cycle} is a non-empty 2-regular connected graph or subgraph. A \emph{tour} of a graph is a sequence
$\alpha=u_0e_1u_1e_2u_2...u_{n-1}e_n$ where $e_i$ is an edge incident
with the vertices $u_{i-1}$ and $u_i$ (where subscripts are read modulo
$n$) and where the edges are distinct. Tours are deemed equivalent if one can be
transformed into another by successive reversals and cyclic permutations of the sequence. We work
on the one hand with eulerian graphs, which may have multiple edges
and loops; such a graph will be denoted by $\Gamma$. On the other
hand, we will carry out local complementations of simple graphs which
will be denoted by $F,G$ or $H$.

A \emph{transition} at a vertex $u$ is either a pair $\{u,e\}$, where $e$ is a loop incident with $u$, or a set $\{u,e_1,e_2\}$, where $e_1$ and $e_2$ are distinct edges
incident with $u$. A
\emph{transition system} (TS), intuitively, is a partition of the ``half-edges''
of the graph into pairs of adjacent ones. To avoid the difficulties of dealing
with loops, we define transition systems by way of tour decompositions.

A \emph{tour decomposition} of an eulerian graph $\Gamma$ is a set of edge-disjoint tours of
$\Gamma$ whose union exhausts all edges. Note that a tour $u_0e_1u_1e_2u_2...e_n$
induces the transitions $\{u_i,e_i,e_{i\mbox{+}1}\}$, $i=0,...,n-1$. Given a tour
decomposition of $\Gamma$, the associated \emph{transition system} is the union of the
sets of transitions induced by the tours. Since we can reconstruct the tours from the transition system,
this is a 1-1 correspondence. Two transition systems (and by extension
their corresponding tour decompositions) are \emph{orthogonal} if they are
disjoint (we prefer ``orthogonal'' to the rather vague ``compatible''). Of course, such a comparison makes sense only when there are
no vertices of degree two. Since isolated vertices are irrelevant to
questions of orthogonality, we will only consider graphs of minimum
degree $\delta > 2$. In this setting, a transition system completely determines the
graph and we can talk about tours and cycles of transition systems without
ambiguity.

Settling a question raised by Nash-Williams, Kotzig [12] showed that for any transition system without loops there exists an orthogonal
Euler tour. When loops are present, a necessary and sufficient
condition for the existence of an orthogonal Euler tour is for all
loop transitions (transitions of the form $\{u,e\}$) incident with
vertices of degree 4 to be in the transition system. In [9], Jackson characterizes graphs having three pairwise
orthogonal Euler tours.

A natural question arises: when can we find
a cycle decomposition orthogonal to a given transition system? A trivial necessary
condition is that the transition system must not contain transitions of the form
$\{u,e_1,e_2\}$, where $u$ is of degree two within
the block (maximal 2-connected subgraph) containing $\{e_1,e_2\}$ (in particular, loop transitions
are excluded), because any cycle decomposition of the graph must contain these. A transition system
satisfying this condition is said to be \emph{admissible}. Many admissible
transition systems admit orthogonal cycle decompositions. A notable exception is the
transition system of $K_5$ decomposed into two 5-cycles. Sabidussi conjectures
that if the transition system is induced by an Euler tour (in which case it is
necessarily admissible), there always exists an orthogonal cycle
decomposition. Much of what follows uses or is inspired by his work on
the subject (see [14]). The results on complementation used in this article can be found in [8].

\section{Transition graphs}

\begin {figure} [htb!]
\begin {center}
\psfrag{u}[][][1]{$u$}
\psfrag{v}[][][1]{$v$}
\psfrag{w}[][][1]{$w$}
\psfrag{e1}[][][1]{$e_1$}
\psfrag{e2}[][][1]{$e_2$}
\psfrag{e3}[][][1]{$e_3$}
\psfrag{e4}[][][1]{$e_4$}
\psfrag{e5}[][][1]{$e_5$}
\psfrag{e6}[][][1]{$e_6$}
\psfrag{e7}[][][1]{$e_7$}
\psfrag{e8}[][][1]{$e_8$}
\epsfig {file=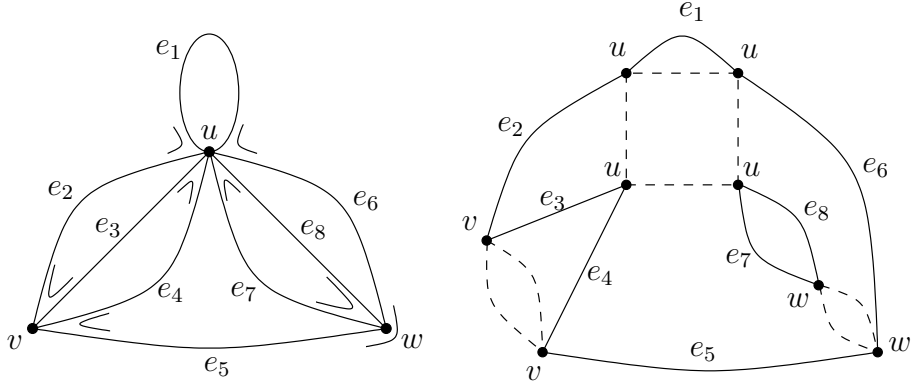}
\end{center}
\caption{\label{TS_TG}To the left, a graph with TS indicated by arcs. To the right, one of its TG's.}
\end {figure}

\begin{definition} Let $\Gamma$ be an eulerian graph of minimum degree
  $\delta >2$
  with transition system $S=\{t_{u,i}|u\in V(\Gamma ), i\in
  \{1,...,\frac{1}{2}d_\Gamma (u)\}$, where each $t_{u,i}$ contains $u$. A \emph{transition graph} (TG) of $S$ is a
  graph $\Gamma_S$ satisfying $V(\Gamma_S)=S$, $E(\Gamma_S)=E(\Gamma )\cup E'$, where $E'=\{e_{u,i}'|u\in V(\Gamma ),
i=1,...,\frac{1}{2}d_\Gamma (u)\}$, with incidences defined as follows:\\

$e\in E(\Gamma )$ is a non-loop edge incident with $t_{u,i}$ and $t_{v,j}$ if $t_{u,i}$ and $t_{v,j}$ are distinct and $e\in
t_{u,i}\cap t_{v,j}$; $e$ is a loop incident with $t_{u,i}$ if
$t_{u,i}=\{u,e\}$,\\

$e_{u,i}'$ is incident with $t_{u,i}$ and $t_{u,i\mbox{+}1}$,
$i=1,...,\frac{1}{2}d_\Gamma (u)-1$; $e_{u,d_\Gamma (u)/2}'$ is incident
with $t_{u,1}$ and $t_{u,d_G(u)/2}$.
\end{definition}
In Figure~\ref{TS_TG}, edges of the original graph are drawn as
solid lines while the edges in $E'$ are broken. To identify a
transition in the transition graph, simply consider the vertex label together
with the incident solid edges.

Note that, in general, different orderings of the transitions of $\Gamma $ will give
rise to different incidences for the edges in $E'$. In
fact, the transition graph will be unique (up to isomorphism) if and only if all vertices of $\Gamma $ are of degree 4 or
6. Observe also that a transition graph of $S$ has a natural transition system provided by the 2-factors
induced by $E(\Gamma )$ and $E'$. The idea behind constructing a transition graph is that
if we can find a cycle decomposition orthogonal to this new transition system, then,
contracting the edges in $E'$, we get a cycle decomposition orthogonal
to $S$. Since all loop transitions of $S$ are in the transition system
of $\Gamma_S$, Kotzig's result tells us we can find an Euler tour
orthogonal to the transition system of $\Gamma_S$. Notice that this tour will have successive
edges alternating between the two 2-factors (solid and broken edges). Accordingly, it
is called an \emph{alternating tour}. When writing it, we omit the
edges and, since $\Gamma_S$ is 4-regular, we are left with a \emph{double occurrence
word} (each letter appears exactly twice). At this point, we can construct an alternance (or circle)
graph. Alternance graphs have been studied in [3][7][15], for example.

\begin{definition}
Given a double occurrence word $w=a_1...a_{2n}$ made up of the letters
$u_1,...,u_n$, the \emph{alternance graph} of $w$ is the simple graph with vertex set $\{u_1,...,u_n\}$, and with an edge between $u_i$ and
$u_j$ if they alternate in $w$ (i.e. if $a_k=a_l=u_i$ and
$a_r=a_s=u_j$ with $k<l,r<s$ then either $k<r<l<s$ or $r<k<s<l$).
\end{definition}

\begin{definition}
An \emph{anticlique} $A$ of a graph $G$ is a maximal independent set of
vertices (i.e. no two vertices of $A$ are adjacent while every vertex in
$V(G)\setminus A$ has a neighbour in $A$).
\end{definition}

\begin{theorem}[Sabidussi {[14]} Theorem 5.1]
Let $\Gamma$ be a transition graph with $E(\Gamma)=E\cup E'$ being the partition into the two
$2$-factors. Let $w$ be (the double occurrence word of) an alternating
tour of $\Gamma$. Every anticlique $A$ of the alternance graph of $w$
determines a cycle decomposition of $\Gamma$ into $|A|+1$ cycles. Moreover,
if $A$ is odd (i.e. consists of vertices of odd degree), the
decomposition is alternating.
\end{theorem}

\begin{proof}
Since any cyclic permutation of the letters of $w$ gives rise to the
same alternance graph, we can suppose that
$w=ua_1a_2...a_rua_{r+1}...a_{2n-2}$ with $u\in A$. We can further
suppose that $a_1,...,a_r \notin A$. Since $A$ is a covering, we must
have that $a_1,...,a_r$ are all distinct. Split the tour $w$ into the
cycle $ua_1...a_ru$ and the tour $ua_{r+1}...a_{2n-2}$. If $r$ is
odd, then both the cycle and the new tour have an even number of edges
and remain alternating. While $u$ is not covered by $A'=A\setminus
\{u\}$, it only appears once in the new tour and (replacing $A$ with
$A'$) the result follows by induction on $|A|$.
\end{proof}

Given an arbitrary alternating tour, the corresponding alternance graph may or may not admit an odd anticlique. However, Sabidussi [14], extending Kotzig's work on
$\kappa$-transformations [12], showed that any other alternating
Euler tour can be obtained by a sequence of twists and switches.

\begin{definition}
The $u$-\emph{twist} (or \emph{twist} at $u$) of the Euler tour
$$w=ua_1...a_rua_{r+1}...a_s$$ is the
Euler tour $$w'=ua_r...a_1ua_{r+1}...a_{s}.$$
\end{definition}

Note that if $w$ is alternating and $r$ is even (i.e. $u$ is of even
 degree in the alternance graph), then $w'$ is also alternating.

\begin{definition}
Given vertices $u$ and $v$ that alternate in the Euler tour
$$w=ua_1...a_kva_{k+1}...a_lua_{l+1}...a_{r}va_{r+1}...a_s,$$ the
$uv$-\emph{switch} of $w$ is the Euler tour
$$w'=ua_1...a_kva_{r+1}...a_sua_{l+1}...a_rva_{k+1}...a_l$$
\end{definition}

If $w$ is alternating and $l$ and $r-k$ are odd (i.e. $u,v$ are of odd
degree in the alternance graph) then $w'$ is also alternating. Also, a
$uv$-switch has the same effect as successive twists at $u,v$ and $u$. Looking at the alternance graph, the result of a $u$-twist is a local
complementation with respect to $u$.\\

\section{Local complementation and pure graphs}

In the following, by a \emph{bicoloured} graph is meant a simple graph
with a $\{$black, white$\}$-colouring of its vertices. A closer look at local complementation and related results can be found in [8].  

\begin{definition}
The (\emph{local}) \emph{complement} of a bicoloured graph $G$ with respect to a
vertex $u$ is a bicoloured graph $Gu$ such that
 $$V(Gu)=V(G)$$
$$E(Gu)=E(G)\triangle E(K_{N_G(u)})\mbox{ (symmetric difference)}$$
where $K_{N_G(u)}$ is the complete
graph on the neighbourhood $N_G(u)$. That is to say, the adjacency
relation of $Gu$ coincides with that of $G$ except on $N_G(u)$, where
it is replaced by its complement. The bicolouring of $Gu$ is defined
to be the same as that of $G$ if $u$ is black in $G$; if $u$ is white
in $G$, then it is obtained from the bicolouring of $G$ by reversing the colours of
the vertices in $N_G(u)$.
\end{definition}

Complementation is extended recursively to
words in the alphabet $V(G)$ by putting $Gsu=(Gs)u$, where $s\in
V(G)^{*}$ (= the set of all words on $V(G)$) and $u\in V(G)$.

The symbol $[uv]$ is short for $uvu$.\\

For the purposes of the next definition, call a set a words $W\subset
V(G)^*$ $parity$ $closed$ if\\

\noindent (i) $W$ contains the empty word;

\noindent (ii) if
$s\in W$ and $u$ is a white vertex of $Gs$, then $su\in W$;

\noindent (iii) if
$s\in W$ and $u,v$ are adjacent black vertices of $Gs$, then
$s[uv],s[vu]\in W$.\\

Clearly the intersection of parity closed sets is
parity closed, hence there is a smallest parity closed set, denoted by
$W^{\circ}(G)$.

\begin{definition}
The \emph{parity class} of a bicoloured graph $G$ is $$[G]=\{Gs|s\in W^{\circ}(G)\}.$$
\end{definition}

\begin{figure}[htb!]
\begin{center}
\epsfig {file=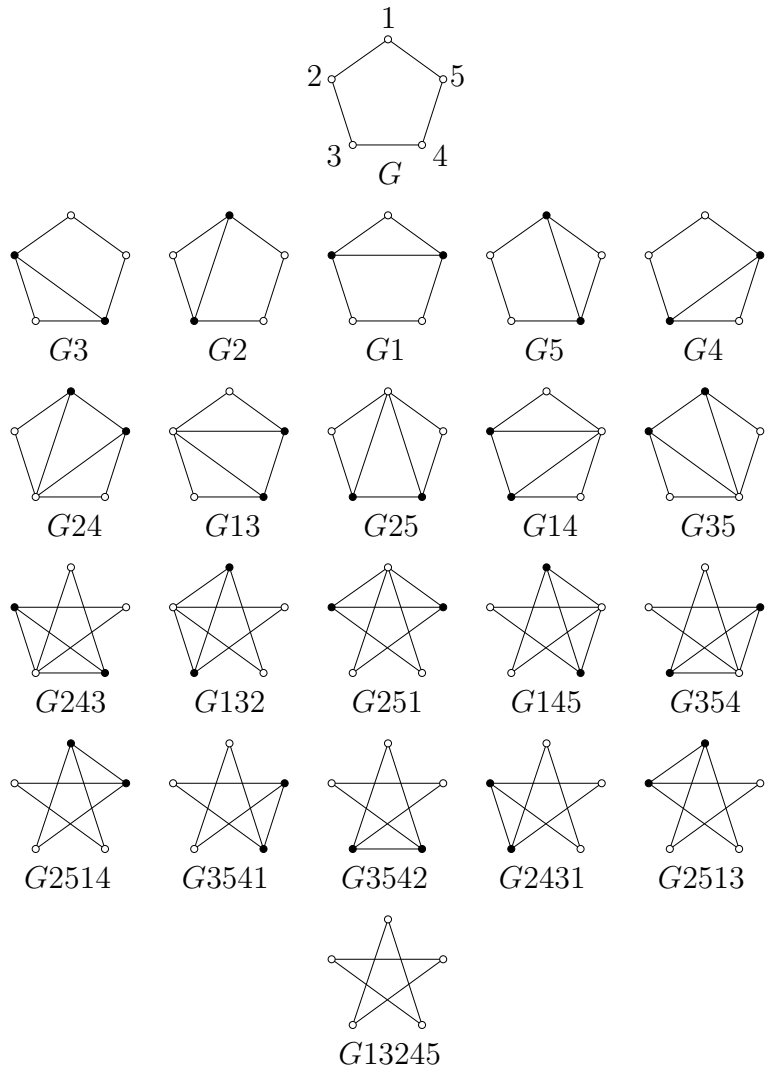, height=14cm}
\end {center}
\caption{\label{classe_c5}The parity class of the pentagon (with its natural colouring).}
\end{figure}

An example of a parity class is given in Figure~\ref{classe_c5}. There, $G$ is a
cycle of length 5 on the vertex set $\{1,2,3,4,5\}$ with all vertices
coloured white. The other graphs in $[G]$ are identified by words
in $W^{\circ}(G)$. These words are not unique: for example, $G13=G31$.\\

\psfrag{G}[][][1]{$G$}
\psfrag{1}[][][1]{$1$}
\psfrag{2}[][][1]{$2$}
\psfrag{3}[][][1]{$3$}
\psfrag{4}[][][1]{$4$}
\psfrag{5}[][][1]{$5$}
\psfrag{G1}[][][1]{$G1$}
\psfrag{G2}[][][1]{$G2$}
\psfrag{G3}[][][1]{$G3$}
\psfrag{G4}[][][1]{$G4$}
\psfrag{G5}[][][1]{$G5$}
\psfrag{G13}[][][1]{$G13$}
\psfrag{G14}[][][1]{$G14$}
\psfrag{G24}[][][1]{$G24$}
\psfrag{G25}[][][1]{$G25$}
\psfrag{G35}[][][1]{$G35$}
\psfrag{G243}[][][1]{$G243$}
\psfrag{G132}[][][1]{$G132$}
\psfrag{G251}[][][1]{$G251$}
\psfrag{G145}[][][1]{$G145$}
\psfrag{G354}[][][1]{$G354$}
\psfrag{G2514}[][][1]{$G2514$}
\psfrag{G3541}[][][1]{$G3541$}
\psfrag{G3542}[][][1]{$G3542$}
\psfrag{G2431}[][][1]{$G2431$}
\psfrag{G2513}[][][1]{$G2513$}
\psfrag{G13245}[][][1]{$G13245$}

If a simple graph is given its \emph{natural} colouring, where vertices of
even degree are white and vertices of odd degree are black, we can see
that the colouring of each graph in its parity class is also natural. We now have
that to each transition graph there corresponds a unique parity class, and our
search for an alternating cycle decomposition amounts to finding a
black anticlique in some graph of the parity class. 

\begin{definition}
A parity class is \emph{pure} if it
contains no graph with an anticlique of black vertices. By extension,
its member graphs are also said to be \emph{pure}.
\end{definition}

Up to local complementation, the two smallest connected pure graphs
are the one-point graph and the pentagon with their natural colourings.

\begin{definition}
A \emph{split} in a graph $G$ is a partition of its vertex set into
$V=V_1\cup V_2$, with at least two vertices in each $V_i$, such that all
vertices in $N(V_1)\cap V_2$ are adjacent to all vertices in
$N(V_2)\cap V_1$. A splitless graph is said to be \emph{prime}.
\end{definition}

It can easily be checked that a split remains a split after local
complementations (see [1]), so this is a feature of the whole parity class.

\begin{definition}
A \emph{rooted} bicoloured graph with root $z$ is a couple $(G,z)$, where
$z\in V(G)$, $G$ being a simple graph with a bicolouring defined on $V(G)\setminus \{z\}$.
\end{definition}

Whenever we mention rooted graphs in this paper, it will be understood
that we are talking about rooted bicoloured graphs.

Note that from a split $V=V_1\cup V_2$ of some bicoloured graph $G$,
we can construct two rooted graphs as follows: let $G_1$
(respectively $G_2$)
be obtained from $G$ by
identifying the vertices of $V_2$ (respectively $V_1$) to a single
new vertex $z_1$ (respectively $z_2$) which will be the root vertex of $G_1$
(respectively $G_2$) and removing loops. $G_1$ and $G_2$ are said to be \emph{induced}
by the split.

\psfrag{u1}[][][1]{$u_1$}
\psfrag{u2}[][][1]{$u_2$}
\psfrag{u3}[][][1]{$u_3$}
\psfrag{u4}[][][1]{$u_4$}
\psfrag{u5}[][][1]{$u_5$}
\psfrag{u6}[][][1]{$u_6$}
\psfrag{u7}[][][1]{$u_7$}
\psfrag{u8}[][][1]{$u_8$}
\psfrag{u9}[][][1]{$u_9$}
\psfrag{m1}[][][1]{$z_1$}
\psfrag{m2}[][][1]{$z_2$}
\psfrag{G}[][][1]{$G$}
\psfrag{G'}[][][1]{$G'$}
\psfrag{G''}[][][1]{$G''$}
\psfrag{G1}[][][1]{$G_1$}
\psfrag{G2}[][][1]{$G_2$}
\psfrag{G1'}[][][1]{$G_1'$}
\psfrag{G2'}[][][1]{$G_2'$}
\psfrag{G1''}[][][1]{$G_1''$}
\psfrag{G2''}[][][1]{$G_2''$}
\psfrag{rw(G1)}[][][1]{$rw(G_1)$}
\psfrag{lw(G2)}[][][1]{$lw(G_2)$}
\psfrag{rw(G2)}[][][1]{$rw(G_2)$}
\psfrag{lw(G1)}[][][1]{$lw(G_1)$}
\psfrag{rb(G1')}[][][1]{$rb(G_1')$}
\psfrag{lb(G2')}[][][1]{$lb(G_2')$}
\psfrag{rc(G1'')}[][][1]{$rc(G_1'')$}
\psfrag{lcu9(G2'')}[][][1]{$lc_{u_9}(G_2'')$}

Let $G$ be a bicoloured graph with a split $V=V_1\cup V_2$
inducing the rooted graphs $G_1$ and $G_2$. We will construct
six pairs of parity classes from $G_1$ and $G_2$ and
show that if any one pair consists of two pure parity classes then
$[G]$ is also a pure class. Because these constructions and proof are technical and hardly readable by themselves, we will, after giving the necessary definitions, discuss them by means of an example that illustrates all possibilities that may arise.

\begin{definition}
Given a vertex $v$ adjacent to a black vertex $u$ of a bicoloured
graph $G$, the \emph{complementation} of $G$ with respect to the set
$\{u,v\}$ is

$$G\{u,v\}=
\begin{cases}
  Gvu& \text{if $v$ is white},\\
  G{[vu]}& \text{if $v$ is black}.
\end{cases}$$

\end{definition}

This is well defined because in the case where $v$ is black,
$G[vu]=G[uv]$. Complementation with respect to a couple is a special
case of set complementation (see [8]). 

\begin{definition} \label{root}
Given a rooted graph $G$ with root $z$, we define the following bicoloured graphs:\\

$rw(G)$ is obtained by colouring $z$ white in $G$,

$rb(G)$ is obtained by colouring $z$ black in $G$,

$rc(G)$ is obtained by complementing $G$ at $z$ as if $z$ were white and colouring $z$ black,

$lw(G)=G-z$,

$lb(G)$ is obtained by complementing $G$ at $z$ as if $z$ were white and then removing $z$,

and $lc_v(G)$, given $v$ adjacent to $z$, is obtained by colouring $z$ black and letting $lc_v(G)=G\{z,v\}-z$.\\

\noindent $rw(G),rb(G),rc(G)$ are called \emph{root} graphs of $G$ and $lw(G),lb(G),lc_v(G)$ are \emph{leaf} graphs of $G$.
\end{definition}

Consider the following bicoloured graph:

\begin{center}
\epsfig{file=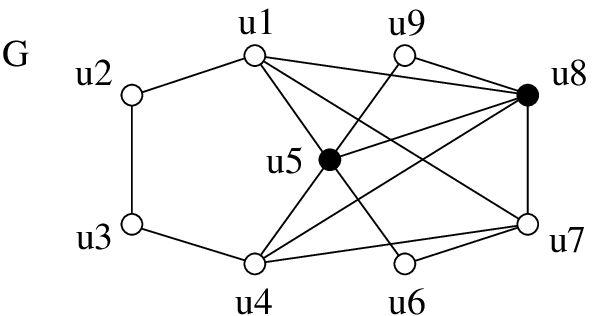,height=3cm}
%
\end{center}

$G$ has only one split (i.e. $V_1=\{u_1,u_2,u_3,u_4\}$, $V_2=\{u_5,u_6,u_7,u_8,u_9\}$), which induces:

\begin{center}
\epsfig{file=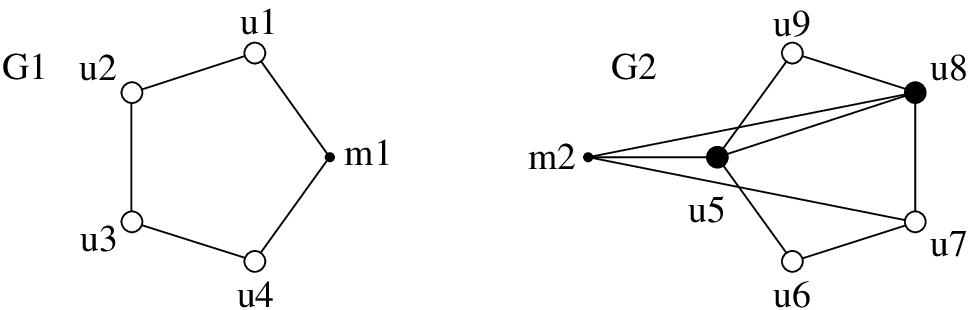,height=3cm}

%
\end{center}

From these rooted graphs, we construct the bicoloured graphs
$rw(G_1)$, $lw(G_2)$, $lw(G_1)$ and $rw(G_2)$:

\begin{center}
\epsfig{file=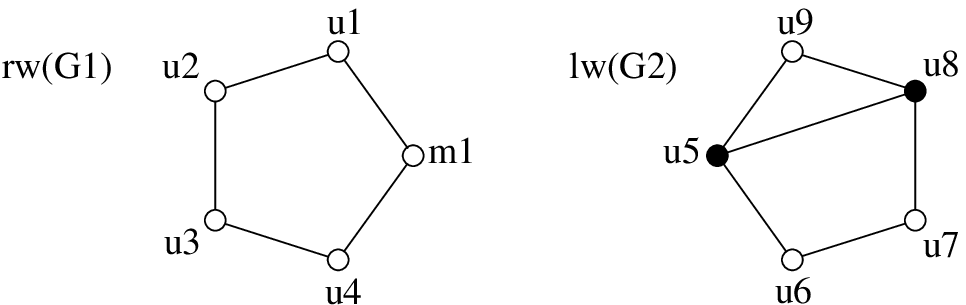,height=3cm}
\end{center}
\begin{center}
\epsfig{file=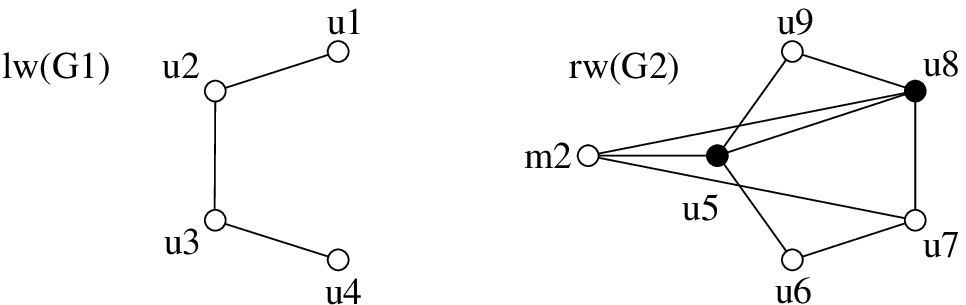,height=3cm}
%
\end{center}

We see that $rw(G_1)$ and $lw(G_2)$ form a pair of pure root and leaf graphs. However, consider now what happens when we start with $G'=Gu_1$:

\begin{center}
\epsfig{file=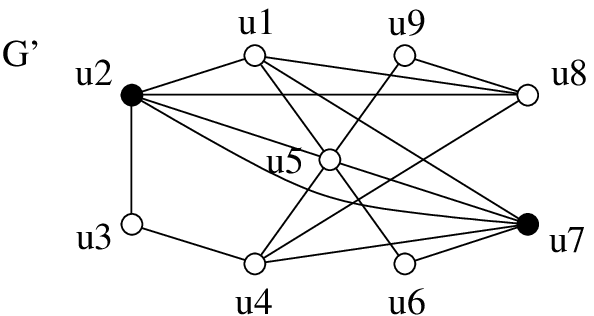,height=3cm}
%
\end{center}

In this case, the split induces:

\begin{center}
\epsfig{file=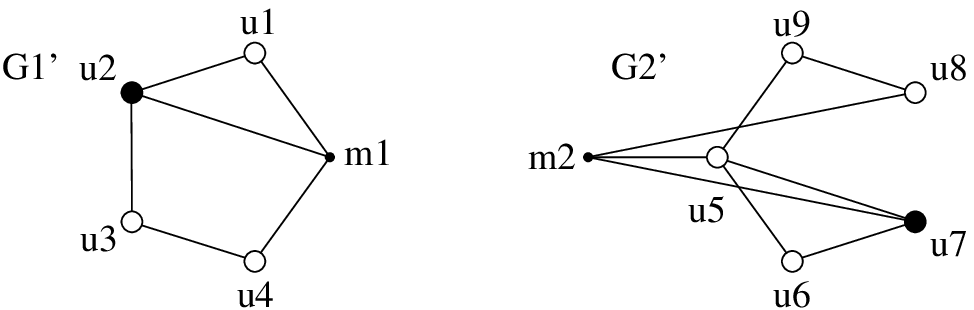,height=3cm}
%
\end{center}

and none of $rw(G_1'),lw(G_2'),rw(G_2'),lw(G_1')$ are pure. To find
pure graphs again, we construct an $\{rb(G_i'), lb(G_j')\}$ pair, in this
case the root graph $rb(G_1')$ and the leaf graph $lb(G_2')$:

\begin{center}
\epsfig{file=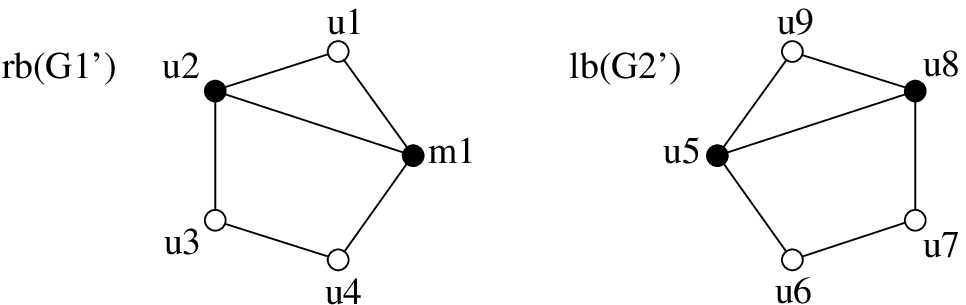,height=3cm}
%
\end{center}

Another situation presents itself when we start with $G''=Gu_1u_5$:

\begin{center}
\epsfig{file=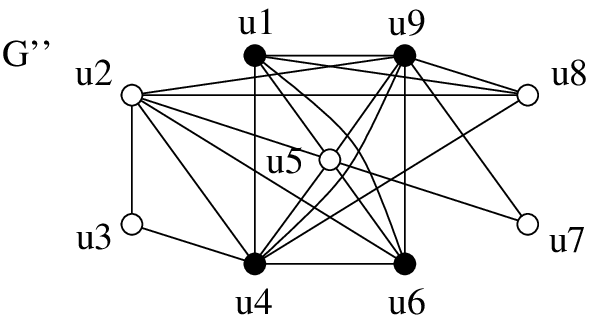,height=3cm}
%
\end{center}

Now the split induces:

\begin{center}
\epsfig{file=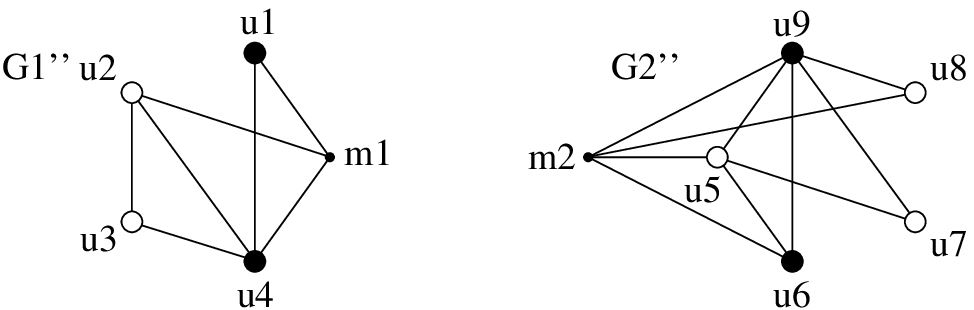,height=3cm}
%
\end{center}

and none of $rw(G_i''),rl(G_j''),rb(G_i''),lb(G_j'')$ are pure, for
$\{i,j\}=\{1,2\}$. We have to resort to the last type of construction:
for example, the root graph $lc(G_1'')$ and the leaf graph
$lc_{u_9}(G_2'')$ form a pure pair:

\begin{center}
\epsfig{file=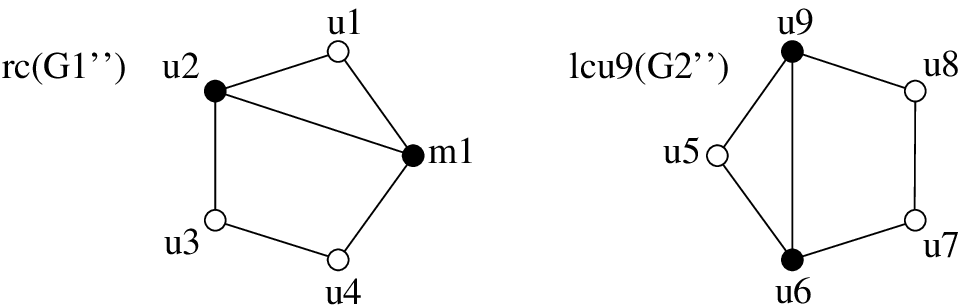,height=3cm}
%
\end{center}

While these constructions are unappealing, the good news is that
we're done: any other graph $H\in[G]$ separated into rooted graphs
$H_1$ and $H_2$ (corresponding to $V_1$ and $V_2$ respectively) by the split will yield a pair of pure graphs of the
form $\{rw(H_1),lw(H_2)\}$, $\{rb(H_1),lb(H_2)\}$ or
$\{rc(H_1),lc_v(H_2)\}$, as will be seen below.

\begin{definition} \label{decomp}
Let $G$ be a simple bicoloured graph with a split $V=V_1\cup V_2$
inducing the rooted graphs $(G_1,z_1)$ and $(G_2,z_2)$. An \emph{essential
decomposition} of $G$ along this split is a pair of root and leaf
graphs of the form $\{rw(G_i),lw(G_j)\}$, $\{rb(G_i),lb(G_j)\}$ or $\{rc(G_i),lc_v(G_j)\}$, where $\{i,j\}=\{1,2\}$, and $v$ is adjacent to $z_j$.
\end{definition}

\begin{lemma}
Let $(G,z)$ be a rooted graph. If $v,w\in N_G(z)$, then $lc_v(G)$ and $lc_w(G)$ are in the same parity class (i.e. $lc_w(G)\in [lc_v(G)]$).
\end{lemma}

\begin{proof}
Let $H$ be the bicoloured graph obtained from $G$ by colouring $z$
black. For $s\in W^\circ(lc_v(G))\subset W^\circ(H\{v,z\})$, we have
$lc_v(G)s=H\{v,z\}s-z$. Using local substitution rules (see [8]), consider the following three cases:\\

if $v\neq w$ are both white in $G$, then $[vw]\in W^{\circ}(lc_v(G))$,
and $wz\sim_H vvwz\sim_H vz[vw]$, so that
$lc_w(G)=Hwz-z=Hvz[vz]-z=lc_v(G)[vz]\in [lc_v(G)]$;\\

if $v$ and $w$ are of different colours, say $v$ is black and $w$ is white, then $wv\in W^{\circ}(lc_v(G))$ and $wz\sim_H wzvv\sim_H [vz]wv$;\\

if $v\neq w$ are both black, then $[vw]\in W^{\circ}(lc_v(G))$ and $[wz]\sim_H [vz][vw]$.
\end{proof}

\begin{definition}
Given a bicoloured graph $G=(V,E)$ with a split $V=V_1\cup V_2$
inducing the rooted graphs $(G_1,z_1)$ and $(G_2,z_2)$, the \emph{essential decompositions} of the parity
class $[G]$ along this split with respect to $G$ are the pairs
$\{[rw(G_i)],[lw(G_j)]\},\{[rb(G_i)],[lb(G_j)]\}$ and, if there is an
edge between $V_1$ and $V_2$ in $G$, $\{[rc(G_i)],[lc_v(G_j)]\}$,
where $\{i,j\}=\{1,2\}$. For convenience, if there is no edge between $V_1$ and $V_2$ in $G$, the pair $\{[\emptyset],[\emptyset]\}$, where $\emptyset$ stands for the empty graph, is also called an essential decomposition of $[G]$.
\end{definition}

Let us write $Rw_G$,$Lw_G$,$Rb_G$,$Lb_G$ for $[rw(G_1)]$,$[lw(G_2)]$,$[rb(G_1)]$,$[lb(G_2)]$ respectively. Let $Rc_G=[rc(G_1)]$ and $Lc_G=[lc_v(G_2)]$ if there is an edge between $V_1$ and $V_2$ in $G$, or else set them equal to the trivial class $[\emptyset]$. To indicate the essential decompositions obtained by interchanging $V_1$ and $V_2$, write $Rw_G',Lw_G',Rb_G',Lb_G',Rc_G'$ and $Lc_G'$. The proof of the following is left to the reader:

\begin{lemma} \label{recolour}
Let $G_1,G_2$ be identical bicoloured graphs except for a vertex $u$
which is of different colour in $G_1$ and $G_2$. Let $v$ be a black vertex
adjacent to $u$. Consider $H_1=G_1\{u,v\}$ and $H_2=G_2\{u,v\}$. Then $u$ is
black in $H_1$ and $H_2$, and each of $H_1,H_2$ can be obtained from the other by complementing at
$u$ and reversing the colouring on $N_{H_1}(u)$ (i.e. by complementing
as if $u$ were white). Conversely, if $H_1$ and $H_2$
are as described, then for $v\in N_G(u)$, $H_1\{u,v\}$ and $H_2\{u,v\}$
differ only by the colour of $u$.
\end{lemma}

\begin{theorem} \label{invariant}
Let $G$ be a bicoloured graph with a split $V=V_1\cup
V_2$. The essential decompositions of $[G]$ along this split with respect to $H\in [G]$ are independent of the choice of $H$. To be more precise, the pairs $\{Rw_H,Lw_H\},\linebreak[5] \{Rb_H,Lb_H\},\{Rc_H,Lc_H\}$ are identical to $\{Rw_G,Lw_G\},\{Rb_G,Lb_G\},\{Rc_G,Lc_G\}$ up to permutation, and the pairs $\{Rw_H',Lw_H'\},\{Rb_H',Lb_H'\},\{Rc_H',Lc_H'\}$ are identical to $\{Rw_G',Lw_G'\},\{Rb_G',Lb_G'\},\{Rc_G',Lc_G'\}$ up to permutation.
\end{theorem}

\begin{proof}
When there is no edge between $V_1$ and $V_2$ in $G$, the result is
trivial, so suppose $V_1$ and $V_2$ are joined by an edge. We have that $H=Gs$ for some $s\in W^\circ (G)$. By induction,
it suffices to consider the cases $s=u$, where $u$ is a white vertex of $G$, and
$s=[xy]$, where $x,y$ are adjacent black vertices of $G$. Let the split induce the rooted graphs
$(G_1,z_1),(G_2,z_2)$ from $G$, and the rooted graphs $(H_1,z_1),(H_2,z_2)$ from $H$. Let $v\in N_{G_2}(z_2)$ and $w\in N_{H_2}(z_2)$ so that $Lc_G=[lc_v(G_2)],Lc_H=[lc_w(H_2)]$. We consider the following eleven cases.\\

Case 1) white $u\in V_1\setminus N_{G_1}(z_1)$.\\
$rw(H_1)=rw(G_1)u$, $rb(H_1)=rb(G_1)u$, $rc(H_1)=rc(G_1)u$, and since $H_2=G_2$, $lw(H_2)=lw(G_2)$, $lb(H_2)=lb(G_2)$, and (choosing
without loss of generality
$w=v$) $lc_w(H_2)=lc_v(G_2)$.\\

Case 2) white $u\in N_{G_1}(z_1)$.\\
$rw(H_1)=rb(G_1)u$, $lw(H_2)=lb(G_2)$, $rb(H_1)=rw(G_1)u$,
$lb(H_2)=lw(G_2)$. Lemma~\ref{recolour} gives $rc(H_1)=rc(G_1)[uz_1]$ and, choosing
$w=v$, $lc_w(H_2)=lc_v(G_2)$.\\

Case 3) white $u\in N_{G_2}(z_2)$.\\
$rw(H_1)=rw(G_1)z_1$, $lw(H_2)=lw(G_2)u$, $rb(H_1)=rc(G_1)$,
$rc(H_1)=rb(G_1)$. Choose $w=v=u$, so that $lc_w(H_2)=lb(G_2)$ and $lb(H_2)=lc_v(G_2)$.\\

Case 4) white $u\in V_2\setminus N_{G_2}(z_2)$.\\
$rw(H_1)=rw(G_1)$, $lw(H_2)=lw(G_2)u$, $rb(H_1)=rb(G_1)$, $lb(H_2)=lb(G_2)u$, $rc(H_1)=rc(G_1)$. We can
suppose that $w=v$. To
get $lc_w(H_2)=lc_v(G_2)u$, let $F$ be $G_2$ with $z_2$ white and verify the following, using local substitution
rules: if $v$ is black and $[u,v]\notin E(F)$, $u[vz_2]\sim_{F}
[vz_2]u$; if $v$ is white and $[u,v]\notin E(F)$, $uvz_2\sim_{F}
vz_2u$; if $v$ is black and $[u,v]\in E(F)$, $uvz_2 \sim_{F}
[vz_2]u$; if $v$ is white and $[u,v]\in E(F)$,
$u[vz_2]\sim_{F}vz_2u$.\\

Case 5) adjacent black $x,y\in V_1\setminus N_{G_1}(z_1)$.\\
$rw(H_1)=rw(G_1)[xy]$, $lw(H_2)=lw(G_2)$, $rb(H_1)=rb(G_1)[xy]$, $lb(H_2)=lb(G_2)$,
$rc(H_1)=rc(G_1)[xy]$, and choosing $w=v$, $lc_w(H_2)=lc_v(G_2)$.\\

Case 6) adjacent black $x,y\in N_{G_1}(z_1)$.\\
$rw(H_1)=rw(G_1)[xy]$, $lw(H_2)=lw(G_2)$, $rb(H_1)=rb(G_1)[xy]$,
$lb(H_2)=lb(G_2)$, $lc_w(H_2)=lc_v(G_2)$, and letting $F$ be $G_1$
with $z_1$ white, we get $rc(H_1)=rc(G_1)uw$ by verifying that $[xy]z_1\sim_{F}
z_1xy$.\\

Case 7) adjacent black $x,y\in N_{G_2}(z_2)$.\\
$rw(H_1)=rw(G_1)$, $lw(H_2)=lw(G_2)[xy]$, $rb(H_1)=rb(G_1)$,
$rc(H_1)=rc(G_1)$. Let $F$ be $G_2$ with $z_2$ white and $F'$ be $G_2$
with $z_2$ black. Since $[xy]z_2\sim_{F} z_2xy$,
$lb(H_2)=lb(G_2)xy$. Choosing $w=x$, since $[xy][xz_2]\sim_{F'}
[yz_2]$, $lc_w(H_2)=lc_v(G_2)$.\\

Case 8) adjacent black $x,y\in V_2 \setminus N_{G_2}(z_2)$.\\
$rw(H_1)=rw(G_1)$, $lw(H_2)=lw(G_2)[xy]$, $rb(H_1)=rb(G_1)$, $lb(H_2)=lb(G_2)[xy]$, $rc(H_1)=rc(G_1)$,
and choosing $w=v$, $lc_w(H_2)=lc_v(G_2)[xy]$.\\

Case 9) adjacent black $x\in V_1\setminus N_{G_1}(z_1)$, $y\in N_{G_1}(z_1)$.\\
$rw(H_1)=rw(G_1)[xy]$, $lw(H_2)=lw(G_2)$, $rb(H_1)=rb(G_1)[xy]$,
$lb(H_2)=lb(G_2)$, $rc(H_1)=rc(G_1)yx$, and choosing $w=v$, $lc_w(H_2)=lc_v(G_2)$.\\

Case 10) adjacent black $x\in V_1$, $y\in V_2$.\\
By Lemma~\ref{recolour}, $rw(H_1)=rc(G_1)xz_1$, $lb(H_2)=lb(G_2)y$ and
$rc(H_1)=rw(G_1)z_1x$. $rb(H_1)=rb(G_1)[xz_1]$. Choosing $w=v=y$,
$lc_w(H_2)=lw(G_2)$ and $lw(H_2)=lc_v(G_2)$.\\

Case 11) adjacent black $x\in N_{G_2}(z_2)$, $y\in V_2\setminus N_{G_2}(z_2)$.\\
$rw(H_1)=rw(G_1)$, $lw(H_2)=lw(G_2)[xy]$, $rb(H_1)=rb(G_1)$,
$rc(H_1)=rc(G_1)$. Let $F$ be $G_2$ with $z_2$ white. Since
$[xy]z_2\sim_{F}z_2xy$, $lb(H_2)=lb(G_2)xy$. Choose $w=y$ and
let $F'$ be $G_2$ with $z_2$ black. Since $[xz_2]\sim_{F'}
[xy][z_2y]$, $lc_w(H_2)=lc_v(G_2)$.\\

\noindent To sum up, complementing at a white $u\in V_1$ adjacent to $V_2$ interchanges the sets
$\{Rw_G,Lw_G\}$ and $\{Rb_G,Lb_G\}$, complementing at a white $u\in
V_2$ adjacent to $V_1$ interchanges $\{Rb_G,Lb_G\}$ and $\{Rc_G,Lc_G\}$, and
complementing with respect to an edge with incident black vertices in
$V_1$ and $V_2$ interchanges $\{Rw_G,Lw_G\}$ and $\{Rc_G,Lc_G\}$.
All other basic (vertex or edge) parity complementations leave the
classes unpermuted.
\end{proof}

The proof of the following is left to the reader:

\begin{lemma} \label{anticlique}
Let $A$ be an independent subset of the black vertices of a bicoloured
graph $G$. Let $v$ be adjacent to $u\in A$. Then $B=A\triangle \{u,v\}$ is an independent subset of the black vertices of
$G\{u,v\}$. Furthermore, $N_G(A)\cup A=N_{G\{u,v\}}(B)\cup B$.
\end{lemma}

\begin{theorem} \label{produit}
Let $G$ be a simple bicoloured graph with a split $V=V_1\cup V_2$ inducing the rooted graphs $(G_1,z_1)$ and $(G_2,z_2)$. If any one of the six essential decompositions of $[G]$ along this split is a pair of pure parity classes, then $[G]$ is pure.
\end{theorem}

\begin{proof}
If $G$ has no edge between $V_1$ and $V_2$, the result is trivial, so let $v\in N_{G_2}(z_2)$.
By way of contradiction, suppose that $G$ is not pure. Without loss of
generality, $G$ admits a black anticlique $A$ and also one of $\{rw(G_1),lw(G_2)\}$,
$\{rb(G_1),lb(G_2)\}$ or $\{rc(G_1),lc_v(G_2)\}$ is a pair of pure graphs. If $A\cap N_{G_1}(z_1)=\emptyset$, then $A\cap V_2$ is a black anticlique of
$lw(G_2)$ and $(A\cap V_1)\cup \{u_1\}$ is a black anticlique of $rb(G_1)$ and $rc(G_1)$. If
not, then $A\cap N_{G_2}(z_2)=\emptyset$, $A\cap V_1$ is a black
anticlique of $rw(G_1)$ and of $rb(G_1)$, and by Lemma~\ref{anticlique},  $(A\cap
V_2)\cup \{v\}$ is a black anticlique of $lc_v(G_2)$.
\end{proof}

Theorem~\ref{produit} allows us to determine numerous pure graphs, for example:

\psfrag{v1}[][][1]{$v_1$}
\psfrag{v2}[][][1]{$v_2$}
\psfrag{v3}[][][1]{$v_3$}
\psfrag{v4}[][][1]{$v_4$}
\psfrag{v5}[][][1]{$v_5$}
\psfrag{v6}[][][1]{$v_6$}
\psfrag{v7}[][][1]{$v_7$}
\psfrag{v8}[][][1]{$v_8$}
\psfrag{v9}[][][1]{$v_9$}
\psfrag{(i)}[][][1]{$(i)$}
\psfrag{(ii)}[][][1]{$(ii)$}
\psfrag{(iii)}[][][1]{$(iii)$}

\begin{center}
\epsfig {file=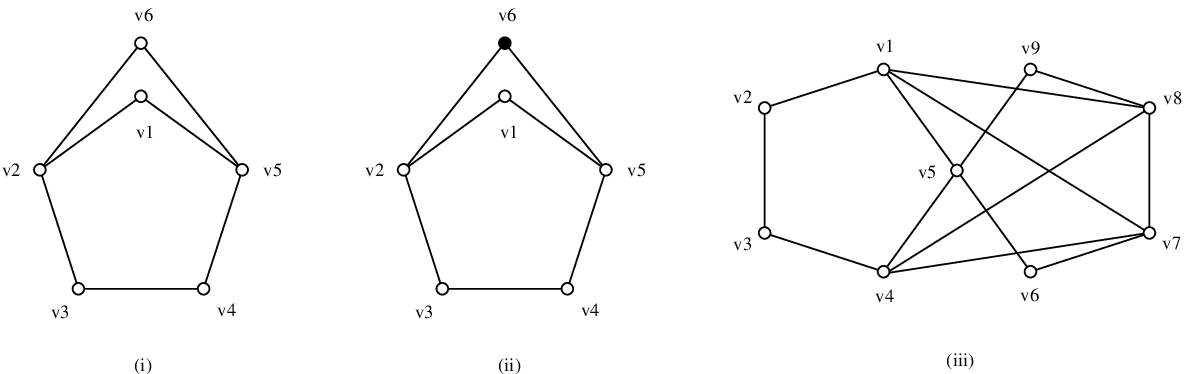}
\end{center}

For $(i)$ and $(ii)$, consider the split $V_1=\{v_2,v_3,v_4,v_5\},
V_2=\{v_1,v_6\}$. For $(iii)$, look at
$V_1=\{v_1,v_2,v_3,v_4\},V_2=\{v_5,v_6,v_7,v_8,v_9\}$.

From example $(i)$, we find that graph $(iv)$ is pure (take
$V_1=\{v_1,v_3,v_4,\linebreak v_5,v_6\}$ and $V_2=\{v_2,v_7\}$) and, by induction, so
is graph $(v)$.

\psfrag{v1}[][][1]{$v_1$}
\psfrag{v2}[][][1]{$v_2$}
\psfrag{v3}[][][1]{$v_3$}
\psfrag{v4}[][][1]{$v_4$}
\psfrag{v5}[][][1]{$v_5$}
\psfrag{v6}[][][1]{$v_6$}
\psfrag{v7}[][][1]{$v_7$}
\psfrag{(iv)}[][][1]{$(iv)$}
\psfrag{(v)C5K2}[][][1]{$(v)\quad C_5\circ \overline{K_2}$}
\psfrag{(iii)}[][][1]{$(iii)$}

\begin{center}
\epsfig {file=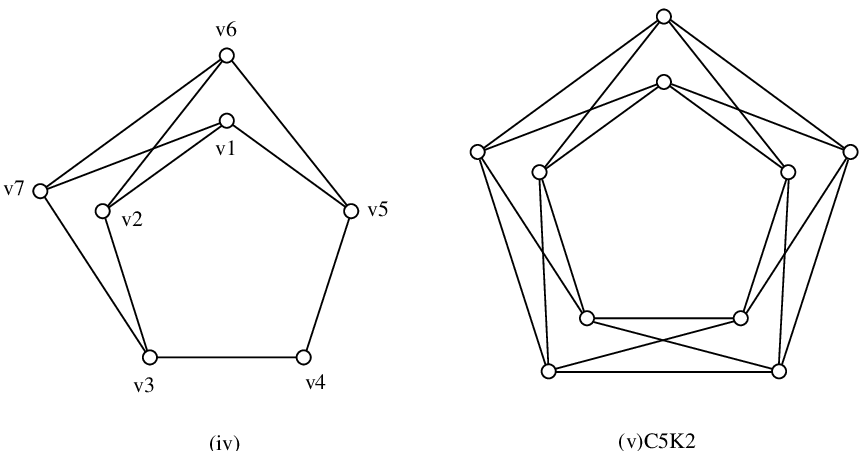}
\end{center}

The last of these naturally coloured pure graphs is of special interest, as we
will see that it is an intermediate step of a correspondence between
$C_5$ and the transition system of $K_5$ consisting of two 5-cycles. The following is easily verified:

\begin{proposition}
Given $H$ with split $V=V_1\cup V_2$ inducing $G_1$ and $G_2$, $H$ is
an alternance graph if and only if $G_1$ and $G_2$ both are.
\end{proposition}

\begin{figure}[htb!]
\begin{center}
\epsfig {file=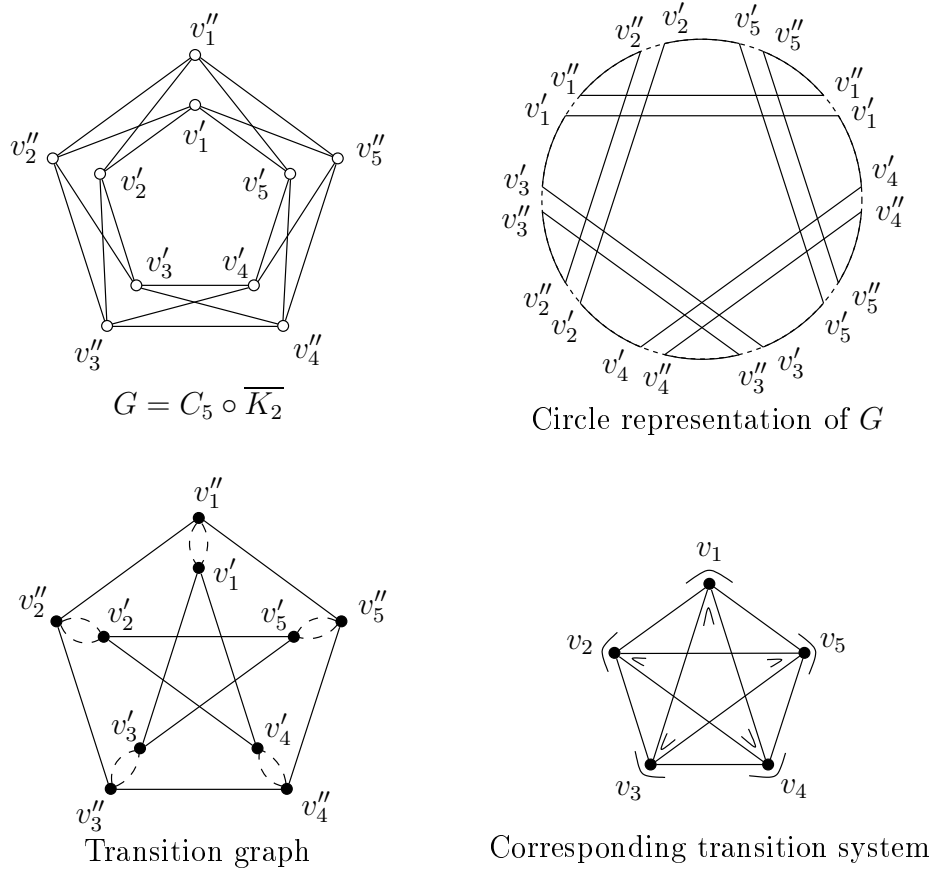}
\end{center}

\caption{\label{C5_K5}Correspondence between $C_5$ and $K_5$ with a transition system made of two 5-cycles.}
\end{figure}

Since the pentagon and $\overline{K_2}$ (two isolated vertices) are alternance graphs, so
is the lexicographic product $C_5\circ \overline{K_2}$. Considering the double occurrence word
giving rise to $C_5\circ \overline{K_2}$ as an Euler tour of a
transition graph, one
of the 2-factors induced by the tour (taking every other edge) is made
up of digons (2-cycles), each corresponding to a vertex of
$C_5$. Contracting this 2-factor results in $K_5$ and determines a transition system
made of two 5-cycles. The different steps of the correspondence
are depicted in Figure~\ref{C5_K5}.

\psfrag{v1}[][][1]{$v_1$}
\psfrag{v2}[][][1]{$v_2$}
\psfrag{v3}[][][1]{$v_3$}
\psfrag{v4}[][][1]{$v_4$}
\psfrag{v5}[][][1]{$v_5$}
\psfrag{v1'}[][][1]{$v_1'$}
\psfrag{v2'}[][][1]{$v_2'$}
\psfrag{v3'}[][][1]{$v_3'$}
\psfrag{v4'}[][][1]{$v_4'$}
\psfrag{v5'}[][][1]{$v_5'$}
\psfrag{v1''}[][][1]{$v_1''$}
\psfrag{v2''}[][][1]{$v_2''$}
\psfrag{v3''}[][][1]{$v_3''$}
\psfrag{v4''}[][][1]{$v_4''$}
\psfrag{v5''}[][][1]{$v_5''$}
\psfrag{G=C5K2}[][][1]{$G=C_5\circ \overline{K_2}$}
\psfrag{Circle representation of G}[][][1]{Circle representation of
  $G$}
\psfrag{TG}[][][1]{Transition graph}
\psfrag{TS}[][][1]{Corresponding transition system}

In the \emph{circle representation} of $C_5\circ \overline{K_2}$,
the vertices are the chords and the edges are chord
intersections. Note that we can read the double occurrence word along the perimeter of the circle. The third graph is obtained by contracting the
chords of the circle representation.\\

The transformation just described, starting from the pentagon and
ending with a transition system of $K_5$, generalizes naturally to a correspondence between parity classes of alternance bicoloured graphs and
transition systems of 4-regular connected graphs.

\begin{definition}
The \emph{double} of a bicoloured graph $G=(V,E)$ is the bicoloured graph
with vertex set $V\times V(\overline{K_2})$, consisting of $G\circ
\overline{K_2}$ and the additional edges
$\{[v',v'']|v \mbox{ black in } G\}$, where $v'$ and $v''$ stand for
$(v,u_1)$ and $(v,u_2)$ and $V(\overline{K_2})=\{u_1,u_2\}$. The
colour of $v'$ and $v''$ is chosen to be the same as that of $v$.
\end {definition}

\begin{proposition} \label{double}
A bicoloured graph is pure if and only if its double is pure.
\end{proposition}

\begin{proof}
We leave it to the reader to verify that if $H$ is the double of $G$
and $u$ is white in $G$, or $v$ and $w$ are adjacent black vertices of
$G$, then $Hu'$ is the double of $Gu$, and $H[v'w']$ is the double of
$G[vw]$, respectively. Also, if $v$ is black in $G$ then
$H[v'v'']=H$. By the symmetry of $\overline{K_2}$, this shows that $[G]$ is in bijection with
$[H]$. Furthermore, a black anticlique of $H$ determines a black anticlique of $G$ using the
projection $v',v''\mapsto v$, and a black anticlique $A$ of $G$ can be used to find a black anticlique
of $H$ by taking $v'$ for each $v\in A$.
\end{proof}

Note that a double is necessarily naturally coloured, and that the
double of an alternance graph is also an alternance graph. The
correspondence we mentioned goes like this:


\begin{center}
parity class of a bicoloured alternance graph\\
$\updownarrow$\\
parity class of its double\\
$\uparrow$\\
TG with a 2-factor consisting of digons\\
$\updownarrow$\\
TS of a 4-regular connected graph
\end{center}

Note that to generalize this to a 1-1 correspondence, we have to associate to each bicolored
alternance graph one of its circle representations. Using Proposition~\ref{double}, we can now state the following fundamental
relationship:

\begin{theorem}A transition system of a connected $4$-regular graph
admits no orthogonal cycle decomposition if and only if the corresponding
parity class is pure.
\end{theorem}

\begin{figure}[htb!]
\begin{center}
\[ \epsfig{file=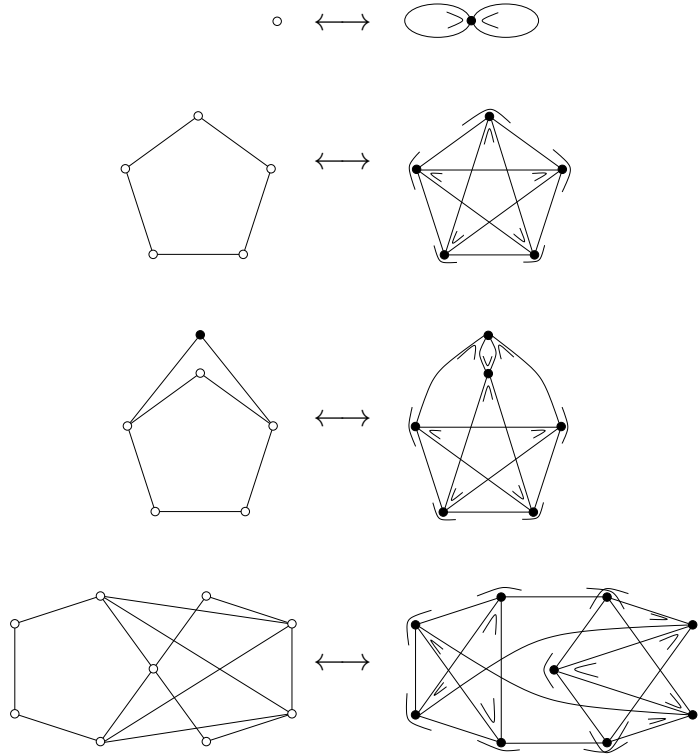,height=10cm}
\]

\caption{\label{purealt}Some pure alternance graphs and corresponding transition systems.}
\end{center}
\end{figure}

Figure~\ref{purealt} shows some examples of pairs (parity class, TS), where the parity class is represented by one of its members. Note that the correspondence induces a bijection
between the vertex set of the 4-regular graph and the vertex set of the bicoloured
graph. The graphs are drawn to reflect this relationship.

\psfrag{equiv}[][][1]{$\longleftrightarrow$}

In the last two pure graphs, hiding on
one side of the split, we recognize the pentagon. In the corresponding
transition system, this pentagon becomes Fleischner's tetrapus, and to each split,
there corresponds a non-trivial edge-cut of size $<6$. This is no
coincidence.

\begin{proposition}
Let $[G]$ be the parity class of a bicoloured alternance graph, with corresponding transition system defined on $\Gamma$. A partition $V(G)=V_1\cup V_2$ is a
split of $G$ if and only if it determines a non-trivial edge-cut
$\mathcal{C}$ of $\Gamma$ of size $2$ or $4$, and
$|V_1|,|V_2|\geq 2$. In this case, $|\mathcal{C}|=2$ if and only if there is no
edge between $V_1$ and $V_2$.
\end{proposition}

\begin{proof}
$\Gamma$ being connected and 4-regular, the partition determines an edge
cut $\mathcal{C}$ containing an even number of edges. $\mathcal{C}$ is also an edge cut of
the transitition graph, partitioning $S$ into $S_1\cup S_2$, where $t\in S_1 \iff t\cap V \in
V_1$. Let $w$ be an alternating tour of the transition graph.
Then the following statements are equivalent:\\

\noindent (1) $|\mathcal{C}|=2$;\\
(2) the tour is of the form $w=w_1w_2$ where $w_1\in S_1^*,
w_2\in S_2^*$ (words in the alphabets $S_1$ and $S_2$ respectively);\\
(3) there is no edge between $V_1$ and $V_2$ in $G$.\\

\noindent We can also verify the equivalence of:

\noindent (4) $|\mathcal{C}|=4$;\\
(5) $w=w_1w_2w_3w_4$ where $w_1,w_3\in
S_1^*\setminus\{\emptyset \},w_2,w_4\in S_2^*\setminus\{\emptyset \}$;\\
(6) $N_G(V_1)\cap V_2, N_G(V_2)\cap V_1$ are non-empty and there is
every possible edge between the two sets.\\

\noindent Furthermore, we have the equivalence of:

\noindent (7) $\mathcal{C}$ is a trivial edge cut of $\Gamma$ (say
$\mathcal{C}=N_{\Gamma}(u)$);\\
(8) without loss of generality $S_1=\{t_1,t_2\}$
where $u\in t_1\cap t_2$;\\
(9) each of $w_1,w_3$ is either $t_1t_2$ or $t_2t_1$;\\
(10) $V_1 = \{u\}$.
\end{proof}

\section{An interpretation of Theorem~\ref{produit}}

At this point, let us revisit Definition~\ref{root} and give an interpretation of the root and leaf graph constructions in terms of transition systems when the rooted graph has the structure of an alternance graph:

\begin{definition}
Given an eulerian graph $\Gamma$ of minimum degree $\delta >2$, a set
of transitions $S_z$ of $\Gamma$ is a \emph{transition system rooted
  at} $z\in V(\Gamma)$ if no transition of $S_z$ contains $z$ and $S_z$ can be completed to a transition system of $\Gamma$ by adding transitions at $z$. Such a set is also called a $z$-\emph{transition system}.
\end{definition}

\begin{figure}[htb!]
\begin{center}
\epsfig{file=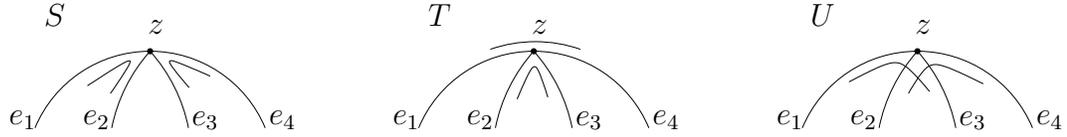,width=14cm}
%
\caption{\label{transitions}The possible pairs of transitions at $z$.}
\end{center}
\end{figure}

Note that a $z$-transition system completely determines $\Gamma$ so
that we can talk about tours and cycles of rooted transition
systems. Let $(G,z)$ be a rooted graph such that the
underlying simple graph is an alternance graph. Let $S$ be the
transition system corresponding to $rw(G)$. First consider the
situation where $z$ is not isolated in $G$. Let $t_1=\{z,e_1,e_2\}$
and $t_2=\{z,e_3,e_4\}$ be the transitions of $S$ containing $z$. We
then have that $S_z=S\setminus \{t_1,t_2\}$ is a $z$-transition
system. $S_z$ can be completed to a transition system in three different ways, as shown in Figure~\ref{transitions}.

\psfrag{m}[][][1]{$z$}
\psfrag{e1}[][][1]{$e_1$}
\psfrag{e2}[][][1]{$e_2$}
\psfrag{e3}[][][1]{$e_3$}
\psfrag{e4}[][][1]{$e_4$}
\psfrag{S1}[][][1]{$S$}
\psfrag{S2}[][][1]{$T$}
\psfrag{S3}[][][1]{$U$}

The associated transition graphs differ as follows:
\newpage
\psfrag{t1}[][][1]{$t_1$}
\psfrag{t2}[][][1]{$t_2$}
\psfrag{t3}[][][1]{$t_3$}
\psfrag{t4}[][][1]{$t_4$}
\psfrag{t5}[][][1]{$t_5$}
\psfrag{t6}[][][1]{$t_6$}
\psfrag{e'}[][][1]{$e'$}
\psfrag{e''}[][][1]{$e''$}
\psfrag{T1}[][][1]{$\Gamma_S$}
\psfrag{T2}[][][1]{$\Gamma_T$}
\psfrag{T3}[][][1]{$\Gamma_U$}

\begin{figure}[htb!]
\begin{center}
\epsfig{file=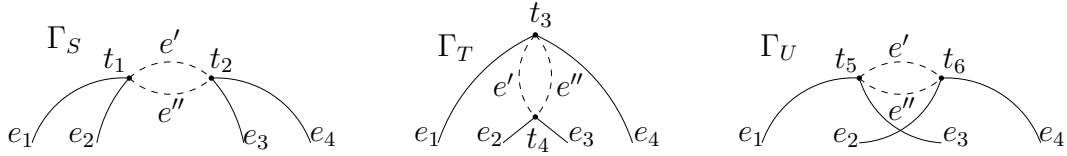,width=14cm}
%
\caption{The corresponding transition graphs.}
\end{center}
\end{figure}

Let $\alpha$ be an alternating Euler tour of $\Gamma_S$ corresponding to $rw(G)$. Since $z$ is white in $rw(G)$, $\alpha$ is of the form $t_1e't_2e_3\sigma e_4t_2e''t_1e_2\tau e_1$, where $\sigma$ and $\tau$ are appropriate sequences of vertices and edges. We then have that the alternating Euler tours $\beta=t_3e't_4e_3\sigma e_4t_3e''t_4e_2\tau e_1$ of $\Gamma_T$ and $\gamma =t_5e't_6e_4\sigma^{-1}e_3t_5e''t_6e_2\tau e_1$ of $\Gamma_U$ (where $\sigma^{-1}$ is $\sigma$ in reverse order) correspond to the root graphs $rb(G)$ and $rc(G)$, respectively.

\begin{figure}[htb!]
\begin{center}
\epsfig{file=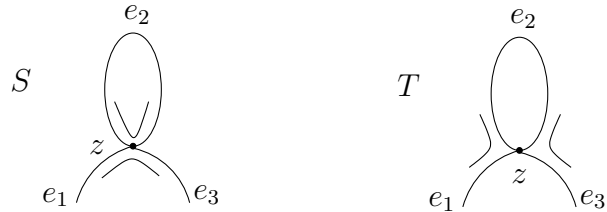,width=8cm}
%
\caption{\label{isolated}When $z$ is isolated in the bicoloured graph.}
\end{center}
\end{figure}

In the case where $z$ is isolated in $G$, let $t_1=\{z,e_1,e_3\}$ and $t_2=\{z,e_2\}$ be the transitions of $S$ at $z$. The $z$-transition system $S_z=S\setminus \{t_1,t_2\}$ can be completed to a transition system in only one other way, as seen in Figure~\ref{isolated}.

\begin{figure}
\begin{center}
\epsfig{file=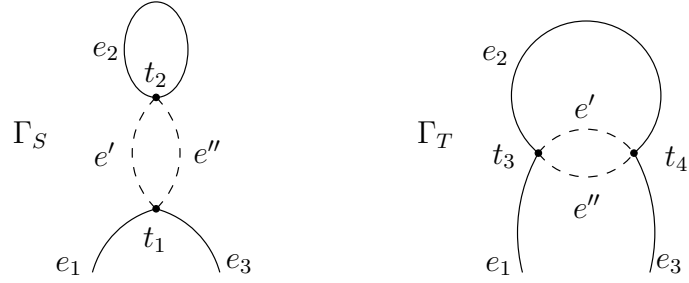,width=9cm}
%
\caption{The corresponding transition graphs.}
\end{center}
\end{figure}

Given an Euler tour $\alpha =t_1e't_2e_2t_2e''t_1e_3\tau e_1$ of $\Gamma_S$ corresponding to $rw(G)$, $\beta=t_3e't_4e_2t_3e''t_4e_3\tau e_3$ is an Euler tour of $\Gamma_T$ corresponding to $rb(G)$.

Thus we see that to each rooted alternance graph $(G,z)$ there corresponds
a unique connected 4-regular rooted transition system $S_z$, and that
the root graphs of $G$ are in correspondence with the
different ways the rooted transition system can be completed. A
similar analysis reveals that the leaf graphs of $G$ are in
correspondence with the different transition systems that can be
obtained from $S_z$ by identifying the edges incident with $z$ two by
two.

We are now in a position to interpret Theorem~\ref{produit} in
terms of transition systems. Let a bicoloured alternance graph $G$ have an essential decomposition into a pair of pure root and leaf graphs $H_1$ and $H_2$, where $V(G)=V_1\cup V_2$ is a split of $G$ inducing the rooted graphs $(G_1,z_1)$ and $(G_2,z_2)$, such that $H_1$ is a root graph of $G_1$, and $H_2$ is a leaf graph of $G_2$. Let $S_{z_1}$ and $S_{z_2}$ be the rooted transition systems corresponding to $G_1$ and $G_2$. It can be verified that if the transition system corresponding to $H_1$ is $S_{z_1}\cup \{\{z_1,e_1,e_2\},\{z_1,e_3,e_4\}\}$, then the transition system corresponding to $H_2$ is obtained from $S_{z_2}$ by identifying $e_1$ with $e_2$, and $e_3$ with $e_4$.

\section{Primitivity of pure bicoloured graphs}

\begin{definition}
A family $\mathcal{F}$ of parity classes is \emph{closed}, if for each bicoloured graph $G$, whenever\\

i) $G$ admits an essential decomposition $\{H_1,H_2\}$ such that $[H_1],[H_2]\in \mathcal{F}$, or \\

ii) $G$ has a component $G_1$ such that $[G_1]\in \mathcal{F}$,\\

\noindent then $[G]\in \mathcal{F}$.
\end{definition}

The intersection of closed families of parity classes is closed, hence given any family $\mathcal{F}$ of parity classes, we can define the \emph{closure} $\overline{\mathcal{F}}$ of $\mathcal{F}$ to be the smallest closed family of parity classes containing $\mathcal{F}$ as a subset.

Note that similar definitions can be made when we restrict our attention to parity classes of bicoloured alternance graphs. Accordingly, some of the following results have an equivalent formulation in the context of alternance graphs.

\begin{definition}
A pure graph $G$ (or its class $[G]$) is \emph{primitive} if for any family $\mathcal{F}$ of pure
classes, $[G]\in \overline{\mathcal{F}}$ implies $[G]\in
\mathcal{F}$.
\end{definition}

\begin{conjecture}\label{primitive_prime}
All primitive pure graphs are prime.
\end{conjecture}

Recall that from a rooted graph we can construct up to three different
root graphs (two if the root vertex is isolated). Of course, any bicoloured
graph $G$ can be viewed as a root graph: for a given vertex $u$ of $G$,
let $(G,u)$ be the rooted graph obtained by restricting the colouring of $G$ to $V(G)\setminus \{u\}$, then $G$ is either $rw(G')$ or $rb(G')$. We call the root graphs constructed from $(G,u)$ the root graphs \emph{induced} by $u$ from $G$.

\begin{definition}
Given a pure graph $G$, a vertex $u$ of $G$ is \emph{critical} if one of the root graphs induced by $u$ is not pure.
\end{definition}

Given a bicoloured graph $G$ and $u\in V(G)$ inducing the root graphs
$G_1,G_2$ and $G_3$, we have, as a corollary of Theorem~\ref{invariant},
that $S=\{[G_1],[G_2],[G_3]\}$ is invariant over $[G]$ (i.e. if $u$ induces the root graphs $H_1,H_2$ and $H_3$ from $H\in [G]$, then $\{[H_1],[H_2],[H_3]\}=S$). This gives us:

\begin{proposition} \label{critical}
The set of critical vertices of a pure graph $G$ is invariant over $[G]$.
\end{proposition}

\begin{proposition}
Let $u$ be a critical vertex of a pure graph $G$. Of the root graphs induced by $u$, only $G$ is pure.
\end{proposition}

\begin{proof}
If $u$ is isolated, the result is trivial. Suppose, by way of contradiction,
that $u$ is not isolated and that only one root graph induced by $u$
is not pure. Let $G'$ be the rooted graph obtained by restricting the
bicolouring of $G$ to $V(G)\setminus \{u\}$. By
Proposition~\ref{critical}, we can suppose that one of $rw(G')$,
$rb(G')$ or $rc(G')$ has a black anticlique $A$. If $A$ is a black
anticlique of $rw(G')$, then it is also a black anticlique of
$rb(G')$. If $A$ is a black anticlique of $rb(G')$ but not of
$rw(G')$, then it is one of $rc(G')$. If $A$ is a black anticlique of
$rc(G')$ but not of $rb(G')$, then it is one of $rw(G')u$. All
cases lead to a contradiction.
\end{proof}

\begin{definition}
A vertex $u$ of a pure graph $G$ is said to be \emph{tight} if there exists $H\in [G]$ such that $H-u$ has a black anticlique. A pure graph is \emph{tight} if every critical vertex is tight. The parity class is then also said to be \emph{tight}.
\end{definition}

Obviously, a tight vertex is also critical. It is easy to verify that
the one-point graph and the pentagon with their natural colourings are tight.

\begin{proposition}
A vertex $u$ of a pure graph $G$ is tight if and only if $G-u$ is not pure.
\end{proposition}

\begin{proof}
One side is clear. Let $S$ be a complementation set of $G$ (see [8]) such that $GS-u$ has a black anticlique $A$. If $u\notin S$, then $S$ is a complementation set of $G'=G-u$ and $A$ is a black anticlique of $G'S$ ($=GS-u$). If $u\in S$, since $u$ is white in $GS$ (by the purity of $G$), $S'=S\setminus \{u\}$ is a complementation set of $G$, and then also a complementation set of $G'=G-u$. Since $u\notin N_{GS}(A)$ (again by the purity of $G$), $A$ is a black anticlique of $G'S'$.
\end{proof}

Consider the closure $\overline{\mathcal{F}}$ of a family
 $\mathcal{F}$ of parity classes. Given a parity class $[G]\in
 \overline{\mathcal{F}}$, the following possibilities may arise:\\

(1) $[G]$ is in $\mathcal{F}$;\\

(2) $G$ is disconnected, in which case $G$ has a connected component whose
 parity class is in $\overline{\mathcal{F}}$; or\\

(3) $G$ has an essential decomposition into two graphs $H_1$ and $H_2$ such that
 $[H_1],[H_2]\in \overline{\mathcal{F}}$.\\

\noindent This analysis can be applied recursively to the resulting smaller graphs (``factors'') with parity classes in
 $\overline{\mathcal{F}}$. The whole process is called an
 $\mathcal{F}$-\emph{factorization} of $G$ (or more precisely, an \emph{essential} $\mathcal{F}$-\emph{factorization}). For the purposes of the following
 proposition, given $[G]\in \overline{\mathcal{F}}$, where
 $\mathcal{F}$ is a family of pure graphs, call a vertex $u$ of $G$
 \emph{factor-essential} if every pure factor $H$ of an
 $\mathcal{F}$-factorization such that $u\in H$ is connected. 

\begin{proposition} \label{facteur}
Let $\mathcal{F}$ be a family of parity classes that are primitive, pure and tight. Let $[G]\in \overline{\mathcal{F}}$ be such that
all primitive pure graphs of smaller order than $G$ have their parity
class in $\mathcal{F}$. If $u\in V(G)$ is factor-essential, then $u$ is a tight vertex of $G$.
\end{proposition}

\begin{proof}
First note that any critical vertex is factor-essential, so that the
truth of the statement for a given graph implies that it is tight. We
proceed by induction on the order of $G$. The statement is true if
$[G]\in \mathcal{F}$. If $G$ is not connected, then, by the induction
hypothesis, $u$ is a tight vertex of the component that contains
it. Since the other components are not pure, $u$ is also tight in
$G$. If $G$ is connected but $[G]\notin \mathcal{F}$, then it has an
essential decomposition into the pure root graph $H_1$ and the pure
leaf graph $H_2$ where $[H_1],[H_2]\in
\overline{\mathcal{F}}$. Suppose that the corresponding split is
$V(G)=V_1\cup V_2$, inducing the rooted graphs $(G_1,z_1)$ and $(G_2,z_2)$. Without loss of generality, $H_1=rw(G_1)$
and $H_2=lw(G_2)$ (because of the presence of white vertices in $H_1$
and $H_2$, we can always locally complement $G$ to reach this kind of
decomposition).

Consider the case where $u\in V_2$. Since $u$ is
factor-essential in $G$, $z_1$ is factor-essential in $H_1$. By the
induction hypothesis, $z_1$ is a tight vertex of $H_1$ and we
can suppose that $H_1-z_1$ has a black anticlique. Since $u$ has to be
factor-essential in $H_2$, it is a tight vertex of $H_2$ and we can
suppose that $H_2-u$ has a black anticlique. Taking the union of the
anticliques, we see that $u$ is a tight vertex of $G$.

There remains the case where $u\in V_1$. Since $u$ is factor-essential in $G$, it is
factor-essential in $H_1$ and thus tight in $H_1$. Let $S$ be a
complementation set of $H_1-u$ such that $H_1S-u$ has a black
anticlique $A$. We can suppose that $z_1\notin A$: if it is, take
$v\in N_{H_1S}(z_1)$ (which is not empty, since $H_1$ is connected, and
does not contain $u$, by the purity of $H_1$) and replace $S$ with
$S'=S\triangle \{z_1,v\}$, then $A'=A\triangle \{z_1,v\}$ is a black
anticlique of $H_1S'-u$ such that $z_1\notin A'$. Suppose that either
one of $rw(G_2)$ or $rb(G_2)$ is pure. Then the purity of $H_2$
implies that $z_2$ is not tight (therefore not critical) and both
$rw(G_2)$ and $rb(G_2)$ are pure. However, this would mean that $z_2$
is not factor-essential in $rw(G_2)$, which would contradict the fact
that $u$ is factor-essential in $G$. Thus neither one of $rw(G_2)$ and
$rb(G_2)$ is pure. If $z_1\in S$, let $S'$ be a complementation set of
$F=rw(G_2)$ such that $FS'$ has a black anticlique $A'$. Without loss
of generality (following an argument similar to the $z_1$ case),
$z_2\notin S'$. By the purity of $H_2$, $z_2\in A'$ (i.e. it has
changed colour). Then $S''=S\cup S'$ is a complementation set of $G$
such that $A\cup A'\setminus \{z_2\}$ is a black anticlique of $GS''-u$, and $u$ is a
tight vertex of $G$. If $z_1\notin S$, let $S'$ be a
complementation set of $F=rb(G_2)$ such that $FS'$ has a black
anticlique $A'$. Without loss of generality, $z_2\notin S'$. By the
purity of $H_2$, $z_2\in A'$ (its colour has not changed). Then
$S''=S\cup S'$ is a complementation set of $G$ such that $A\cup
A'\setminus \{z_2\}$ is a black anticlique of $GS''-u$, and $u$ is a tight vertex of $G$.
\end{proof}

\begin{corollary} \label{prim}
Let $\mathcal{F}$ be a family of parity classes that are primitive, pure and tight. Let $G$ be a smallest primitive pure graph such that $[G]\notin \mathcal{F}$. Then $G$ has no essential decomposition into a root graph $H_1$ and a leaf graph $H_2$ such that $H_1$ is pure.
\end{corollary}

\begin{proof}
Suppose, by contradiction, that $V(G)=V_1\cup V_2$ is a partition inducing the rooted graphs $(G_1,z_1)$ and $(G_2,z_2)$ and such that $H_1$ is a pure root graph of $G_1$ and $H_2$ is the corresponding leaf graph of $G_2$. The minimality of $G$ forces $H_1$ to be in $\overline{\mathcal{F}}$. If $z_1$ is not critical in $H_1$, then it has an essential factorization such that at least one pure graph contains $z_1$ but has a pure component not containing $z_1$. But then $G\in \overline{\mathcal{F}}$, contradicting the primitivity of $G$. Therefore $z_1$ is a tight vertex of $H_1$, and $G$ is pure if and only if $H_2$ is pure, contradicting either the purity or the primitivity of $G$.
\end{proof}

A corresponding statement is obtained by considering only alternance graphs:

\begin{corollary} \label{altern}
Let $\mathcal{F}$ be a family of alternance parity classes (parity classes of bicoloured alternance graphs) that are primitive, pure and tight. Let $G$ be a smallest primitive pure alternance graph such that $[G]\notin \mathcal{F}$. Then $G$ has no essential decomposition into a root graph $H_1$ and a leaf graph $H_2$ such that $H_1$ is pure.
\end{corollary}

\begin{theorem} \label{primitif}
A smallest primitive pure alternance graph different from the white one-point graph and the
graphs in the parity class of the white pentagon must be prime.
\end{theorem}

\begin{proof}
Suppose, by way of contradiction, that such a bicoloured graph $G$ has a
split $V=V_1\cup V_2$ inducing the rooted graphs $(G_1,z_1)$ and $(G_2,z_2)$. Let $S_{z_1}$ be the
$z_1$-transition system corresponding to $G_1$ and $S_{z_2}$ the
$z_2$-transition system corresponding to $G_2$. We have seen that the
transition system $S$ corresponding to $G$ has an edge-cut of size 2
or 4 separating $V_1$ from $V_2$. In fact, since $G$ is primitive
(and thus connected), this cut contains four edges. Clearly, these are
the edges, say $e_1,e_2,e_3$ and $e_4$, that are incident with $z_1$ in $S_{z_1}$ and incident with $z_2$ in $S_{z_2}$. Consider the number $n(S_z)$ of partitions into pairs of the edges incident with $z$ that can be extended to an orthogonal cycle decomposition of the rooted transition system $S_z$. By Corollary~\ref{altern}, $n(S_{z_1}),n(S_{z_2})>1$. By the pigeonhole principle, at least one partition of the set $\{e_1,e_2,e_3,e_4\}$ into pairs of edges can be extended to an orthogonal cycle decomposition of $S$.
\end{proof}

The following is essentially Conjecture 12 in [6], due to Fleischner
and {\nobreak Jackson}:

\begin{conjecture} \label{edge_cut}
If $\Gamma \neq K_5$ (with $\delta > 2$) has no non-trivial edge cut
of size $<6$, then any transition system of $\Gamma$ admits an orthogonal cycle decomposition.
\end{conjecture}

\begin{conjecture} \label{nada}
The only connected prime pure alternance graphs are the white one-point graph and those in the parity class of the white pentagon.
\end{conjecture}

The two conjectures are equivalent: a connected prime pure alternance graph gives rise to a 4-regular
transition system without non-trivial edge cuts of size $<6$, while a counter-example
to Conjecture~\ref{edge_cut} would imply the existence of at least one new primitive pure
alternance graph, the smallest of which must be prime. In addition to
providing a characterization of transition systems with orthogonal cycle
decompositions, the truth of Conjectures~\ref{edge_cut} and \ref{nada}
would imply the following:

\begin{conjecture}[Cycle Double Cover Conjecture]
For every bridgeless graph (i.e. without edge-cuts of size one) there
is a family of cycles of this graph such that every edge belongs to exactly two cycles of the family.
\end{conjecture}

The idea is that a counter-example $G$, minimal
with respect to the number of edges, would be 3-regular and cyclically 6-edge-connected. Contracting a 1-factor (which exists by Petersen's Theorem)
and choosing the transition system induced by the remaining 2-factor of $G$, we get a 4-regular graph without non-trivial edge cuts of size
$<6$. By Conjecture~\ref{edge_cut}, we obtain a cycle decomposition orthogonal to
the transition system. Completing it with edges of the 1-factor, this decomposition
extends to a family of cycles of $G$ covering the 1-factor twice and
the 2-factor once. Throwing in the cycles of the 2-factor, we get a
cycle double cover of $G$. Alternatively, as pointed out by Jaeger
(see [11]), we can take
the line graph of $G$ with transition system inducing the triangles
corresponding to the vertices of $G$. This is another transition system of a
4-regular graph without non-trivial edge cuts of size $<6$. A cycle
decomposition orthogonal to this transition system extends naturally to a cycle double
cover of $G$, and vice versa.

\begin{figure}[htb!]
\begin{center}
\epsfig{file=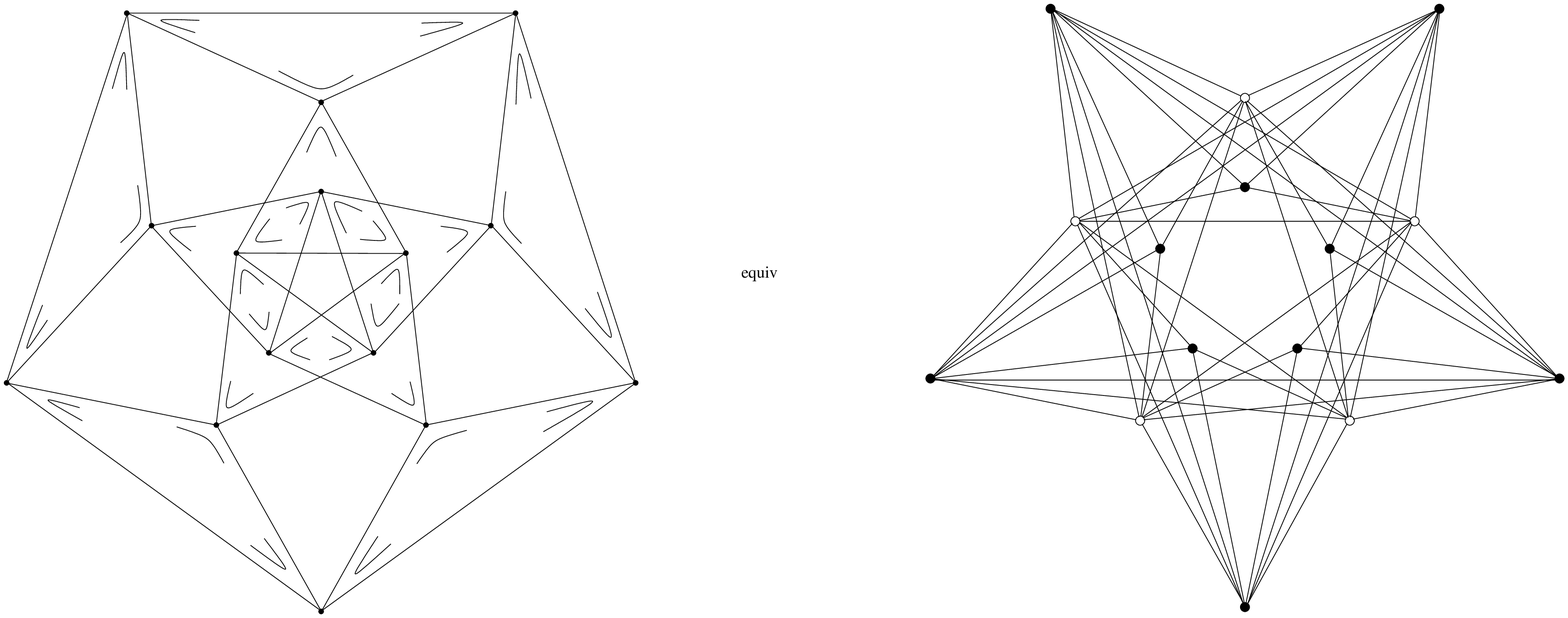,height=5.5cm}

\epsfig{file=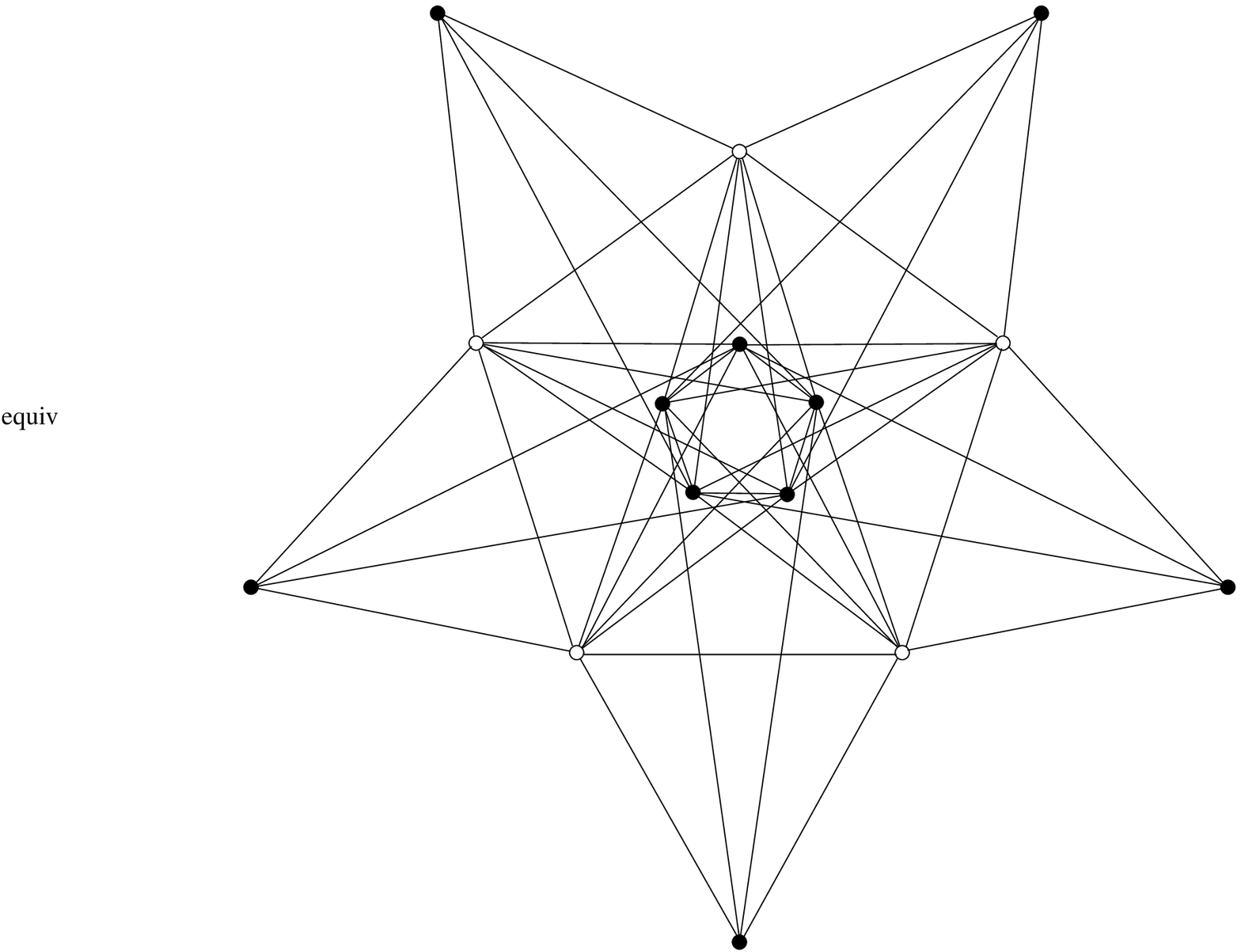,height=5.5cm}
\caption{\label{Petersen}A transition system constructed from the Petersen graph and some associated bicoloured graphs.}
\end{center}
\end{figure}

To illustrate the preceding argument, here is how we can rule out the Petersen
graph as a counter-example to the Cycle Double Cover Conjecture:
contracting a 1-factor yields our familiar $K_5$ with a 5-cycle
decomposition, corresponding to the pure class of the pentagon. However, taking the line graph
approach, we obtain a parity class which admits black anticliques.

In Figure~\ref{Petersen}, we give two highly symmetric representatives of
this class. Black anticliques can be found by taking the inner five vertices of the first graph or the outer five vertices of the second one.

It is easy to see from the correspondence that it suffices to prove Conjecture~\ref{edge_cut} for 4-regular graphs. Furthermore,
any counter-example must contain a $K_5$-minor (Fan and Zhang [4]).

Two more primitive pure classes are known. We have that $[Z_{13}]$ and
$[Z_{17}]$ are pure classes, where $Z_{13}$ and $Z_{17}$ are the naturally coloured Cayley graphs Cay$(\mathbb{Z}_{13},\pm \{1,3,4\})$ and Cay$(\mathbb{Z}_{17},\pm \{1,5\})$. However, while prime, the members of $[Z_{13}]$ and $[Z_{17}]$ are not alternance graphs.

\begin{figure}[htb!]
\begin{center}
\epsfig{file=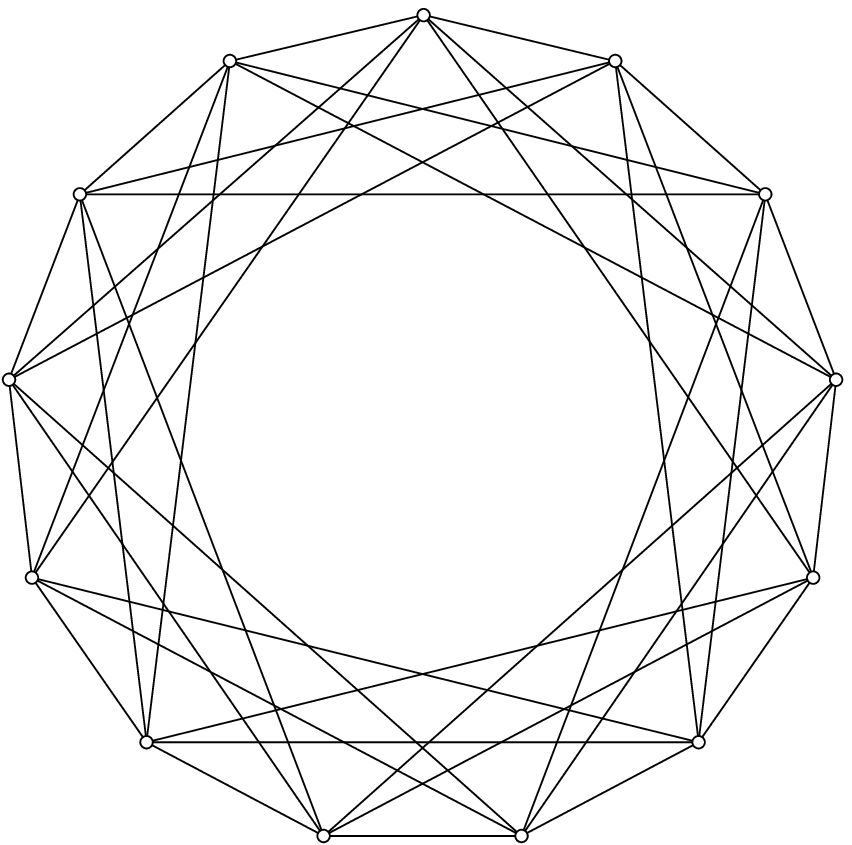,height=5cm}
\qquad
\epsfig{file=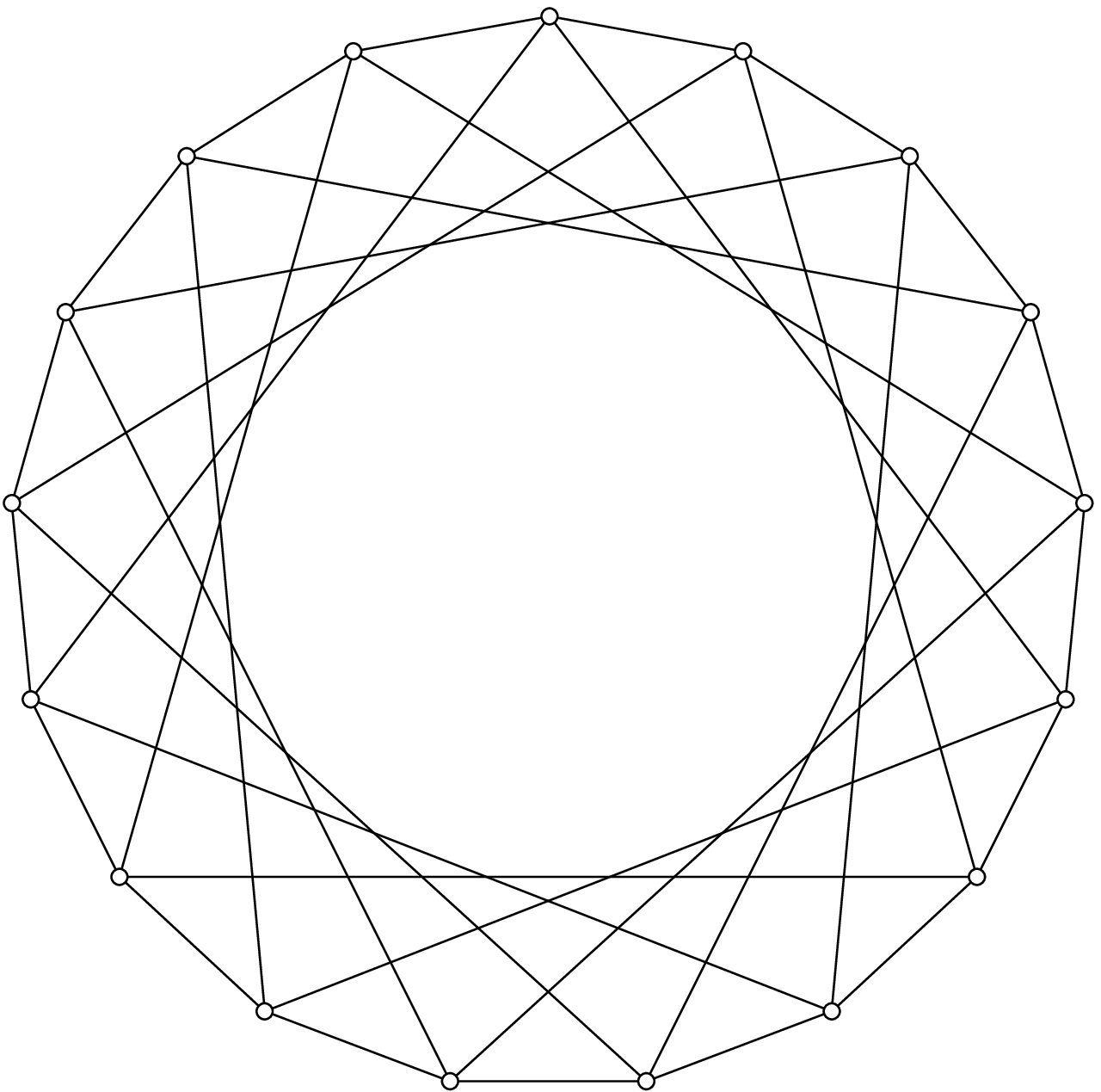,height=5cm}

\caption{The primitive pure graphs $Z_{13}$ and $Z_{17}$.}
\end{center}
\end{figure}

$[Z_{13}]$ and $[Z_{17}]$ contain 39 and 1069 graphs respectively, up to
isomorphism. Many
more graphs were tested by computer, with the focus on
vertex-transitive graphs. I tend to believe these four primitive pure classes are exhaustive but
this is definitely an open question. Parity classes are subsets of
complementation classes, in which complementation is allowed at black
vertices as well as at white vertices. Isotropic systems are a
generalization of 4-regular graphs and dual binary matroids. Bouchet [2]
showed how the set of isotropic systems can be put in bijection with
the set of complementation classes of naturally coloured
graphs. Perhaps a new algebraic object akin to isotropic systems can be put
in bijection with parity classes. This might help to characterize
primitive pure graphs and settle the Cycle Double Cover Conjecture.

\section{References}

\setlength{\parindent}{0pt}

[1] Bouchet, A., Reducing prime graphs and recognizing circle graphs,
\emph{Combinatorica} {\bf 7} (1987), 243-254.\\

[2] Bouchet, A., Graphic presentations of isotropic systems,
\emph{J. Combin. Theory Ser. B} {\bf 45} (1988), 58-76.\\

[3] Bouchet, A., Circle graph obstructions, \emph{J. Comb. Theory Series B}
{\bf 60} (1994), 107-144.\\

[4] Fan, G.; Zhang, C.-Q., Circuit decompositions of eulerian graphs, \emph{J. Combin. Theory
  Ser. B.} {\bf 78} (2000), 1-23.\\

[5] Fleischner, H., Cycle decompositions, $2$-coverings, removable
cycles, and the four-color-disease, \emph{Progress in graph theory
  (Waterloo, 1982)}, Academic Press, Toronto (1984), 233-246.\\

[6] Fleischner, H., Some blood, sweat, but no tears in eulerian graph
Theory, \emph{Congr. Numer.} {\bf 63} (1988), 8-48.\\

[7] de Fraysseix, H., A characterization of circle graphs,
\emph{European J. Combin.} {\bf 5} (1984), 223-238.\\

[8] Genest, F., Word and set complementation of graphs,
  invertible graphs, submitted.\\

[9] Jackson, B., A characterization of graphs having three pairwise
compatible Euler tours, \emph{J. Combin. Theory Ser. B} {\bf 53}
(1991), 80-92.\\

[10] Jackson, B., On circuit covers, circuit decompositions and Euler
tours of graphs, \emph{Surveys in Combinatorics (Keele, 1993)}, London Math. Soc. Lecture Note Ser. {\bf 187}, Cambridge Univ. Press,
Cambridge (1993), 191-210.\\

[11] Jaeger, F., A survey of the cycle double cover conjecture, \emph{Cycles in graphs (Burnaby, 1982)}, North-Holland Math. Stud. {\bf 115}, North-Holland, Amsterdam (1985), 1-12.\\

[12] Kotzig, A., Eulerian lines in finite $4$-valent graphs and their
transformations, \emph{Theory of Graphs (Tihany, 1966)}, Proc. Colloqium on Graph Theory Tihany 1966,
Academic Press (1968), 219-230.\\

[13] Kotzig, A., \emph{Quelques remarques sur les transformations $\kappa$}, séminaire Paris (1977).\\

[14] Sabidussi, G., \emph{Eulerian walks and local complementation},
D.M.S. 84-21, Dép. de math. et stat., Université de Montréal (1984).\\

[15] Spinrad, J., Recognition of circle graphs, \emph{J. Algorithms}
{\bf 16} (1994), 264-282.\\


%% file: conclusion.tex
Le théorême de décomposition essentielle (théorême~\ref{produit} du deuxième
article) est le théorême fondamental de la théorie des graphes
purs. Grâce à lui, l'étude des graphes purs se résume à l'étude des
graphes purs primitifs. Nous proposons que ces derniers sont premiers
(conjecture~\ref{primitive_prime} du même article).

Étant donné la correspondance canonique entre les systèmes de transitions
4-réguliers connexes et les classes de parité de graphes de cordes
bicoloriés, l'étude des systèmes de transitions sans décompositions en
cycles est équivalente à l'étude des graphes de cordes purs primitifs.

Une meilleure compréhension de la pureté passe probablement par une
généralisation des systèmes de transitions 4-réguliers connexes. Le
principal intérêt d'une caractérisation des graphes purs (ne serait-ce
que celle des graphes de cordes purs) réside dans son éventuelle
contribution à la résolution de la conjecture de double recouvrement.

%% file: annexe1.tex
\psfrag{0}[][][.7]{$0$}
\psfrag{1}[][][.7]{$1$}
\psfrag{2}[][][.7]{$2$}
\psfrag{3}[][][.7]{$3$}
\psfrag{4}[][][.7]{$4$}
\psfrag{5}[][][.7]{$5$}
\psfrag{6}[][][.7]{$6$}
\psfrag{7}[][][.7]{$7$}
\psfrag{8}[][][.7]{$8$}
\psfrag{9}[][][.7]{$9$}
\psfrag{G}[][][.7]{$G$}
\psfrag{G0}[][][.7]{$G0$}
\psfrag{G1}[][][.7]{$G1$}
\psfrag{G2}[][][.7]{$G2$}
\psfrag{G02}[][][.7]{$G02$}
\psfrag{G13}[][][.7]{$G13$}
\psfrag{G14}[][][.7]{$G14$}
\psfrag{G15}[][][.7]{$G15$}
\psfrag{G16}[][][.7]{$G16$}
\psfrag{G26}[][][.7]{$G26$}
\psfrag{G021}[][][.7]{$G021$}
\psfrag{G041}[][][.7]{$G041$}
\psfrag{G051}[][][.7]{$G051$}
\psfrag{G026}[][][.7]{$G026$}
\psfrag{G132}[][][.7]{$G132$}
\psfrag{G135}[][][.7]{$G135$}
\psfrag{G136}[][][.7]{$G136$}
\psfrag{G145}[][][.7]{$G145$}
\psfrag{G146}[][][.7]{$G146$}
\psfrag{G0213}[][][.7]{$G0213$}
\psfrag{G0216}[][][.7]{$G0216$}
\psfrag{G0314}[][][.7]{$G0314$}
\psfrag{G0651}[][][.7]{$G0651$}
\psfrag{G0245}[][][.7]{$G0245$}
\psfrag{G5802}[][][.7]{$G5802$}
\psfrag{G1324}[][][.7]{$G1324$}
\psfrag{G1326}[][][.7]{$G1326$}
\psfrag{G1357}[][][.7]{$G1357$}
\psfrag{G1457}[][][.7]{$G1457$}
\psfrag{G1458}[][][.7]{$G1458$}
\psfrag{G1576}[][][.7]{$G1576$}
\psfrag{G02134}[][][.7]{$G02134$}
\psfrag{G02136}[][][.7]{$G02136$}
\psfrag{G02145}[][][.7]{$G02145$}
\psfrag{G03416}[][][.7]{$G03416$}
\psfrag{G03651}[][][.7]{$G03651$}
\psfrag{G06514}[][][.7]{$G06514$}
\psfrag{G14580}[][][.7]{$G14580$}
\psfrag{G24680}[][][.7]{$G24680$}
\psfrag{G13245}[][][.7]{$G13245$}
\psfrag{G13246}[][][.7]{$G13246$}
\psfrag{G13576}[][][.7]{$G13576$}
\psfrag{G14576}[][][.7]{$G14576$}
\psfrag{G132460}[][][.7]{$G132460$}
\psfrag{G132580}[][][.7]{$G132580$}
\psfrag{G021457}[][][.7]{$G021457$}
\psfrag{G245801}[][][.7]{$G245801$}
\psfrag{G021576}[][][.7]{$G021576$}
\psfrag{G132457}[][][.7]{$G132457$}
\psfrag{G132458}[][][.7]{$G132458$}
\psfrag{G132576}[][][.7]{$G13576$}
\psfrag{G0213465}[][][.7]{$G0213465$}
\psfrag{G0215834}[][][.7]{$G0215834$}
\psfrag{G0213657}[][][.7]{$G0213657$}
\psfrag{G0341576}[][][.7]{$G0341576$}
\psfrag{G1324576}[][][.7]{$G1324576$}
\psfrag{G13245760}[][][.7]{$G13245760$}
\psfrag{G02145768}[][][.7]{$G02145768$}
\psfrag{G13245768}[][][.7]{$G13245768$}
\psfrag{G132457680}[][][.7]{$G132457680$}

\begin{center}
\epsfig{file=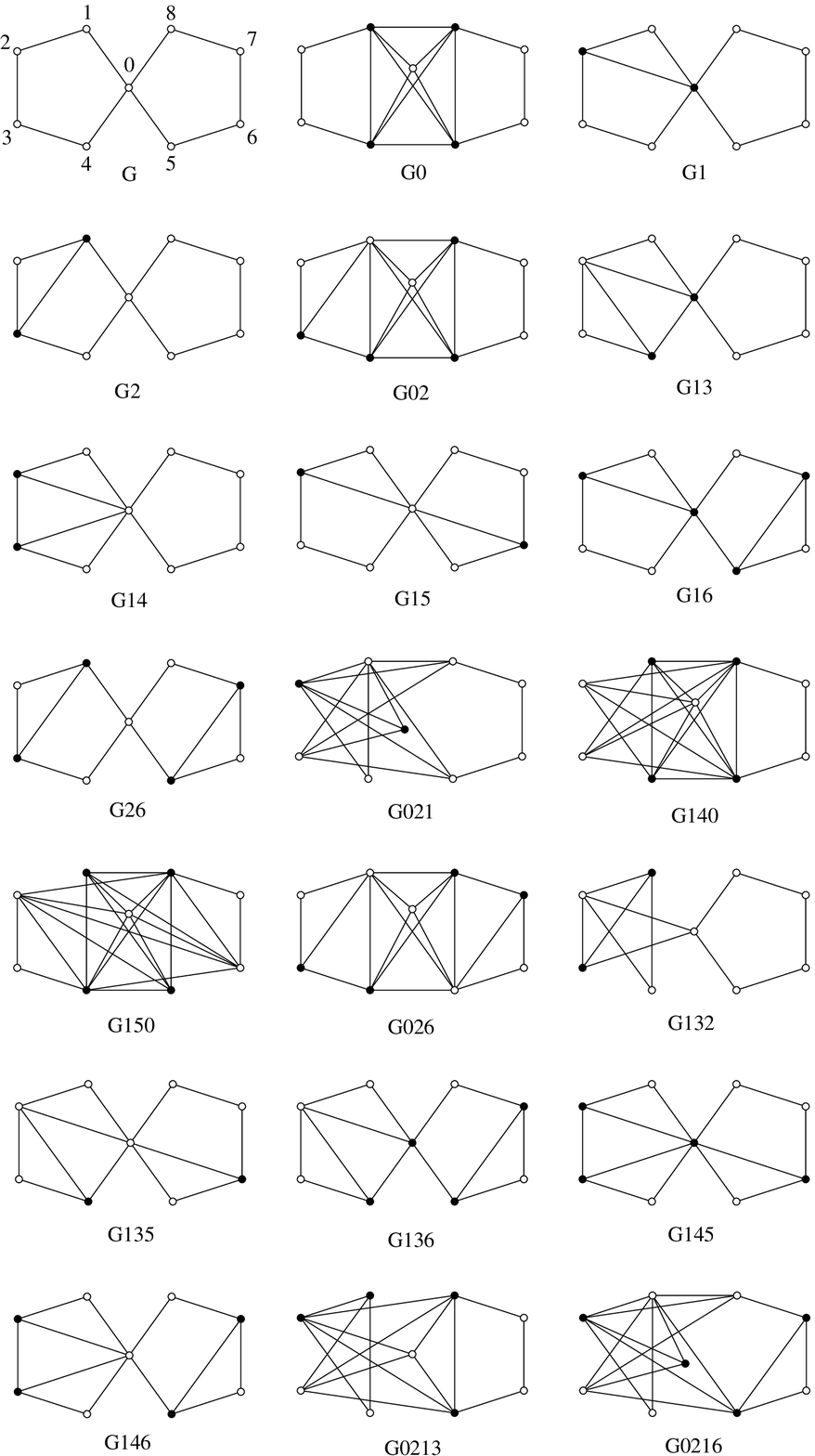, height=23cm}
\end{center}
\begin{center}
\epsfig{file=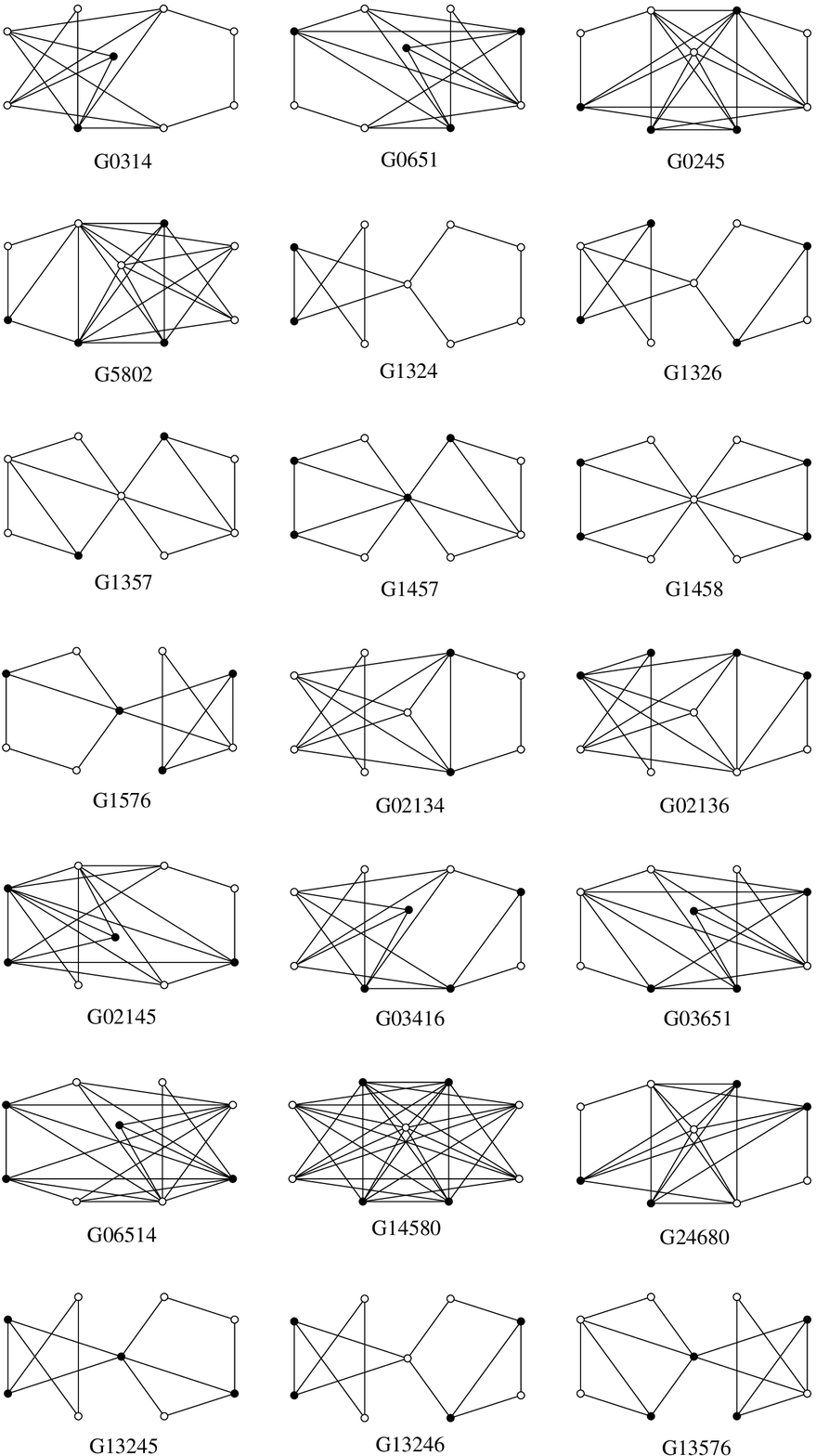, height=23cm}
\end{center}
\begin{center}
\epsfig{file=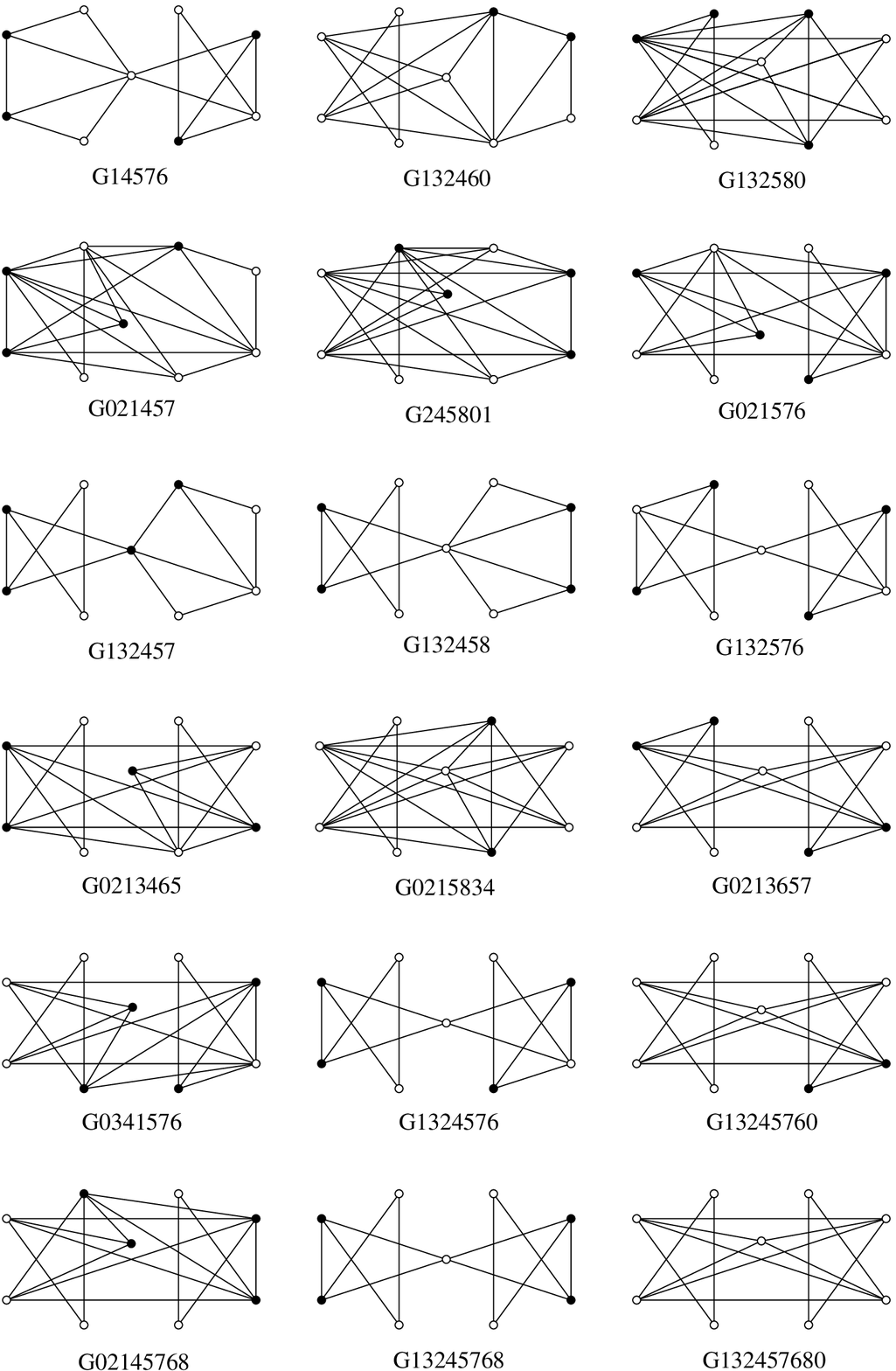, height=19.6cm}
\end{center}
Tout autre graphe de la classe de parité de $G$ est isomorphe à un de
ces 60 représentants non isomorphes. On peut vérifier que la classe
est pure en s'assurant qu'aucun de ces graphes n'a d'anticlique
noire.

%% file: annexe2.tex
\psfrag{0}[][][.7]{$0$}
\psfrag{1}[][][.7]{$1$}
\psfrag{2}[][][.7]{$2$}
\psfrag{3}[][][.7]{$3$}
\psfrag{4}[][][.7]{$4$}
\psfrag{5}[][][.7]{$5$}
\psfrag{6}[][][.7]{$6$}
\psfrag{7}[][][.7]{$7$}
\psfrag{8}[][][.7]{$8$}
\psfrag{9}[][][.7]{$9$}
\psfrag{A}[][][.7]{A}
\psfrag{B}[][][.7]{B}
\psfrag{C}[][][.7]{C}

\psfrag{G}[][][.7]{$G$}
\psfrag{G0}[][][.7]{$G$0}
\psfrag{G02}[][][.7]{$G$02}
\psfrag{G130}[][][.7]{$G$130}
\psfrag{G1C0}[][][.7]{$G$1C0}
\psfrag{G027}[][][.7]{$G$027}
\psfrag{G1C02}[][][.7]{$G$1C02}
\psfrag{G051C}[][][.7]{$G$051C}
\psfrag{G1C06}[][][.7]{$G$1C06}
\psfrag{G2450}[][][.7]{$G$2450}
\psfrag{G4A02}[][][.7]{$G$4A02}
\psfrag{G2960}[][][.7]{$G$2960}
\psfrag{G1C024}[][][.7]{$G$1C024}
\psfrag{G051C2}[][][.7]{$G$051C2}
\psfrag{G02169}[][][.7]{$G$02169}
\psfrag{G134C0}[][][.7]{$G$134C0}
\psfrag{G139C0}[][][.7]{$G$139C0}
\psfrag{G1C069}[][][.7]{$G$1C069}
\psfrag{G4A023}[][][.7]{$G$4A023}
\psfrag{G0274A}[][][.7]{$G$0274A}
\psfrag{G4A028}[][][.7]{$G$4A028}
\psfrag{G4A02C}[][][.7]{$G$4A02C}
\psfrag{G139C02}[][][.7]{$G$139C02}
\psfrag{G1302BC}[][][.7]{$G$1302BC}
\psfrag{G02168C}[][][.7]{$G$02168C}
\psfrag{G134C07}[][][.7]{$G$134C07}
\psfrag{G134C0B}[][][.7]{$G$134C0B}
\psfrag{G139C05}[][][.7]{$G$139C05}
\psfrag{G1C0693}[][][.7]{$G$1C0693}
\psfrag{G139C07}[][][.7]{$G$139C07}
\psfrag{G139C08}[][][.7]{$G$139C08}
\psfrag{G139C0B}[][][.7]{$G$139C0B}
\psfrag{G051C76}[][][.7]{$G$051C76}
\psfrag{G054A62}[][][.7]{$G$054A62}
\psfrag{G054A92}[][][.7]{$G$054A92}
\psfrag{G0BA462}[][][.7]{$G$0BA462}
\psfrag{G0274A8}[][][.7]{$G$0274A8}
\psfrag{G0BA492}[][][.7]{$G$0BA492}
\psfrag{G06A592}[][][.7]{$G$06A592}
\psfrag{G021354697A8BC}[][][.7]{$G$021354697A8BC}

\begin{center}
\epsfig{file=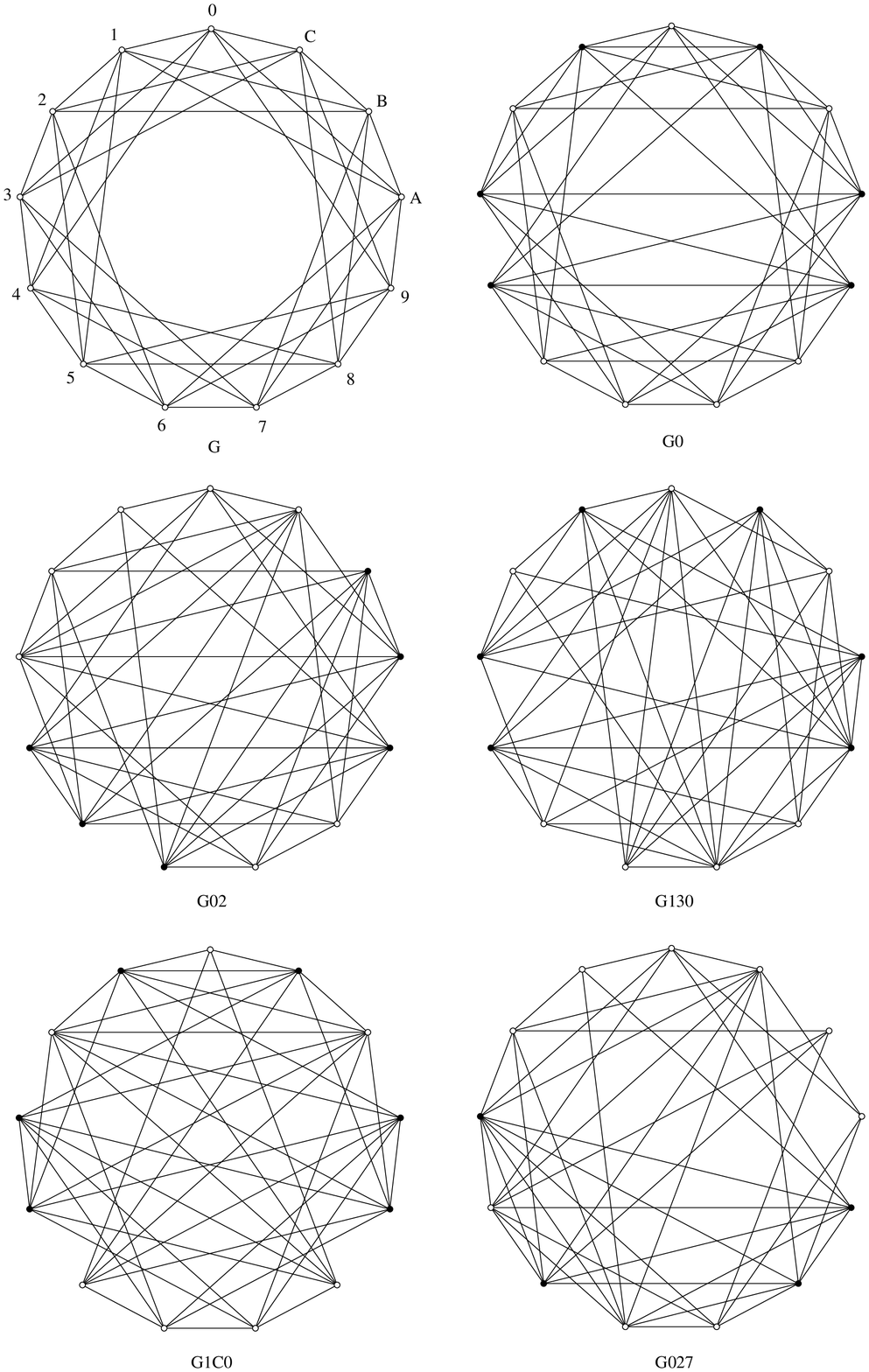, height=22.5cm}
\end{center}

\begin{center}
\epsfig{file=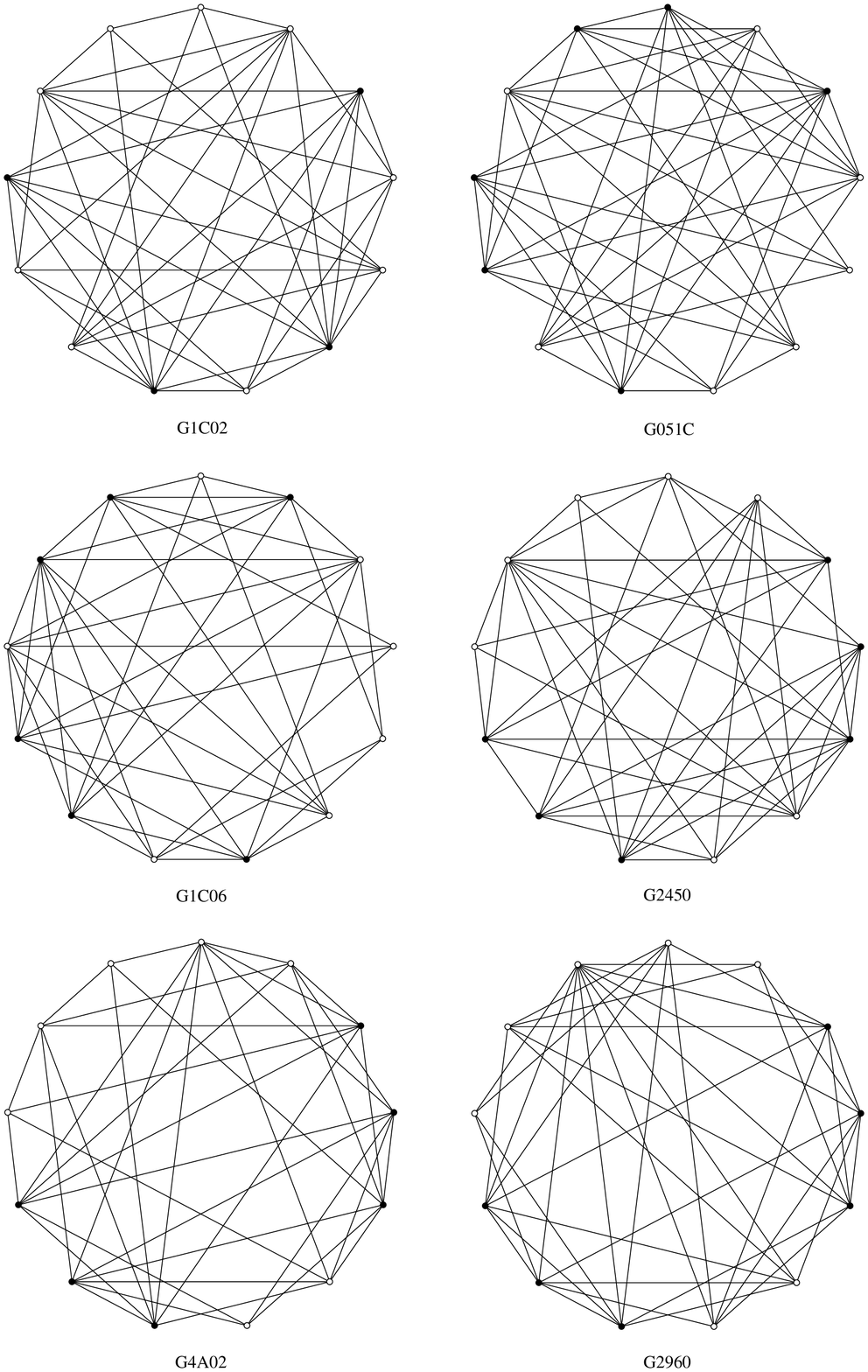, height=22.5cm}
\end{center}

\begin{center}
\epsfig{file=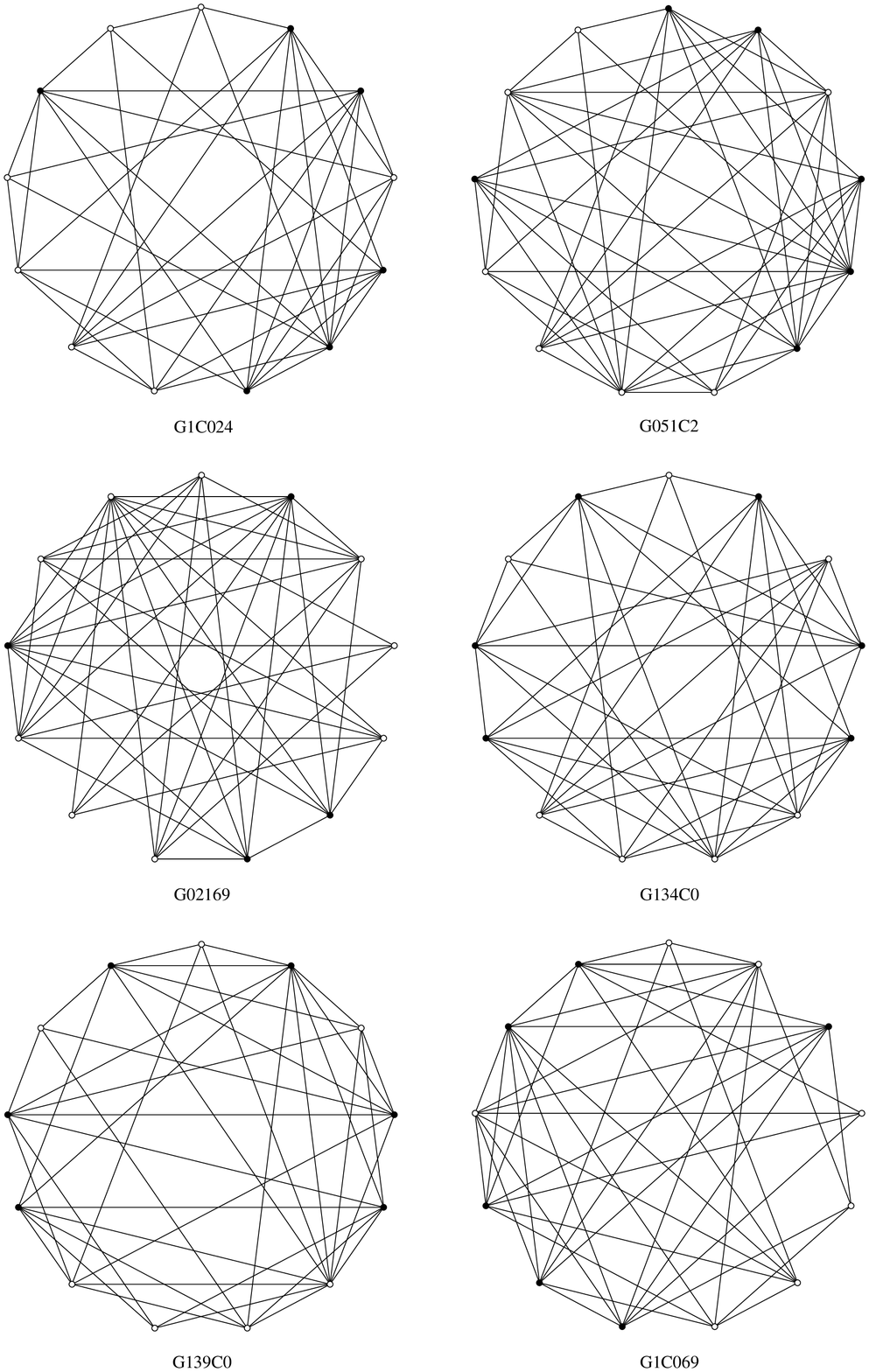, height=22.5cm}
\end{center}

\begin{center}
\epsfig{file=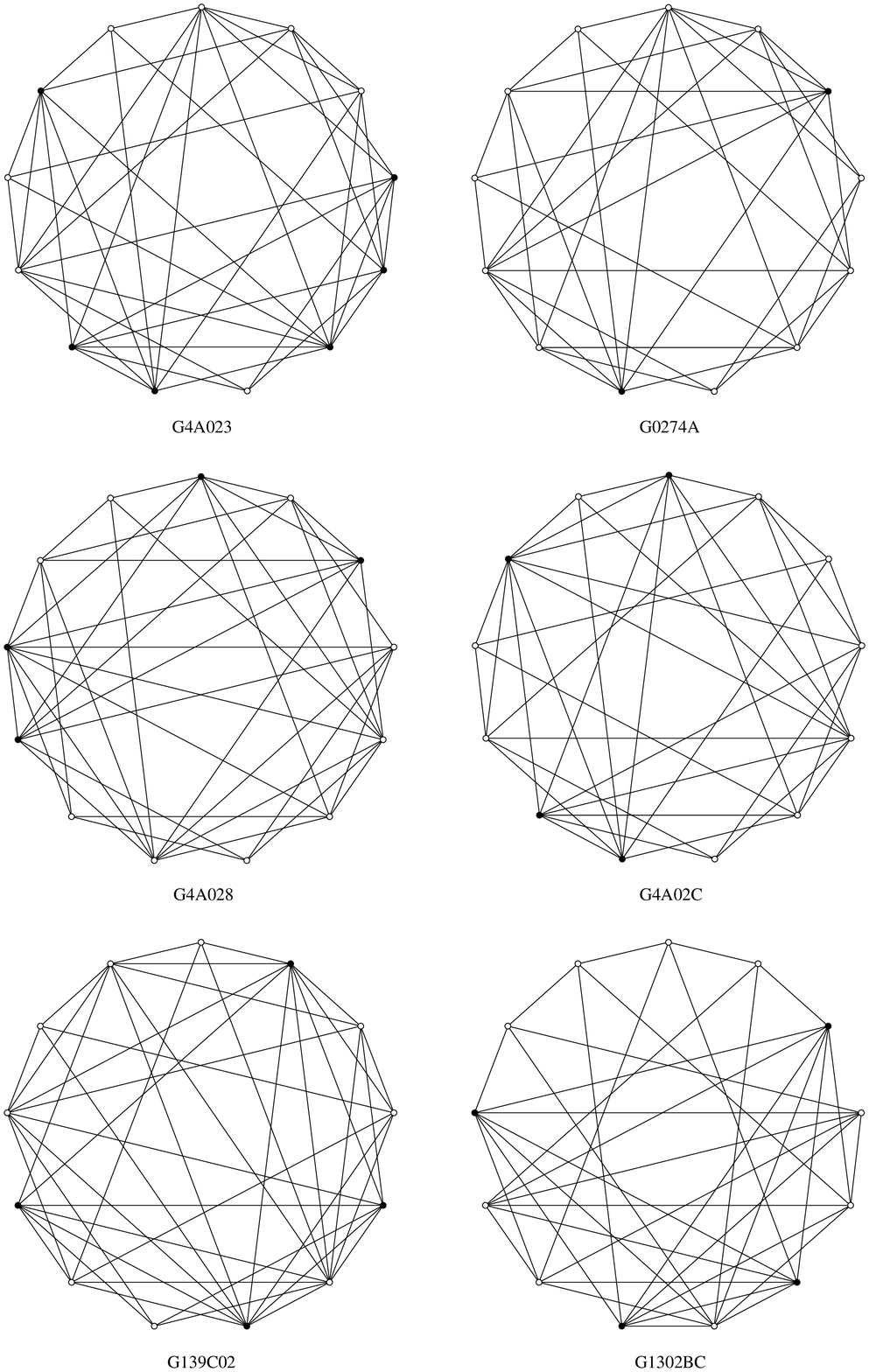, height=22.5cm}
\end{center}

\begin{center}
\epsfig{file=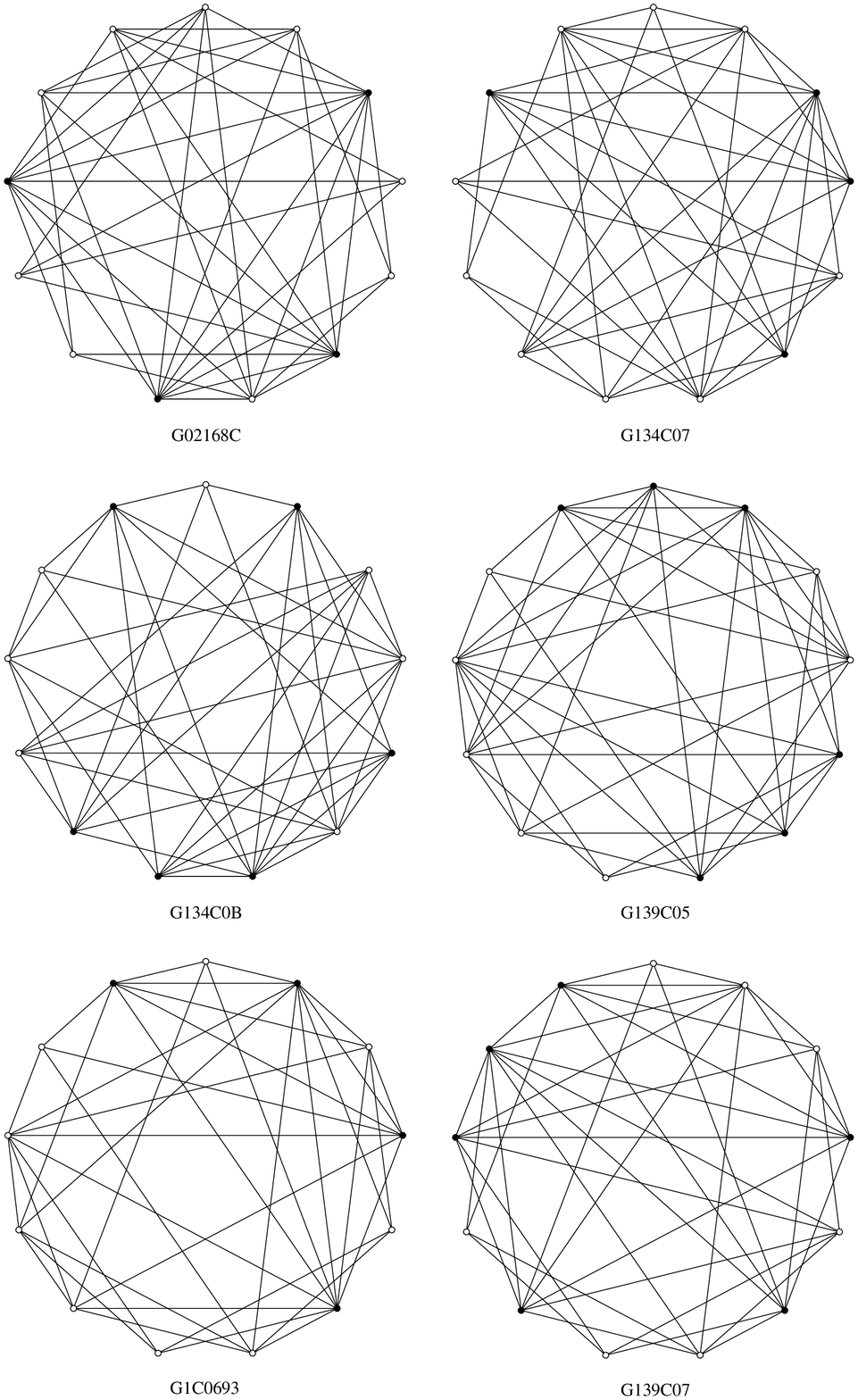, height=22.5cm}
\end{center}

\begin{center}
\epsfig{file=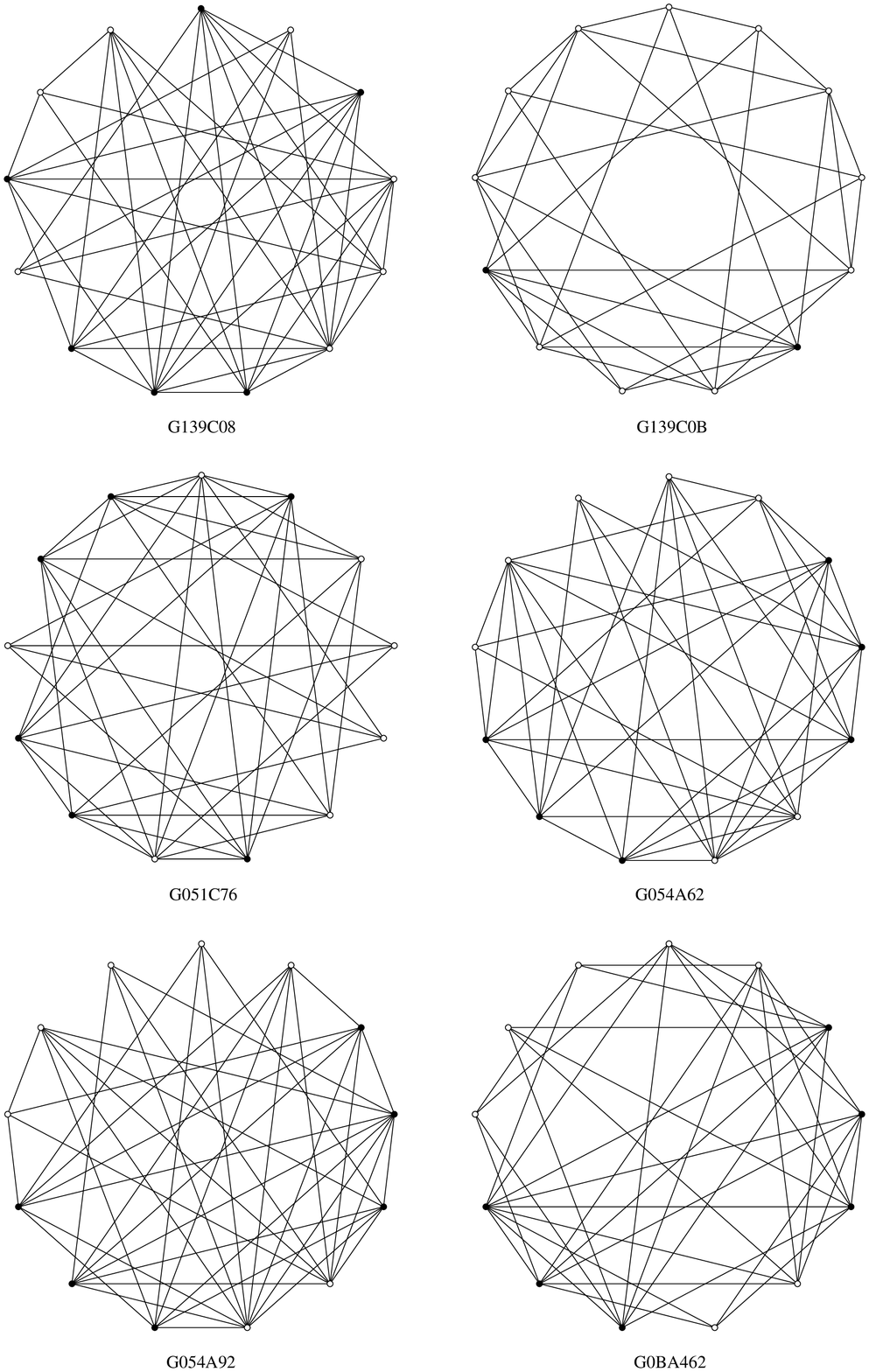, height=22.5cm}
\end{center}

\begin{center}
\epsfig{file=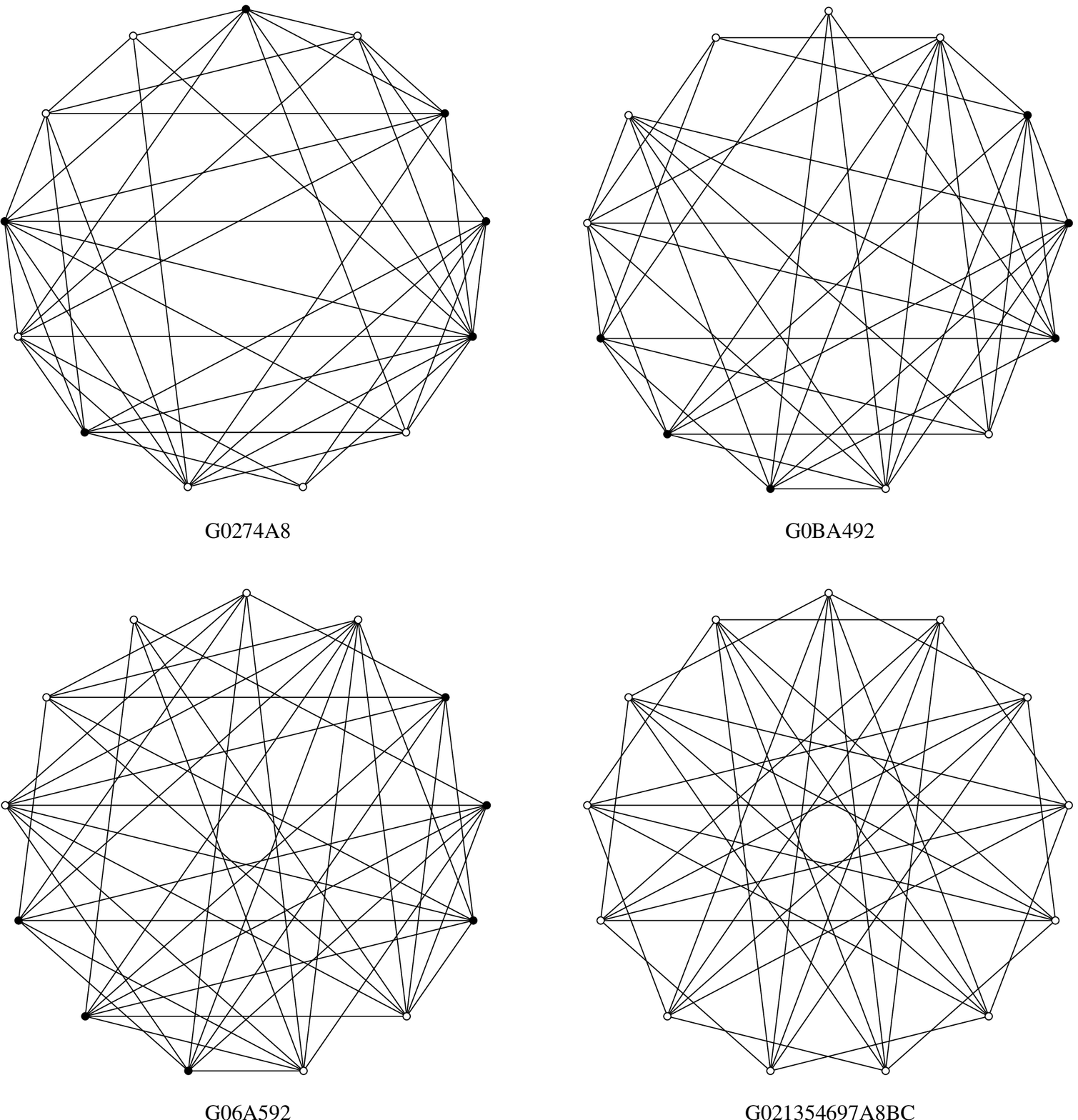, height=14.78cm}
\end{center}

Les trente-neuf premiers graphes sont des représentants non isomorphes de la classe de parité de $G$, le graphe pur sommet transitif naturellement colorié d'ordre treize. Le dernier graphe est l'inverse de $G$.